\documentclass[reqno]{amsart}
\usepackage{setspace}
\usepackage[utf8]{inputenc}
\usepackage{graphicx}
\graphicspath{ {./images/} }
\usepackage[pagewise]{lineno}
\usepackage{amsmath, amsthm, amssymb, amstext}
\usepackage{enumitem}
\usepackage{amssymb}
\usepackage[foot]{amsaddr}
\usepackage[all,arc]{xy}
\usepackage{mathrsfs}
\usepackage[margin=1in]{geometry}
\usepackage{amsmath,amsthm}
\usepackage{mathtools}
\usepackage[T1]{fontenc}
\usepackage{fourier}
\usepackage{cite}
\usepackage{stmaryrd}
\usepackage{amsfonts}
\usepackage{hyperref}
\usepackage[capitalise,nameinlink]{cleveref}
\DeclareMathOperator*{\esssup}{ess\,sup}
\DeclareMathOperator*{\essinf}{ess\,inf}
\usepackage{xcolor}
\hypersetup{
colorlinks,
linkcolor={blue!100!black},
citecolor={red!100!black},
urlcolor={green!100!black}
}

 \newtheorem{theorem}{Theorem}[section]
 
 \newtheorem{lm}[theorem]{Lemma}
\newtheorem{prop}[theorem]{Proposition}
 \theoremstyle{definition}
 \newtheorem{dfn}{Definition}[section]
  \newtheorem{dfns}{Definitions}[section]
\newtheorem{remark}[theorem]{Remark}
\newtheorem{remarks}[theorem]{Remarks}

 \numberwithin{equation}{section}

 \newcommand{\h}{\mathcal{H}}

 \renewcommand{\a }{\alpha }

\newcommand{\vi}{v_i}

\newcommand{\e }{\epsilon }
\newcommand{\g }{\gamma}

 \newcommand{\foral }{\forall\, }

\renewcommand{\l }{\lambda }

\renewcommand{\o }{\omega }
\renewcommand{\O }{\Omega }

\newcommand{\ds}{\displaystyle}   
  \newcommand{\cqfd}

\newcommand{\R}{\mathbb R}
\def\ds{\displaystyle}
\def\w{W^{1,\mathcal{H}}_0(\Omega)} \def\ww{W^{1,\mathcal{H}}(\Omega)} \def\lh{L^{\mathcal{H}}(\Omega)} \def\vi{\rho_{\mathcal{H}}}
\date {}

\DeclarePairedDelimiterX{\inp}[2]{\langle}{\rangle}{#1, #2}

\makeatletter
\def\@tocline#1#2#3#4#5#6#7{\relax
\ifnum #1>\c@tocdepth 
\else
\par \addpenalty\@secpenalty\addvspace{#2}%
\begingroup \hyphenpenalty\@M
\@ifempty{#4}{%
\@tempdima\csname r@tocindent\number#1\endcsname\relax
}{%
\@tempdima#4\relax
}%
\parindent\z@ \leftskip#3\relax \advance\leftskip\@tempdima\relax
\rightskip\@pnumwidth plus4em \parfillskip-\@pnumwidth
#5\leavevmode\hskip-\@tempdima
\ifcase #1
\or\or \hskip 1em \or \hskip 2em \else \hskip 3em \fi%
#6\nobreak\relax
\dotfill\hbox to\@pnumwidth{\@tocpagenum{#7}}\par
\nobreak
\endgroup
\fi}
\makeatother
\usepackage{hyperref}

\numberwithin{theorem}{section}
\numberwithin{equation}{section}

\usepackage{epigraph}

\begin{document}
\title{A new class of anisotropic double phase problems:\\ exponents depending on solutions and their gradients}
\author{Ala Eddine Bahrouni, Anouar Bahrouni and Hlel Missaoui }
\maketitle

\begin{abstract} In this work, we introduce two novel classes of
quasilinear elliptic equations, each driven by the double phase
operator with variable exponents. The first class features a new
double phase equation where exponents depend on the gradient of the
solution. We delve into proving various properties of the
corresponding Musielak-Orlicz Sobolev spaces, including the
$\Delta_2$ property, uniform convexity, density and compact
embedding. Additionally, we explore the characteristics of the new
double phase operator, such as continuity, strict monotonicity, and
the (S$_+$)-property. Employing both variational and nonvariational
methods, we establish the existence  of solutions for this inaugural
class of double phase equations. In the second category, the
treatment of exponents is dependent on the solution itself. This
class differs from the first one due to the unavailability of
suitable Musielak-Orlicz Sobolev spaces. For this reason, we employ
a perturbation argument that leads to the classical double phase
class. These two new classes highlight how different physical
processes like the movement of special fluids through porous
materials, phase changes, and fluid dynamics interact with each
other. Our results are novel in this context and includes a
self-contained techniques.


\end{abstract}

\tableofcontents


\section{Introduction}
In this paper, we present two novel extensions of the double phase
equations with variable exponents. Our first extension adjusts
exponents based on the gradient of the solution, while the second
addresses  exponents depend on the solution itself. This dual
dependence of double phase equations on both solution and solution
gradients emphasizes the complex interactions among different
physical phenomena. It furnishes a flexible framework ideally suited
for modeling an array of complex systems, such as the flow of
non-Newtonian fluids through porous media, phase transitions, and
fluid dynamics...\\

More precisely, our focus is directed towards the following two
equations:
\begin{equation}\label{ext1}
 \left\lbrace
\begin{array}{rll}
-\operatorname{div}\left(|\nabla u|^{p(x,|\nabla u|)-2} \nabla u+ \mu(x)| \nabla u|^{q(x,|\nabla u|)-2}  \nabla u\right) & =f(x,u),  & \text { in } \Omega, \\
\ & \ & \\
u & =0, & \text { on } \partial \Omega,
\end{array}\right.
\end{equation}
and
\begin{equation}\label{ext2}
 \left\lbrace
\begin{array}{rll}
-\operatorname{div}\left(|\nabla u|^{p(x,u)-2} \nabla u+ \mu(x)|\nabla u|^{q(x,u)-2} \nabla u\right)& =f(x,u),  & \text { in } \Omega, \\
\ & \ & \\
u &=0, & \text { on } \partial \Omega,
\end{array}\right.
\end{equation}
where \( \Omega \subset \mathbb{R}^d, d \geq 2 \) is a bounded
domain with a Lipschitz boundary \( \partial \Omega \),
 $f : \Omega \times \mathbb{R}\rightarrow \mathbb{R}$ is a Carath\'eodory function satisfying certain assumptions,  \(\mu \in L^\infty (\Omega,\mathbb{R}_+)\)
   and  $p, \ q :\Omega\times \mathbb{R} \rightarrow [2,+\infty) $ are two
   continuous
   functions.
 \par  When $p(x,t)=p$ and $q(x,t)=q$, then equations \eqref{ext1} and \eqref{ext2}  become the so-called double phase equation
   \begin{equation}\label{doubphas1}
   \left\lbrace
\begin{array}{rll}
-\operatorname{div}\left(|\nabla u|^{p-2} \nabla u+ \mu(x)| \nabla u|^{q-2}  \nabla u\right) &=f(x,u),  & \text { in } \Omega, \\
 \ & \ & \\
u &=0, & \text { on } \partial \Omega.
\end{array} \right.
\end{equation}
The nonlinear and nonhomogeneous differential operator involved in
equation \eqref{doubphas1} is the following weighted double phase
operator
\begin{equation}\label{doubop}
-\operatorname{div}\left(|\nabla u|^{p-2} \nabla u+ \mu(x)| \nabla
u|^{q-2}  \nabla u\right).
\end{equation}
It is easily demonstrable that the integral expression for
\eqref{doubop} is defined by
\begin{equation}\label{doubint}
\int_{\mathbb{R}^{d}}\left(|\nabla u|^{p}+\mu(x)|\nabla
u|^{q}\right)dx.
\end{equation}
This last functional is a basic prototype of a nonautonomous one
featuring nonstandard polynomial growth conditions and a soft kind
of nonuniform ellipticity.  The integral functional \eqref{doubint}
was first introduced by Zhikov \cite{z1,zk1,zk2} to provide models
for strongly anisotropic materials in the framework of
homogenization.  A first mathematical treatment of such problem has
been done in \cite{mg1,mg2,mg3} and can be thought as a model for
composite media with different hardening exponents $p$ and $q$. The
geometric structure resulting from the combination of the two
materials is in fact described by the zero set $\{\mu(x) = 0\}$ of
the coefficient $\mu(\cdot)$, where the transition from $q$-growth
to $p$-growth takes place. Moreover, the functional in
\eqref{doubint} belongs to a class of functionals with non-standard
growth conditions introduced by Marcellini in
\cite{marc1,marc2,mr3}. For more details in the physical backgrounds
and other applications, we refer to \cite{bah1,bah2} (for phenomena
associated with transonic flows) and to \cite{benci} (for models
arising in quantum physics). There is by now a large number of
papers and an increasing interest about double phase equations
\eqref{doubphas1}. With no hope of being complete, let us mention
some pioneering works on existence, nonexistence, and uniqueness of
solutions
\cite{Ambrosio-Essebei-2023,Arora-Fiscella-Mukherjee-Winkert-2023,Liu-Dai-2020,Liu-Winkert-2022,mgg,papa,Ge-Pucci-2022}
and some regularity results \cite{fil1,fil2,fil3}. Recently, Crespo
et al. \cite{CrespoBlanco2022} achieved the variable exponent
version of equation \eqref{doubphas1} for the first time. So far,
there are only few results involving the variable exponent double
phase operator, see \cite{abe,cen,zeng}.
\par When $\mu(x)=0$, equations \eqref{ext1} and \eqref{ext2} revert to the
form of a classical nonlinear equation, the specific structure of
which may vary depending on the solution. Such equations are
commonly encountered in mathematical models representing a wide
range of real-life processes. In \cite{ar}, a system of nonlinear
equations was examined, characterizing the stationary
thermoconvective flow of a non-Newtonian fluid. The models of a
thermistor have been explored in prior works \cite{z1,z2}. In
\cite{rz}, research delved into models of electrorheological fluids,
where the nonlinearity of the governing Navier-Stokes equations
varies based on the applied electromagnetic field. Furthermore,
functionals with growth conditions, contingent upon the solution or
its gradient, have proven effective in denoising digital images. For
instance, models utilizing minimization of functionals with
$p(|\nabla u|)-$growth have been discussed in \cite{bl,bolt,ch},
while \cite{ti} provides insights into a model for image denoising
centered around the minimization of a functional with nonlinearity
dependent on the solution $u$. As far as we know, equations with the
$p(u)-$Laplace operator $\left( \Delta_{p(u)}:=\operatorname{div}
\left( |\nabla u|^{p(u)-2}\nabla u\right) \right)$ in bounded domain
have  been firstly explored in papers \cite{and,Chipot2019}. A
challenging feature of the equations that involve $p(u)-$Laplacian
is that they cannot be interpreted as a duality relation in a fixed
Banach space. Hence, the authors of \cite{and} streamline their
investigation by narrowing it down to the $L^1$-setting. They employ
Young measures to derive a solution, viewing it as the limit of a
sequence ${u_n}$ of solutions originating from the regularized
equations featuring $p_n(x)$ and $p^+$-Laplacian operators. The
authors of \cite{Chipot2019} adopt an alternative approach and
overcome this difficulty adapting the idea of \cite{ch} about
passing to the limit in a sequence $\{|\nabla v_k(x)|^{q_k(x)}\}$.
The case of unbounded domain has been considered in \cite{bah}. In
\cite{sha1,sha2}, the authors investigated the evolutive equations
associated with equations \eqref{ext1} and \eqref{ext2}. It is worth
noting that, when $\mu(x)=0$ the equation with a variable exponent,
dependent on the gradient of the solution, remains relatively
underexplored. As far as our current knowledge extends, there is a
dearth of research addressing these issues
through the application of variational methods.\\

This study aims to explore the intersection of the phenomena
described above in the context of two specific problems defined by
equations \eqref{ext1} and \eqref{ext2}. To our knowledge, this
represents the inaugural investigation into equations \eqref{ext1}
and \eqref{ext2} in the literature. This work aims to fill this gap
by studying these important problems. Despite the novelty presented
by the last two problems, we can discern genuine motivation and
practical applications. When the exponents in double phase equations
depend on the gradient of the solution, it suggests a coupling
between the system's dynamics and the spatial variations in the
solution. This coupling can lead to richer and more realistic models
that capture the complex interactions between different components
or phases. The double phase equation with exponents that depend on
the solution itself arises from various physical systems where the
governing equations exhibit nonlinearity and self-interaction. More
precisely, many physical systems exhibit nonlinear behavior, meaning
their dynamics cannot be described simply through linear
relationships between variables. Nonlinear systems often involve
feedback loops or interactions between different components, leading
to complex behaviors. The double phase equation with
solution-dependent exponents can capture such nonlinear dynamics,
allowing for a more realistic representation of the system's
behavior. Compared to the classical PDEs, the problems \eqref{ext1}
and \eqref{ext2} are notably more demanding and present a distinct
set of obstacles that require careful consideration and ingenuity.
More precisely, the presence of $\nabla u$ in the exponent of
equation \eqref{ext1} adds intricacy to selecting the suitable
generalized N-function and establishing the corresponding
Musielak-Orlicz Sobolev spaces. To address these issues, we develop
specific tools tailored to this situation. In particular, we focus
on proving different properties of the Musielak-Orlicz Sobolev
spaces involved, see Subsection \ref{3.1}. These characteristics include  the $\Delta_2$
property (see Propositions \ref{dlt2} and \ref{delta2}), reflexivity (see Proposition \ref{rflxv}),   compact embedding (see Theorem \ref{inject10} and Proposition \ref{prp3}), and density (see Theorem \ref{dnsty}). Note
that the inclusion of the function $\mu(x)$ in equations
\eqref{ext1} and \eqref{ext2} hinders the application of the
standard approach for addressing the density and embedding
properties. Additionally, we look into the characteristics of the
new double phase operator (see Subsection \ref{3.2}), such as continuity, strict monotonicity,
and the (S$_+$)-property (see Proposition \ref{op1} and Theorem \ref{op2}). By combining these findings with
variational and nonvariational arguments, we prove the existence of
nontrivial solutions for equation \eqref{ext1} (see Theorems
\ref{thmA1} and \ref{thmA2} ). Unlike the issue presented in problem
\eqref{ext1}, the situation depicted in the second equation
\eqref{ext2} contrasts significantly, see Section \ref{last}. This occurs as a result of the
solution being present in the exponent, which leads to the
unavailability of suitable Musielak-Orlicz Sobolev spaces. For this
reason, we use a perturbation argument (see \cite{Chipot2019}) that
leads to the classical double phase equation with variable exponents
introduced in \cite{CrespoBlanco2022}. In particular, we establish
the existence and multiplicity of solutions for equations
\eqref{ext2} (see Theorems \ref{th1} and \ref{th2}). Note that the
multiplicity findings presented here are novel concerning equation
\eqref{ext2},
specifically when $\mu(x)=0$.\\

The outline of the paper is the following: in Section $2$, we
collect some preliminary results.  Section 3 is dedicated to proving
abstract results as discussed earlier. We then employ both
variational and nonvariational methods to establish the existence of
solutions to equation \eqref{ext1}. Section 4 focuses on exploring
the existence of multiple solutions to equation \eqref{ext2}.
\section{Preliminaries}
This section is divided into two subsections. In the first part, we
recall some known definitions and results about the
Musielak-Orlicz-Sobolev space. In the second part, we deal with the
double phase space with variable exponent, which is a particular
case of  Musielak-Orlicz-Sobolev  space introduced  in
\cite{CrespoBlanco2022}. Moreover, we give some new topological
properties related to these topics that will be useful in treating
our main problems.

\subsection{Museilak-Orlicz-Sobolev spaces}
In this section, we recall some definitions and properties related
to the Musielak-Orlicz and the  Musielak-Orlicz-Sobolev spaces, see
\cite{Chlebicka2021, Diening2011, Fan2012a, Fan2012b, Musielak1983}.

\begin{dfn} Let $\Omega$ be an open subset of $\mathbb{R}^d$. A function $\h:\Omega\times\mathbb{R}\longrightarrow \mathbb{R}$ is called a generalized N-function if it satisfies the following conditions:
    \begin{enumerate}
        \item[$(1)$] For a.a. $x\in \O$, $\h(x,t)$ is even, continuous, nondecreasing and convex in $t$, and for each $t\in \mathbb{R}, \ \ \h(x,t) $ is measurable in $x$;
        \item[$(2)$]$\displaystyle{\lim\limits_{t\rightarrow 0}\frac{\h(x,t)}{t}=0}$, for a.a. $x\in \Omega;$
        \item[$(3)$]$\displaystyle{\lim\limits_{t\rightarrow \infty}\frac{\h(x,t)}{t}=\infty}$, for a.a. $x\in \Omega;$
        \item[$(4)$] $\displaystyle{\h(x,t)>0}$, for all $t>0$ and all $x\in \Omega$ and $\h(x,0)=0$, for all $x\in \O$.
    \end{enumerate}
\end{dfn}
\vspace{3mm}
\begin{remark}
 We give an equivalent definition of a generalized N-function that admits an integral representation.
For $x \in \O$ and $t \geq 0$, we denote by $h(x,t)$ the right-hand
derivative of $\h(x,\cdot)$ at $t$, and define $h(x,t) =-h(x,-t)$
for $t < 0$. Then for each $x \in \O$, the function $h(x,\cdot)$ is
odd, $h(x,t) \in \R$ for $t \in \R,\ h(x, 0) = 0,\  h(x,t) > 0$ for
$t > 0,\ h(x,\cdot) $ is right-continuous and nondecreasing on
$[0,+\infty),\  h(x,t)\longrightarrow +\infty$ as $t \rightarrow
+\infty ,$ and
$$
\h(x,t)=\int_{0}^{|t|} h(x,s)ds, \text{ for } x\in \O \text{ and } t
\in \R.
$$
\end{remark}
\begin{dfn}
We say that a generalized N-function $\h$ satisfies the
$\Delta_2$-condition if there exist $C_0 > 0$ and a nonnegative
function $m \in L^1(\Omega)$ such that
$$\h(x,2t)\leq C_0\h(x,t)+m(x),\ \ \text{for a.a.}\ x\in \Omega\ \text{and all}\ t\geq 0.$$
\end{dfn}
\begin{dfns}
 Let $\h_1$ and $\h_2$ be two generalized N-functions.
 \begin{itemize}
  \item [(1)]  We say that $\h_1$ increases essentially slower than $\h_2$ near infinity and we write $\h_1 \ll \h_2$ , if for any $k>0$
$$
\lim _{t \rightarrow \infty} \frac{\h_1(x, k t)}{\h_2(x, t)}=0,\
\text {uniformly in}\  x \in \Omega.
$$
  \item [(2)]We say that $\h_1 $ is weaker than $\h_2$, denoted by $\h_1 \preceq \h_2$, if there exist two positive constants $C_1, C_2$ and a nonnegative function $m \in L^1(\Omega)$ such that
$$
\h_1(x, t) \leq C_1 \h_2\left(x, C_2 t\right)+m(x),\ \ \text{for
a.a.}\ x\in \Omega\ \text{and all}\ t\geq 0.
$$
 \end{itemize}
\end{dfns}

\begin{dfn}\label{CF}
For any generalized N-function $\h$, the function
$\widetilde{\h}:\Omega\times\mathbb{R}\longrightarrow \mathbb{R}$
defined by
\begin{equation}\label{Cf}
\widetilde{\h}(x,t):=\sup_{\tau\geq 0}\left(
t\tau-\h(x,\tau)\right),\ \ \text{for all}\ x\in \Omega\ \text{and
all}\ t\geq 0,
\end{equation}
 is called the complementary function of $\h$.
\end{dfn}
In view of the definition of the complementary function
$\widetilde{\h}$, we have the following Young's type inequality:
\begin{equation}\label{Yi}
    \tau\sigma\leq \h(x,\tau)+\widetilde{\h}(x,\sigma),\ \text{for all}\ x\in \Omega\ \text{and all}\ \tau,\sigma\geq0.
\end{equation}

\begin{remarks}\label{Delta}\begin{itemize}
    \item [(1)] Note that the complementary function $\widetilde{\h}$ is also a generalized N-function.
    \item [(2)] If there exist $m,\ell\in \mathbb{R}$ such that
    \begin{equation}\label{D2}
        1\leq m\leq \frac{h(x,t)t}{\h(x,t)}\leq \ell,\ \ \text{for all}\ x\in \Omega\ \text{and all}\ t> 0,
    \end{equation}
     then, the generalized N-function $\h$ satisfies the $\Delta_2$-condition.
   \end{itemize}
\end{remarks}
\begin{lm}\label{lm1}
     Let $\h$ be a generalized $N$-function. We suppose that $t\mapsto h(x,t)$ is continuous and increasing on $\R$, for a.a. $x\in \O$. Moreover, we assume that there exist $m,\ell\in \mathbb{R}$ such that
    \begin{equation}\label{D22}
        1<m\leq \frac{h(x,t)t}{\h(x,t)}\leq \ell,\ \ \text{for all}\ x\in \Omega\ \text{and all}\ t> 0,
    \end{equation}
    then,
    \begin{equation}\label{L1}
       \widetilde{\h}(x,h(x,s))\leq (\ell-1)\h(x,s),\ \ \text{for all}\ s\geq 0\ \text{and all}\ x\in \Omega,
    \end{equation}
    and
   \begin{equation}\label{D3}
        \frac{\ell}{\ell-1}:=\widetilde{m}\leq \frac{\widetilde{h}(x,s)s}{\widetilde{\h}(x,s)}\leq \widetilde{\ell}:=\frac{m}{m-1},\ \ \text{for all}\ x\in \Omega\ \text{and all}\ s> 0,
    \end{equation}
    where $\displaystyle{\widetilde{\h}(x,s)=\int_{0}^s \widetilde{h}(x,\eta)\ d\eta}.$
\end{lm}
\begin{proof}
For $s\geq 0$, let $f(t)=h(x,s)t-\h (x,t)$, for all $x\in \O$ and
all $t\geq0.$ Then, it is immediate to see that $f(t)\leq f(s)$, for
all $t\geq 0,$
 which implies that
\begin{equation*}
 h(x,s)t -\h(x,t) \leq h(x,s)s-\h(x,s), \text{ for all } x\in \O \text{ and all } t \geq 0.
\end{equation*}
Hence, for all $x\in \Omega$ and all $s\geq 0$, we have
\begin{equation}\label{L2}
    \widetilde{\h}(x,h(x,s)) =\sup_{t\geq0} \big( h(x,s)t-\h(x,t)\big)=\max_{t\geq0} \big( h(x,s)t-\h(x,t)\big)=h(x,s)s-\h(x,s).
\end{equation}
It follows, from \eqref{D22} and \eqref{L2}, for all $x\in \Omega$
and all $s\geq 0$, that
\begin{align*}
    \widetilde{\h}(x,h(x,s)) &\leq  (\ell-1)\h(x,s).
\end{align*}
Thus, \eqref{L1} holds. It remains to show \eqref{D3}. According to
Definition \ref{CF}, we see that $\widetilde{\widetilde{\h}}=\h$, so
replacing in \eqref{L2} $\h$ with $\widetilde{\h}$, for all $x\in
\Omega$ and all $s>0$, we get
\begin{equation}\label{L3}
    \h(x,\widetilde{h}(x,s))=\widetilde{\widetilde{\h}}(x,\widetilde{h}(x,s))=\widetilde{h}(x,s)s-\widetilde{\h}(x,s).
\end{equation}
Since $h(x,.)$ is continuous and increasing for a.a. $x\in  \Omega$, 
we infer that
\begin{equation}\label{L4}
    h(x,\widetilde{h}(x,s))=s,\ \text{for all}\ x\in \Omega\ \text{and all}\ s>0.
\end{equation}
Substituting $t$ by $\widetilde{h}(x,s)$ in \eqref{D22} and using
\eqref{L4}, we obtain
\begin{equation*}
    m\h(x,\widetilde{h}(x,s))\leq \widetilde{h}(x,s)s\leq  \ell\h(x,\widetilde{h}(x,s)),\ \text{for all}\ x\in \Omega\ \text{and all}\ s>0.
\end{equation*}
Therefore, from \eqref{L3}, it yields
\begin{equation*}
        \frac{\ell}{\ell-1}:=\widetilde{m}\leq \frac{\widetilde{h}(x,s)s}{\widetilde{\h}(x,s)}\leq \widetilde{\ell}:=\frac{m}{m-1},\ \ \text{for all}\ x\in \Omega\ \text{and all}\ s> 0.
    \end{equation*}
    Thus, the proof is completed.
\end{proof}
\begin{remark}\label{compl}
   We see that the condition \eqref{D22} implies that $\h$ and its complementary function $\widetilde{\h}$ satisfy the $\Delta_2$-condition.
\end{remark}
In the following lemma, we extend the well-known Simon identity (see
\cite{M50}) to the framework of Musielak-Orlicz space.
\begin{lm}\label{simon} Assume that assumption \eqref{D22} holds. Then, we have
$$\left(h(x,\tau)-h(x,\sigma)\right)(\tau -\sigma)\geq 4 \h\left(x,\frac{\vert\tau-\sigma\vert}{2}\right),\ \text{for all}\ \tau,\sigma \in \mathbb{R}\setminus\lbrace 0\rbrace \ \text{and all }\ x\in \Omega.$$
\end{lm}
\begin{proof}
    The proof is similar  to the one in \cite{M50}.
\end{proof}
 Now, we  define the Musielak-Orlicz space as follows:
$$L^{\h}(\Omega):=\left\lbrace u:\Omega\longrightarrow \mathbb{R}\ \text{measurable :}\ \rho_{\h}(\lambda u)<+\infty,\ \ \text{for some}\ \lambda>0\right\rbrace,$$
where
\begin{equation}\label{Mo}
  \rho_{\h}(u):= \int_{\Omega}\h(x, u)dx.
\end{equation}
The space $L^{\h}(\Omega)$ is endowed with the Luxemburg norm
\begin{equation}\label{No}
  \Vert u\Vert_{\lh}:=\inf\left\lbrace \lambda:\ \rho_{\h}\left(x,\frac{u}{\lambda}\right)\leq 1\right\rbrace.
\end{equation}
\begin{prop}\label{HM}
    Let $\h$ be a generalized $N$-function satisfies the $\Delta_2$-condition, then
    $$\lh =\left\lbrace  u: \O \longrightarrow \mathbb{R}\ \text{measurable :}\ \rho_{\h}( u)<+\infty\right\rbrace.$$
\end{prop}
\begin{prop} \label{zoo}
     Let $\h$ be a generalized $N$-function satisfies \eqref{D22}, then the following assertions hold:
     \begin{itemize}
    \item [(1)] $\min \{\l^m, \l^\ell\}\h(x,t)\leq  \h(x,\l t)\leq \max \{ \l^m, \l^\ell\}\h(x,t),\text{ for a.a. } x \in \O$\\ $\text{ and all }  \l, \ t \geq 0. $
 \item [(2)]  $\min \left\{\|u\|_{\lh}^{m},\|u\|_{\lh}^{\ell}\right\} \leq \rho_{\mathcal{H}}(u) \leq\max \left\{\|u\|_{\lh}^{m},\|u\|_{\lh}^{\ell}\right\},$ \ for all $  u\in \lh$.
 \item [(3)] Let $\left\{ u_n\right\}_{n \in \mathbb{N}}\subseteq \lh $ and $  u \in \lh $, then
 $$
 \|u_n-u\|_{\lh }\longrightarrow0 \ \Longleftrightarrow  \rho_\h (u_n-u) \longrightarrow 0,\
\text{ as } n \rightarrow +\infty. $$
\end{itemize}
\end{prop}
As a consequence of \eqref{Yi}, we have the following result:
\begin{lm}[H\"older's type
inequality]\label{H1}
  Let $\Omega$ be an open subset of $\mathbb{R}^d$ and $\h$ be a generalized $N$-function satisfies \eqref{D22}, then
  \begin{equation}\label{Ho}
     \left\vert \int_{\Omega} uv dx \right\vert \leq 2 \Vert u\Vert_{\lh}\Vert v\Vert_{L^{\widetilde{\h}}(\Omega)},\ \text{for all}\ u\in L^{\h}(\Omega)\ \text{and all}\ v\in L^{\widetilde{\h}}(\Omega).
  \end{equation}
\end{lm}
The subsequent proposition deals with some topological properties of
the Musielak-Orlicz space, see \cite[Theorem 7.7 and Theorem
8.5]{Musielak1983}.
\begin{prop}\label{AB}
\begin{itemize}
\item [(1)] Let $\h$ be a generalized N-function and $\Omega$ an open subset of
 $\mathbb{R}^d$. Then,
 \begin{enumerate}
     \item[(i)] the space
    $\left(L^{\h}(\Omega),\|\cdot\|_{\lh}\right)$ is a Banach space,
       \item[(ii)] if $\h$ satisfies \eqref{D2}, then
$L^{\h}(\Omega)$ is a separable space,
     \item[(iii)] if $\h$ satisfies \eqref{D22}, then
$L^{\h}(\Omega)$ is a reflexive space.
 \end{enumerate}
\item [(2)] Let $\h_1$ and $\h_2$  be two generalized N-functions such that $\h_1 \preceq \h_2$. Then,
$$
L^{\h_2}(\Omega) \hookrightarrow L^{\h_1}(\Omega).
$$\end{itemize}
\end{prop}

Now, we are ready to define the Musielak-Orlicz Sobolev space. Let
$\h$ be a generalized N-function and $\Omega$ an open subset of
$\mathbb{R}^d$. The Musielak-Sobolev space is defined as follows
$$W^{1,\h}(\Omega):=\left\lbrace u\in L^{\h}(\Omega):\ |\nabla u| \in L^{\h}(\Omega)\right\rbrace.$$
The space $W^{1,\h}(\Omega)$ is endowed with the norm
\begin{equation}\label{NM}
  \Vert u\Vert_{\ww}:=\Vert u\Vert_{\lh}+\Vert \nabla u\Vert_{\lh},\ \ \text{for all}\ u\in W^{1,\h}(\Omega),
\end{equation}
where $\|\nabla u\|_{\lh} := \| |\nabla u| \|_{\lh}$.\\ We denote by
$\w$
 the completion of $C^\infty _0(\O)$ in $\ww$.
\begin{theorem}\label{thc}\cite[Theorem 10.2] {Musielak1983}\cite[Proposition 1.7]{Fan2012b}~~~~~~~~~~~~~~~~~~~~~~~~~~~~~~~~~~~~~~~~~~~~~~~~~~~~~~~~~~~~~~

 Let $\h$ be a generalized N-function and locally integrable such that
\begin{equation}\label{c1c}
\inf _{x \in \Omega} \h(x, 1)>0.
\end{equation}
Then, the spaces $W^{1, \h}(\Omega)$ and $W_0^{1, \h}(\Omega)$ are
separable Banach spaces which are reflexive if $L^{\h}(\Omega)$ is
reflexive.
\end{theorem}

\begin{remark}\label{Ref}
 If $\h$  satisfies \eqref{D22}, then the space
$W^{1,\h}(\Omega)$ is a reflexive and separable Banach space with
respect to the norm $\Vert \cdot\Vert_{\ww}$.
\end{remark}


 \subsection{Double phase space with variable exponent}
In this section, we focus our attention on a special case of the
generalized N-function, $\displaystyle{\h(x,t):= t^{p(x)} +
\mu(x)t^{q(x)}}$,
which possesses an additional topological properties. This type of generalized N-function is employed in the study of the double phase equations with variable exponents. For more details about this subject, we refer readers to \cite{CrespoBlanco2022,Ho2023}.\\

 Let $\O \subset \R^d$ $(d\geq  2)$ be a bounded domain with Lipschitz boundary $ \partial \Omega$ and  we consider a function $\h : \O \times [0,+\infty)\rightarrow [0,+\infty)$ defined by $$\h(x,t):= t^{p(x)} + \mu(x)t^{q(x)},\ \text{for all}\ x \in \O \ \text{and all } t\geq0. $$
Here, the exponents $p, \ q $  satisfy the following assumption
\begin{itemize}
  \item[$(\mathrm{H}_0)$] $:\ p,q \in C(\overline{\O}), \ 1 <p(x)<d,\ p(x)< q(x), \text{ for a.a. } x \in \O \text{ and } 0 \leq \mu(\cdot) \in  L^1(\O)$.
\end{itemize}
It is clear that $\h$ is a generalized N-function which satisfies
\eqref{c1c}, locally integrable and it fulfills the
$\Delta_2$-condition, precisely
\begin{equation}\label{D21}
1<p^-\leq \frac{h(x,t)t}{\h(x,t)}\leq q^+,\ \ \text{for all}\ x\in
\Omega\ \text{and all}\ t> 0,
\end{equation}
where $p^-:=\ds\inf_{x\in\O} p(x)$ and $\ q^+:=\ds\sup_{x\in\O}
q(x)$. Therefore, the results introduced in the last subsection
remain valid.

The next proposition gives some embedding results, see
\cite[Proposition 2.16]{CrespoBlanco2022}.

 \begin{prop} Let hypotheses $(\mathrm{H}_0)$ be satisfied. Then, the following Sobolev embedding hold:
    \begin{enumerate}
    \item  $L^{\mathcal{H}}(\Omega) \hookrightarrow L^{r(\cdot)}(\Omega)$, $W^{1, \mathcal{H}}(\Omega) \hookrightarrow W^{1, r(\cdot)}(\Omega), W_0^{1, \mathcal{H}}(\Omega) \hookrightarrow W_0^{1, r(\cdot)}(\Omega)$ are continuous, for all $r(\cdot) \in C(\overline{\Omega})$ with $1 \leq r(x) \leq p(x)$, for all $x \in \overline{\Omega}$;
    \item Let $p \in C^{0,1}(\overline{\Omega})$ and $ p^*(x):=\dfrac{d p(x)}{d-p(x)}$, for all $x \in \overline{\Omega}$. Then,
      \begin{itemize}
\item[(i)]  $W^{1, \mathcal{H}}(\Omega) \hookrightarrow L^{r(\cdot)}(\Omega)$ and $W_0^{1, \mathcal{H}}(\Omega) \hookrightarrow L^{r(\cdot)}(\Omega)$ are continuous for $r(\cdot) \in C(\overline{\Omega})$ with $1 \leq r(x) \leq  p^*(x)$, for all $x \in \overline{\Omega}$;
\item[(ii)] $W^{1, \mathcal{H}}(\Omega) \hookrightarrow L^{r(\cdot)}(\Omega)$ and $W_0^{1, \mathcal{H}}(\Omega) \hookrightarrow L^{r(\cdot)}(\Omega)$ are compact for $r(\cdot) \in C(\overline{\Omega})$ with $1 \leq r(x)<p^*(x)$, for all $x \in \overline{\Omega}$;
\end{itemize}
\item Let $p \in W^{1, \g}(\Omega)$ for some $\g>d$ and $p_*(x):=\dfrac{(d-1)p(x)}{d-p(x)}$, for all $x \in \overline{\Omega}$, then,
\begin{itemize}
\item[(i)] $W^{1, \mathcal{H}}(\Omega) \hookrightarrow L^{r(\cdot)}(\partial \O)$ and $W_0^{1, \mathcal{H}}(\Omega) \hookrightarrow$ $L^{r(\cdot)}(\partial \O)$ are continuous for $r(\cdot) \in C(\overline{\Omega})$ with $1 \leq r(x) \leq p_*(x)$, for all $x \in \overline{\Omega}$;
\item[(ii)] $W^{1, \mathcal{H}}(\Omega) \hookrightarrow L^{r(\cdot)}(\partial \O)$ and $W_0^{1, \mathcal{H}}(\Omega) \hookrightarrow L^{r(\cdot)}(\partial \O)$ are compact for $r(\cdot) \in C(\overline{\Omega})$ with $1 \leq r(x)<p_*(x)$, for all $x \in \overline{\Omega}$;
\end{itemize}
\item  If $\mu \in L^{\infty}(\Omega)$, then $L^{q(\cdot)}(\Omega) \hookrightarrow L^{\mathcal{H}}(\Omega)$ is continuous.
\end{enumerate}
 \end{prop}
 \vspace{3mm}
 We present the following embedding results, where the exponents are constants, which are crucial in Section 4.
 \begin{prop}\cite[Proposition 2.15]{Cola2018} \label{prp}
We assume that $p(x)\equiv p$ and $q(x)\equiv q,$ where $p,\ q$ are
two constants such that $1<p<q<\infty.$ Then the following
embeddings hold:
\begin{itemize}

  \item [(i)] If $p \neq d$, then
$ W_0^{1, \mathcal{H}}(\Omega) \hookrightarrow L^r(\Omega)  \text {
is continuous for all } r \in\left[1, p^*\right].$

\item [(ii)] If $p=d$, then
$ W_0^{1, \mathcal{H}}(\Omega) \hookrightarrow L^r(\Omega) \text {
is continuous for all } r \in[1,+\infty]. $
  \item [(iii)] If $p < d$, then
$ W_0^{1, \mathcal{H}}(\Omega)  \hookrightarrow L^r(\Omega)  \text {
is compact for all } r \in\left[1, p^*\right). $
 \item [(iv)]If $p> d$, then
$ W_0^{1, \mathcal{H}}(\Omega)  \hookrightarrow L^{\infty}(\Omega)
\text{ is  compact}. $

\end{itemize}

\end{prop}

Next, we consider the following assumptions:
 \begin{itemize}
\item [$(\mathrm{H_1})$] $
:\ p,q \in C(\overline{\O}), \ 1 <p(x)<d,\ p(x)< q(x)<p^\ast(x)
,\text{ for all } x \in \overline{\O} \text{ and } 0 \leq \mu(\cdot)
\in  L^\infty(\O).
 $
    \item [$(\mathrm{H_2})$] $:\ p,q \in C^{0,1}(\overline{\O}), \ 1 <p(x)<d,\ p(x)< q(x),\text{ for all } x \in \overline{\O},\ \dfrac{q^+}{p^-}\leq 1+\dfrac{1}{d} \text{ and } 0 \leq \mu(\cdot) \in  C^{0,1}(\overline{\O}),$
        where $C^{0,1}(\overline{\O})$ refers to the space of functions that are Lipschitz continuous on $\overline{\O}.$
\end{itemize}
 \begin{prop}\cite[Proposition 2.18]{CrespoBlanco2022}\label{WW}
Let hypothesis $(\mathrm{H}_1))$ be satisfied. Then,
 \begin{itemize}
    \item [(i)] $W^{1, \mathcal{H}}(\Omega) \hookrightarrow L^{\mathcal{H}}(\Omega)$ is a compact embedding;
\item [(ii)] There exists a constant $C>0$ independent of $u$ such that
$$
\|u\|_{\lh} \leq C\|\nabla u\|_{\lh}, \quad \text { for all } u \in
W_0^{1, \mathcal{H}}(\Omega) .
$$
  \end{itemize}
 \end{prop}
\begin{prop}
  Let hypothesis $(\mathrm{H}_2)$ be satisfied. Then, $C^\infty (\O) \cap \ww$  is dense in $\ww$.
\end{prop}
Now, we give a compactness-type result inspired by the work of
Zhikov \cite{Zhikov2011} which will be useful in the sequel.
\begin{lm}\label{lem}
Let $\left\{u_n\right\}_{n \in \mathbb{N}} \subseteq W^{1, 1}_0(\O)$
and  $u \in W^{1, 1}_0(\O)$. Moreover, we assume that
\begin{equation}\label{k1}
\begin{aligned}
&1< p^-  \leq p_n(x) \leq p^+<+\infty,\\& \ 1< q^-  \leq q_n(x) \leq
q^+\ \   n \in \mathbb{N}, \text{ for a.a. } x \in \O,\end{aligned}
\end{equation}
\vspace{0.1cm}
\begin{equation}\label{k2}
p_n\rightarrow p, \ q_n \rightarrow q \text{ a.a. in } \O, \text{ as
} n\rightarrow +\infty
\end{equation}
\begin{equation}\label{k4}
\left( |\nabla u_n|^{p_n(x)}+\mu(x)  |\nabla u_n|^{q_n(x)}\right)
\text{  is bounded in } L^1(\O),
\end{equation}
\vspace{0.1cm} and
\begin{equation}\label{k3}
 \nabla u_n \rightharpoonup \nabla u \text{ in } \left(L^1(\O)\right)^d, \text{ as } n \rightarrow +\infty,
\end{equation}\vspace{0.1cm}


\vspace{0.33cm} Then $\nabla u \in \left(\lh\right)^d$ and
\begin{equation}\label{k8}
\liminf_{n \rightarrow +\infty} \int_{\O}^{} |\nabla u_n|^{p_n(x)} +
\mu(x) |\nabla u_n |^{q_n(x)} dx \geq \int_{\O}^{} |\nabla
u|^{p(x)}+ \mu(x) |\nabla u|^{q(x)}dx.
\end{equation}
\end{lm}
\begin{proof} We use some ideas coming from \cite{Chipot2019}.
 By Young's inequality one has for $\mathbf{a}, \mathbf{b_1}, \mathbf{b_2} \in \mathbb{R}^d$ and $1<r, \ s<\infty$,
\begin{equation}\label{k5}
\begin{aligned}
&\mathbf{a} \cdot \mathbf{b_1} \leq|\mathbf{a}||\mathbf{b_1}| \leq|\mathbf{a}|^r+\frac{|\mathbf{b_1}|^{r^{\prime}}}{r^{\prime} r^{\frac{r^{\prime}}{r}}}, \quad \frac{1}{r}+\frac{1}{r^{\prime}}=1,\\
&\mathbf{a} \cdot \mathbf{b_2} \leq|\mathbf{a}||\mathbf{b_2}|
\leq|\mathbf{a}|^s+\frac{|\mathbf{b_2}|^{s^{\prime}}}{s^{\prime}
s^{\frac{s^{\prime}}{s}}}, \quad \frac{1}{s}+\frac{1}{s^{\prime}}=1.
\end{aligned}
\end{equation}
If now $\mathbf{b_1}$ and $\mathbf{b_2}$ are two functions in
$\left(L^{\infty}(\Omega)\right)^d$ and we take $s=q_n$ and $r=p_n$
in \eqref{k5} and use assumption \eqref{k1}, one derives
\begin{equation}\label{k6}
\int_{\Omega}\mu(x)\left(\nabla u_n \cdot
\mathbf{b_2}-\frac{|\mathbf{b_2}|^{q_n^{\prime}(x)}}{q_n^{\prime}(x)
\left[q_n(x)\right]^{\frac{q_n^{\prime}(x)}{q_n(x)}}}\right) d x
\leq \int_{\Omega}\mu(x)\left|\nabla u_n\right|^{q_n(x)} d x,
\end{equation}
and

\begin{equation}\label{k7}
\int_{\Omega}\left(\nabla u_n \cdot
\mathbf{b_1}-\frac{|\mathbf{b_1}|^{p_n^{\prime}(x)}}{p_n^{\prime}(x)
\left[p_n(x)\right]^{\frac{p_n^{\prime}(x)}{p_n(x)}}}\right) d x
\leq \int_{\Omega}\left|\nabla u_n\right|^{p_n(x)} d x.
\end{equation}
Using assumptions \eqref{k2} and \eqref{k3} and passing to the limit
in \eqref{k6}-\eqref{k7} as $n \rightarrow \infty$, we deduce that
\begin{equation}\label{k70}
\int_{\Omega}\mu(x)\left(\nabla u \cdot
\mathbf{b_2}-\frac{|\mathbf{b_2}|^{q^{\prime}(x)}}{q^{\prime}(x)
\left[q(x)\right]^{\frac{q^{\prime}(x)}{q(x)}}}\right) d x \leq
\liminf _{n \rightarrow \infty} \int_{\Omega}   \mu(x)\left|\nabla
u_n\right|^{q_n(x)} d x:=L_1 ,
\end{equation}
and
\begin{equation}\label{k10}
\int_{\Omega}\left(\nabla u \cdot
\mathbf{b_1}-\frac{|\mathbf{b_1}|^{p^{\prime}(x)}}{p^{\prime}(x)
\left[p(x)\right]^{\frac{p^{\prime}(x)}{p(x)}}}\right) d x \leq
\liminf _{n \rightarrow \infty} \int_{\Omega}   \left|\nabla
u_n\right|^{p_n(x)} d x:=L_2.
\end{equation}
Afterward, we consider the following functions,

$$
\mathbf{b_1}=p(x)  |\nabla u|_k^{\frac{1}{p'(x)-1}} \frac{\nabla
u}{|\nabla u|},$$
and
$$
\mathbf{b_2}=q(x)|\nabla u|_k^{\frac{1}{q'(x)-1}} \frac{\nabla
u}{|\nabla u|},
$$
 with $|\nabla u|_k:=  \min\{|\nabla u|,k\}.$ Inserting $\mathbf{b_2}$ into \eqref{k70}, one gets


 $$
\int_{\Omega}\mu (x)\left(|\nabla u|_k q(x)|\nabla
u|_k^{\frac{1}{q^{\prime}(x)-1}}-|\nabla
u|_k^{\frac{q^{\prime}(x)}{q^{\prime}(x)-1}}
\frac{q(x)}{q^{\prime}(x)}\right) d x \leq L_1
$$
which implies
$$
\int_{\Omega}\mu (x)|\nabla
u|_k^{\frac{q^{\prime}(x)}{q^{\prime}(x)-1}} d x \leq L_1 .
$$
Thus
\begin{equation}\label{k9}
\int_{\Omega}\mu(x)|\nabla u|_k^{q(x)} d x \leq L_1.\end{equation}
Similarly, if we insert $\mathbf{b_1}$ in \eqref{k10}, we achieve
that
\begin{equation}\label{k11}
\int_{\Omega}|\nabla u|_k^{p(x)} d x \leq L_2.\end{equation}
Ultimately, the inequality  \eqref{k8} follows by letting $k
\rightarrow \infty$ in \eqref{k9}-\eqref{k11} , and $\nabla u \in
\left(\lh\right)^d $ due to assumption \eqref{k4}. The proof of the
lemma is complete.
\end{proof}
\section{Double phase problems with exponents dependent on solution gradient}
In this section, we present the first type of the new class of
double phase problem with exponents that depend on the gradient of
the solution:
\begin{equation}\label{t11}
\left\lbrace\begin{array}{rll}
-\operatorname{div}\left(a(x,|\nabla u|)\nabla u\right) &=f(x, u),  & \text { in } \Omega, \\
& \ & \ \\
u &=0, & \text { on } \partial \Omega,
\end{array}\right.
\end{equation}
where \( \Omega \subseteq \mathbb{R}^d, d \geq 2 \) is a bounded domain with a Lipschitz boundary \( \partial \Omega \), $f : \Omega \times \mathbb{R}\rightarrow \mathbb{R}$ is a Carathéodory function satisfying certain assumptions, $a:\O\times \R \longrightarrow (0,+\infty)$ is a function given by $$a(x,t):=|t|^{p(x,|t|)-2}+\mu (x)|t|^{q(x,|t|)-2},$$ for all $t\in \R$ and all $x\in \O,$ with \(\mu \in L^\infty (\Omega)\) such that \(\mu (x) \geq 0\) a.a. \(x \in \Omega\), and  $p, \ q :\O\times \mathbb{R}_+ \rightarrow [2,+\infty) $ are two Carathéodory functions.\\

 Note that the equation \eqref{t11}  takes the form of the following partial differential equation:
\begin{equation}\label{Ph}
\left\lbrace\begin{array}{rll}
-\Delta_\h (u) &=f(x, u),  & \text { in } \Omega, \\
\ & \ & \\
u &=0, & \text { on } \partial \Omega,
\end{array}\tag{$\mathcal{P}_\h$} \right.
\end{equation} where $\Delta_\h (u) :=\operatorname{div}\left(h(x,|\nabla u|)\nabla u\right)$, with $h$ and $\h$ will be defined in Subsection 3.1.  \\

We highlight that the presence of the solution's gradient in the
exponent introduces several challenges beyond the scope of classical
methodologies. Consequently, our initial focus will be on
elucidating the fundamental properties of the associated function
spaces. Specifically, we will delve into the requisite assumptions
on $a(\cdot,\cdot)$ to establish that the space $W^{1,
\mathcal{H}}(\Omega)$ (as detailed in Subsection 3.1) conforms to
the Musielak-Orlicz Sobolev type space structure introduced in
Section 2. Leveraging these abstract findings in conjunction with
both variational and non-variational approaches, we shall explore
the existence of nontrivial solutions to equation \ref{t11} (as
evidenced by Theorems \ref{thmA1} and \ref{thmA2}). Note that the
mentioned theorems are new in the literature.
 It is  worth noting that these theorems represent novel contributions to the field.\\
\vspace{2mm}

We suppose that the functions $p$ and $q$ satisfy the following
assumptions:
  \begin{enumerate} \item[$(\mathrm{H}_1)$] Functions $p$ and $q$ are both bounded, meaning
               $$ 2\leq p^-:= \essinf_{(x,t)\in \O\times \R_+}p(x,t) \leq p^+:=\esssup_{(x,t)\in \O\times \R_+}p(x,t)\leq d,$$
               $$2\leq q^-:= \essinf_{(x,t)\in \O\times \R_+} q(x,t)\leq q^+:=\esssup_{(x,t)\in \O\times \R_+}q(x,t)<+\infty,$$
               and $$p(x,t)\leq q(x,t),\text{ for a.a. }x\in \O \text{ and all }t\geq0.$$
                \item[$(\mathrm{H}_2)$] The functions $t \mapsto p(x,t)$ and $t \mapsto q(x,t)$  are both continuous for a.a. $x\in \O$; furthermore, they exhibit a nonincreasing  behavior for all $t \in [0,1]$ and a nondecreasing behavior for $t \geq 1.$
          \end{enumerate}

\subsection{The functional analysis of the novel space related to the new double phase operator}\label{3.1}

\vspace{0.33cm}

We introduce the function $\mathcal{H}: \O \times\mathbb{R}
\rightarrow \mathbb{R}$ as $\mathcal{H}(x,t) = \ds\int_0^{t}  h(x,s)
\, ds$, where  $h: \O\times\mathbb{R}\rightarrow \mathbb{R}$  is
defined by $$
 h(x,t)=
 \left\lbrace\begin{array}{rl}
a(x,|t|)t, &\text{ for } t\neq 0\\
 \ & \  \\
  0, &\text{ for }  t=0.
 \end{array}\right.
 $$

The following result ensures that the function $\h$ is a generalized
N-function.
\begin{prop}\label{nfct} 
Under hypotheses $(\mathrm{H}_1)$ and $(\mathrm{H}_2)$, $\h$ is a
generalized N-function.
\end{prop}
\begin{proof}
It is clear that, for each $x\in \O$, the function
\begin{equation}\label{B0}
t\mapsto h(x,t)\ \text{is continuous on}\ (0,+\infty).
\end{equation}
On the other hand,  the following inequalities hold
\begin{equation}\label{ineqs}
\begin{aligned}
  (1+\mu (x) )t^{q^+-1} &\leq h(x,t)\leq(1+\mu (x) )t^{p^--1},\  \text{for all } t\in(0,1] \text{ and a.a. } x\in \O,
     \\
   (1+\mu (x) )t^{p^--1}&\leq h(x,t)\leq(1+\mu (x) )t^{q^+-1},\  \text{for all } t\in [1,+\infty) \text{ and a.a. } x\in \O.
   \end{aligned}
\end{equation}
By \eqref{ineqs}, we obtain
\begin{equation}\label{B1}
    \lim_{t\rightarrow 0} h(x,t) = 0 \quad \text{and} \quad \lim_{t\rightarrow +\infty} h(x,t) = +\infty, \text{ for a.a. } x\in \O.
\end{equation}
Moreover, using  $(\mathrm{H}_2)$,  for all $t_2> t_1>0$ and a.a.
$x\in \Omega$, we can observe that
\begin{align*}
h(x,t_2)=t_2^{p(x,t_2)-1}+\mu(x)t_2^{q(x,t_2)-1}\geq
t_1^{p(x,t_1)-1}+\mu(x)t_1^{q(x,t_1)-1}=h(x,t_1)
\end{align*}
Thus, for each $x\in\O$, the function
\begin{equation}\label{B2}
t\mapsto h(x,t)\ \text{is increasing on } (0,+\infty).
\end{equation}
In light of \eqref{B0}, \eqref{B1},  and \eqref{B2}, we see that
$t\mapsto\mathcal{H} (x, t)$ is an N-function, for each $x\in
\Omega$. Since, for all $t\in \mathbb{R}$, the function
$x\mapsto\mathcal{H}(x, t)$ is measurable on $\O$, then,
$\mathcal{H}$ is a generalized N-function.
\end{proof}

Now, we prove an important  property of the function $\h$.

\begin{prop}\label{dlt2}Let the hypotheses $\left(\mathrm{H}_1\right)$ and $\left(\mathrm{H}_2\right)$ be satisfied. Then,
\begin{equation}\label{l2}
p^{-} \leq \frac{ h(x, t)t}{\mathcal{H}(x, t)} \leq q^+, \text{  for
all } t > 0 \text{ and a.a. }x \in \Omega.
\end{equation}
\end{prop}
\begin{proof} We define the functions $\o_1,\  \o_2:\R_+ \times \O\rightarrow \R$ by:
$$\left.
\begin{array}{ll}
\omega_1(t) & :=\h(x,t)-\frac{1}{p^-}h(x,t)t, \\
& \ \\
\omega_2(t) &:=\h(x,t)-\frac{1}{q^+} h(x,t)t,
\end{array}
\right\} \text{ for all }  t\in [0,+\infty) \text{ and a.a. }x \in
\Omega.
$$
 From assumptions $\left(\mathrm{H}_1\right)$ and
$\left(\mathrm{H}_2\right)$, we can see that $t  \mapsto
\dfrac{h(x,t)}{t^{p^--1}}$ is nondecreasing on $(0,+\infty)$. Hence,
for $t>\e>0$, we get
$$\begin{aligned}
\o_1(t)& =\h(x,t)-\frac{1}{p^-}h(x,t)t +\h(x,\e)-\h(x,\e)
= \h(x,\e) -\frac{1}{p^-}h(x,t)t  +\int_{\epsilon}^t h(x,s)ds\\
&=  \h(x,\e) -\frac{1}{p^-}h(x,t)t  +\int_{\epsilon}^t \frac{h(x,s)}{s^{p^--1}}s^{p^--1}ds\\ &\leq \h(x,\e) -\frac{1}{p^-}h(x,t)t  +\frac{h(x,t)}{t^{p^--1}}\int_{\epsilon}^t s^{p^--1}ds\\
&=\h(x,\e) -\frac{1}{p^-}h(x,t)t  +\frac{h(x,t)}{t^{p^--1}}\frac{1}{p^-} \left[ t^{p^-}-\e^{p^-} \right]=  \h(x,\e) -\frac{h(x,t)\e^{p^-}}{t^{p^--1}}\frac{1}{p^-} \\
&\leq \h(x,\e)  -\frac{1}{p^-}h(x,\e)\e=\o_1(\e).
\end{aligned}
$$
Thus, for $t>\e>0$, we infer that
\begin{equation}\label{SB}\o_1(t)\leq \o(\e).\end{equation} Tending
$\e$ to zero in \eqref{SB}, we deduce that $\o_1(t)\leq \o_1(0)=0$,
namely
\begin{equation}\label{SB1}
p^{-} \leq \frac{ h(x, t)t}{\mathcal{H}(x, t)}, \text{  for all } t
> 0 \text{ and a.a. }x \in \Omega.
\end{equation}
Again, from assumptions $\left(\mathrm{H}_1\right)$ and
$\left(\mathrm{H}_2\right)$, we can see that $t  \mapsto
\dfrac{h(x,t)}{t^{q^+-1}}$ is nonincreasing on $(0,+\infty)$. Hence,
for $t>\e>0$, one has
$$\begin{aligned}
\o_2(t)&=\h(x,t)-\frac{1}{q^+}h(x,t)t +\h(x,\e)-\h(x,\e)= \h(x,\e) -\frac{1}{q^+}h(x,t)t  +\int_{\epsilon}^t h(x,s)ds\\
&=  \h(x,\e) -\frac{1}{q^+}h(x,t)t  +\int_{\epsilon}^t \frac{h(x,s)}{s^{q^+-1}}s^{q^+-1}ds\\
&\geq \h(x,\e) -\frac{1}{q^+}h(x,t)t  +\frac{h(x,t)}{t^{q^+-1}}\int_{\epsilon}^t s^{q^+-1}ds\\
&=\h(x,\e) -\frac{1}{q^+}h(x,t)t  +\frac{h(x,t)}{t^{q^+-1}}\frac{1}{q^+} \left[ t^{q^+}-\e^{q^+} \right]=  \h(x,\e) -\frac{h(x,t)\e^{q^+}}{t^{q^+-1}}\frac{1}{q^+} \\
&\geq \h(x,\e)  -\frac{1}{q^+}h(x,\e)\e=\o_2(\e).
\end{aligned}
$$
Therefore,  \begin{equation}\label{SB0}\o_2(t)\geq
\o_2(\e),\end{equation} which implies, by tending $\e$ to zero in
\eqref{SB0},  that $\o_2(t)\geq \o_2(0)=0,$ namely
\begin{equation}\label{SB2}
     \frac{ h(x, t)t}{\mathcal{H}(x, t)} \leq q^+, \text{  for all } t > 0 \text{ and a.a. }x \in \Omega.
\end{equation}
Consequently, by combining \eqref{SB1} and \eqref{SB2}, we conclude
the proof of \eqref{l2}. This finishes the proof.
\end{proof}

 \begin{prop} \label{delta2} 
Assume that the hypotheses $\left(\mathrm{H}_1\right)$ and
$\left(\mathrm{H}_2\right)$ are satisfied. Then, both  $\h$ and
$\widetilde{\h}$ fulfill the $\Delta_2$-condition.
 \end{prop}
\begin{proof}
 The proof follows by combining  Proposition \ref{dlt2} and Remark \ref{compl}.
\end{proof}


 Thanks to Propositions \ref{nfct} and \ref{delta2}, we     are ready to define the  Musielak–Orlicz space $\lh$ as follows:
$$\lh=\big\{ u:\O \rightarrow \R \text{ measurable}; \vi(u)<+\infty \big\},$$  endowed with the following norm
$$\|u\|_{\lh} =\inf \bigg\{ \l >0, \ \vi\big(\frac{u}{\l}\big) \leq1\bigg\}, $$ where $\vi$ is defined as in \eqref{Mo}. Proceeding
similarly to Subsection 2.1, we can introduce the spaces $\ww$ and
$W^{1,\h}_0(\O)$ equipped with the norm
$$
\|u\|_{\ww}:=\|u\|_{\lh}+\|\nabla u\|_{\lh} .
$$
\begin{prop}\label{separble}
   Let hypotheses $\left(\mathrm{H}_1\right)$ and $\left(\mathrm{H}_2\right)$ hold. Then, the space $\lh$ is a separable, reflexive Banach space.
\end{prop}
\begin{proof}
Invoking Propositions \ref{AB}, \ref{nfct} and \ref{delta2},  we get
the proof of our desired result.
\end{proof}
 As a consequence of Propositions \ref{zoo} and \ref{delta2}, we have the following result.
\begin{prop}\label{rel}
 Suppose that condition $\left(\mathrm{H}_1\right)$ and $\left(\mathrm{H}_2\right)$ are satisfied. The $\mathcal{H}$-modular has the following properties:
 \begin{itemize}
                \item[$(\mathrm{i})$] If $\|u\|_{\lh}<1,$ then $\|u\|_{\lh}^{q^+}\leqslant \rho_{\mathcal{H}}(u) \leqslant\|u\|^{p^-}_{\lh}$.
                \item  [$(\mathrm{ii})$] If $\|u\|_{\lh}>1,$ then $\|u\|^{p^-}_{\lh} \leqslant \rho_{\mathcal{H}}(u) \leqslant\|u\|_{\lh}^{q^+}.$
                \item  [$(\mathrm{iii})$]  $u_n \rightarrow u$ in $\lh$ $\Longleftrightarrow$ $\rho_{\mathcal{H}}(u_n-u)\rightarrow 0$.
 \end{itemize}
\end{prop}
According to Theorem \ref{thc}, we have the following result.
\begin{prop}\label{rflxv}
Assume that $(\mathrm{H}_1)$ and $(\mathrm{H}_2)$  hold. Then,
$W^{1, \mathcal{H}}(\Omega)$ and  $W_0^{1, \mathcal{H}}(\Omega)$ are
both separable and reflexive Banach spaces.
\end{prop}
\begin{proof}
  By inequality \eqref{ineqs}, we can prove that $$\inf_{x\in \O} \mathcal{H}(x,1) >0.$$
  It follows, by Proposition \ref{separble} and Theorem \ref{thc}, that  $W^{1, \mathcal{H}}(\Omega)$ and  $W_0^{1, \mathcal{H}}(\Omega)$ are both separable and reflexive Banach spaces.
\end{proof}
\begin{prop} \label{prprp} Suppose that hypotheses $\left(\mathrm{H}_1\right)$ and $\left(\mathrm{H}_2\right)$  are satisfied.
 Let $u, v \in L^{\mathcal{H}}(\Omega)$, then we have
$$
\int_{\Omega}|u|^{p(x,|u|)-2} u v d x<+\infty\ \text{, }\
\int_{\Omega} \mu(x)|u|^{q(x,|u|)-2} u v d x<+\infty,
$$
and
$$
\int_{\Omega} h(x,u)vdx \leq 2q^+
\left(\rho_{\mathcal{H}}(u)+\rho_{\mathcal{H}}(v) \right).$$
\end{prop}
\begin{proof} Let $u, v \in L^{\mathcal{H}}(\Omega)$. Since the function $t \mapsto t^{p(x,t)}$ is nondecreasing on $[0,+\infty)$, then
$$
\begin{aligned}
\left|\int_{\Omega}|u|^{p(x,|u|)-2} u v d x \right|& \leq \int_{\Omega}|u|^{p(x,|u|)-1}|v| d x \\
& =\int_{\{|u| \geq|v|\}}|u|^{p(x,|u|)-1}|v| d x+\int_{\{|u| \leq|v|\}}|u|^{p(x,|u|)-1}|v| d x \\
& \leq \int_{\{|u| \geq|v|\}}|u|^{p(x,|u|)} d x+\int_{|u| \leq|v|}|v|^{p(x,|v|)} d x \\
& \leq \int_{\Omega}|u|^{p(x,|u|)} d x+\int_{\Omega}|v|^{p(x,|v|)} d x \\
& \leq \int_{\Omega}|u|^{p(x,|u|)}+\mu(x)|u|^{q(x,|u|)}  d x+\int_{\Omega}|v|^{p(x,|v|)}+\mu(x)|v|^{q(x,|v|)} d x \\
& \leq
q^+\left(\rho_{\mathcal{H}}(u)+\rho_{\mathcal{H}}(v)\right)<+\infty.
\end{aligned}
$$
Similarly, we prove that
$$\int_{\Omega} \mu(x)|u|^{q(x,|u|)-2} u v d x\leq q^+\left(\rho_{\mathcal{H}}(u)+\rho_{\mathcal{H}}(v)\right)<+\infty.$$
The proof is completed.
\end{proof}
To obtain new properties of $\ww$ (specifically the density and the
embedding), we assume additional conditions on $p$ and $q$:
\vspace{0.33cm}
 \begin{enumerate}
\item[$(\mathrm{H}_2')$] The functions $t \mapsto p(x,t)$ and $t \mapsto q(x,t)$  are both continuous for a.a. $x\in \O$; furthermore, they exhibit a constant behavior equal $p(x)$ and $q(x)$, respectively, for all $t \in [0,1]$ and a nondecreasing behavior for $t \geq 1.$
\item[$(\mathrm{H}_3)$] The function $x\mapsto \mu(x)$ is $C^{0,1}(\overline{\O}).$ Moreover, there exist $c_p, c_q>0$ such that
 $$|p(x,t)-p(y,t)| \leq c_p|x-y| \text{ and } |q(x,t)-q(y,t)| \leq c_q|x-y|, $$ $\text{ for all }  t \geq0  \text{ and for all } x,y \in \O.$
\vspace{2mm}
                  \item[$(\mathrm{H}_4)$] $\frac{q^+}{p^-}<1+\frac{1}{d}.$
               \end{enumerate}
               \vspace{3.3mm}
The hypothesis $(\mathrm{H}_2)$ is weaker than
$\left(\mathrm{H}_2^{\prime}\right)$. Then, if we replace
$(\mathrm{H}_2)$ by $\left(\mathrm{H}_2^{\prime}\right)$, the
previous results remain valid. \vspace{2mm}

Now, we discuss the compact  embedding of $\ww$ and $\w$ into the
Museilak space $L^A(\O)$ by applying  the new optimal theorem
obtained in \cite{Cianchi2023}. To this end, we introduce some
definitions and conditions found in \cite{Harjulehto2019}.
\vspace{0.33cm}
\begin{dfn}~~~~~~~~~~~~~~~~~~~~~~~~~~~~~~~~~~~~~~~~~~~~~~~~~~~~~~~~~

\begin{itemize}
\item  [(i)] We call a function $g:(0, \infty) \rightarrow \mathbb{R}$ almost increasing if there exists a constant $a \geq 1$ such that $g(s) \leq a g(t)$ for all $0<s<t$.
Similarly, we define almost decreasing functions. \vspace{0.33cm}
\item  [(ii)] We say that $\varphi: \Omega \times[0, \infty) \rightarrow[0, \infty]$ is a $\Phi$-prefunction if $x \mapsto \varphi(x,|g(x)|)$ is measurable for every
 measurable function $g: \Omega \rightarrow \mathbb{R}, \varphi(x, 0)=0$, $t\mapsto \varphi (x,t)$ is nondecreasing for a.a. $x\in \O,$
$$
\lim _{t \rightarrow 0^{+}} \varphi(x, t)=0 \text { and } \lim _{t
\rightarrow \infty} \varphi(x, t)=\infty \quad \text { for a.a. } x
\in \Omega .
$$
If in addition the condition
$$
\frac{\varphi(x, t)}{t} \text { is almost increasing for a.a. } x
\in \Omega
$$
is satisfied, then $\varphi$ is called a weak $\Phi$-function and
the class of all weak $\Phi$-functions is denoted by
$\Phi_{\mathrm{w}}(\Omega)$. \vspace{0.33cm}
\item  [(iii)] We say that $\varphi \in \Phi_{\mathrm{w}}(\Omega)$ satisfies $(\mathrm{A}_0)$, if there exists a constant $\beta \in(0,1]$ such that $\beta \leq \varphi^{-1}(x, 1) \leq \frac{1}{\beta}$ for a.a. $x \in \Omega$.
\vspace{0.33cm}
\item [(iv)] We say that $\varphi \in \Phi_{\mathrm{w}}(\Omega)$ satisfies $(\mathrm{A}_1)$, if there exists $\beta \in(0,1)$ such that
$$
\beta \varphi^{-1}(x, t) \leq \varphi^{-1}(y, t)
$$
for every $t \in\left[1, \frac{1}{|B|}\right]$, for a.a. $x, y \in B
\cap \Omega$ and for every ball $B$ with $|B| \leq 1$.
\vspace{0.33cm}
\item [(v)] We say that $\varphi \in \Phi_{\mathrm{w}}(\Omega)$ satisfies $\left(\mathrm{A}^{\prime}_ 1\right)$, if there exists $\beta \in(0,1)$ such that
$$
\varphi(x, \beta t) \leq \varphi(y, t)
$$
for every $\varphi(y, t) \in\left[1, \frac{1}{|B|}\right]$, for a.a.
$x, y \in B \cap \Omega$ and for every ball $B$ with $|B| \leq 1$.
\vspace{0.33cm}
\item [(vi)] We say that $\varphi \in \Phi_{\mathrm{w}}(\Omega)$ satisfies $(\mathrm{A}_2)$, if for every $s>0$ there exist $\beta \in(0,1]$ and $h \in L^1(\Omega) \cap L^{\infty}(\Omega)$ such that
$$
\beta \varphi^{-1}(x, t) \leq \varphi^{-1}(y, t),
$$
for a.a. $x,y \in \Omega$ and for all $t \in[h(x)+h(y), s]$.
\vspace{0.33cm}
\item [(vii)] We say that $\varphi: \Omega \times(0, \infty) \rightarrow \mathbb{R}$ satisfies $(\mathrm{aDec})$ if there exists $\a \in(0, \infty)$ such that
$$
\frac{\varphi(x, t)}{t^{\a}} \text { is almost decreasing for a.a. }
x \in \Omega \text {. }
$$
\end{itemize}

\vspace{0.33cm} In the sequel we will use $g_1 \approx g_2$ and $g_1
\lesssim g_2$ if there exist constants $c_1, c_2>0$ such that $c_1
g_1 \leq g_2 \leq c_2 g_1$ and $g_1 \leq c_2 g_2$, respectively.
\end{dfn}
Next, we give a sufficient conditions for $\h$ to satisfy conditions
$(\mathrm{A}_0)$, $(\mathrm{A}_1)$, $(\mathrm{A}_2)$ and
$(\mathrm{aDec})$ (see \cite{CrespoBlanco2022} and
\cite{Harjulehto2019}).
\begin{prop}\label{eqs}
  Let $\varphi \in \Phi_{\mathrm{w}}(\Omega)$, we have the following assertions:
  \vspace{0.33cm}
  \begin{itemize}
    \item [$(\mathrm{i})$] If there exists $c>0$ such that $\varphi (x,c) \approx1$, then $\varphi$ satisfies $(\mathrm{A}_0),$ (see \cite[Corollary 3.7.5]{Harjulehto2019}).
    \vspace{0.33cm}
    \item [$(\mathrm{ii})$] If  $\varphi$ satisfies $(\mathrm{A}_0)$. Then $\varphi$ satisfies $(\mathrm{A}_1)$ if and only if $\varphi$ satisfies $(\mathrm{A}^\prime_1),$ (see
     \cite [Corollary 4.1.6]{Harjulehto2019}).
    \vspace{0.33cm}
    \item [$(\mathrm{iii})$]If $\O$ is bounded, then $\varphi$ satisfies $(\mathrm{A}_2),$ (see \cite[Lemma 4.2.3]{Harjulehto2019}).
    \vspace{0.33cm}
    \item [$(\mathrm{iv})$] The $\Delta_2$-condition is equivalent to $(\mathrm{aDec}),$ (see \cite[ Lemma 2.2.6]{Harjulehto2019}).
  \end{itemize}
\end{prop}
\begin{lm}\label{lm0}
Let hypotheses $\left(\mathrm{H}_1\right)$,
$\left(\mathrm{H}_2'\right)$, $\left(\mathrm{H}_3\right)$ and
$\left(\mathrm{H}_4\right)$ be satisfied. Then, $\mathcal{H}$
satisfies $\left(\mathrm{A}_1^{\prime}\right)$.

\end{lm}
\begin{proof}

 We employ the methodologies used in the proof of  \cite[Proposition 2.23]{CrespoBlanco2022}.
From Proposition \ref{nfct} and \cite{Harjulehto2019}, we can
observe that $\mathcal{H} \in \Phi_{\mathrm{w}}(\O).$
We divide the rest of the proof into three steps.\\

    \begin{itemize}

      \item [\textbf{Claim 1.}] The condition $\mathcal{H}(y, t) \in\left[1, \frac{1}{|B|}\right]$ implies that
\begin{equation}\label{aq2}
t
\in\left[\left(1+\|\mu\|_{\infty}\right)^{-\frac{1}{p^{-}}},\left(q^+\right)^{\frac{1}{p^{-}}}|B|^{-\frac{1}{p^{-}}}\right],
\end{equation}
for a.a. $y \in B \cap \Omega$ and for every ball $B$ with $|B| \leq
1$. To this end, let $B \subset \R^d$ such that $0<|B|<1$ and
$\mathcal{H}(y, t) \in\left[1, \frac{1}{|B|}\right]$ for a.a. $y \in
B$. Then, by Proposition \ref{dlt2}, we have $$\mathcal{H}(y, t)
\leq p^{-} \mathcal{H}(y, t) \leq h(y, t) t \leq q^+ \mathcal{H}(y,
t), \text{  for a.a. } y \in B.$$ Thus, it follows that
$$
t^{p^{-}}\left(1+\|\mu\|_{\infty}\right) \geq \begin{cases}t^{p(y, t)}+\mu(y) t^{q(y, t)} \geq \mathcal{H}(y, t) \geq 1, & \text { if } t \leq 1 \\
\ & \\ 1, & \text { if } t \geq 1\end{cases}
$$
and
$$
t^{p^{-}} \leq \begin{cases}t^{p(y, t)} \leq t^{p(y, t)}+\mu(y)
t^{q(y, t)} \leq q^+\mathcal{H}(y, t) \leq \dfrac{q^+}{|B|}, & \text
{ if } t \geq 1 \\ \ & \\ \dfrac{q^+}{|B|}, & \text { if } t \leq 1
.\end{cases}
$$
Thus, according to $\mathcal{H}(y, t) \in\left[1,
\frac{1}{|B|}\right]$, we conclude \eqref{aq2}.

      \item [\textbf{Claim 2.}] If $t$ verifies \eqref{aq2}, we have \begin{equation}\label{aq4}
(\xi t) h(x, \xi t) \leq  h(y, t)t,
\end{equation}
for some $\xi \in(0,1)$, for a.a. $x, y \in B \cap \Omega$ and for every ball $B$ with $|B| \leq 1$.\\
 Indeed: Let us fix a ball $B \subset \R^d$ of radius $r > 0$ such that $|B| \leq 1$ (in particular, $r < 1)$. Notice that $|B| = \alpha(d) r^d$, where $\alpha(d) > 1$ is a constant that depends only on the dimension $d$.\vspace{2mm}\\
Firstly, we prove that: For all $t$ verifies \eqref{aq2}, all $\xi
\in(0,1)$, and a.a. $x, y \in \Omega$ such that $|x-y| \leq 2 r$,
one has
\begin{equation}\label{h4}
t^{p(x,\xi t)} \leq M \cdot t^{p(y,t)} \text { and } t^{q(x,\xi t)}
\leq M \cdot t^{q(y,t)},
\end{equation}
for some constant $M = M(d, p^-, q^+, \mu) \geq 1$ not depending on $x, y, t$.\\
Fix $t$ satisfying \eqref{aq2} and $x, y \in \Omega$ such that
$|x-y| \leq 2 r$. By $\left(\mathrm{H}_3\right)$ , we have
\begin{equation}\label{h5}
|p(x,t)-p(y,t)| \leq c_p|x-y| \leq 2 r c_p.
\end{equation}
We consider three cases:
\begin{itemize}
  \item [\textbf{Case }I:] If $t \leq 1$ and $p(x, \xi t) \geq p(y, t)$ or $t \geq 1$ and $p(x, \xi t) \leq p(y,t )$, then the first inequality in \eqref{h4} holds with $M=1$. \\
 \item [\textbf{Case }II:] If $t \leq 1$ and $p(x, \xi t) \leq p(y ,  t)$. Then, by applying \eqref{aq2}, it follows that
$$\begin{aligned}
t^{p(x, \xi t)}=t^{p(x, \xi t)-p(y,  t)} t^{p(y,  t)}
&\leq\left(\left(1+\|\mu\|_{\infty}\right)^{\frac{1}{p-}}\right)^{p(y,
t)-p(x, \xi t)} t^{p(y,  t)}\\&
\leq\left(1+\|\mu\|_{\infty}\right)^{\frac{p^{+}}{p-}} t^{p(y,t)} .
\end{aligned}$$
Thus, the first inequality in \eqref{h4} holds with $M=\left(1+\|\mu\|_{\infty}\right)^{\frac{p^{+}}{p^{-}}}$.\\
 \item [\textbf{Case }III:]
  If $t \geq  1$ and $p(x, \xi t) \geq p(y, t)$. By using \eqref{aq2}, \eqref{h5} and assumption $(\mathrm{H}_2')$, we
  get
$$\begin{aligned}
t^{p(x, \xi t)}=t^{p(x , \xi t)-p(y,  t)} t^{p(y,  t)}&\leq t^{p(x ,
t)-p(y,  t)} t^{p(y,  t)}
\\&\leq\left((q^+)^{\frac{1}{p^-}}\alpha(d)^{-\frac{1}{p^{-}}}
r^{-\frac{d}{p^{-}}}\right)^{2 r c_p} t^{p(y, t)} \\
&\leq\left((q^+)^{\frac{2c_p}{p^-}}\alpha(d)^{-\frac{2
c_p}{p^{-}}}\right)^r\left(r^r\right)^{-\frac{2 d c_p}{p^{-}}}
t^{p(y,t)}.
\end{aligned}$$
Note that the function
$\delta(r)=\left((q^+)^{\frac{2c_p}{p^-}}\alpha(d)^{-\frac{2
c_p}{p^-}}\right)^r\left(r^r\right)^{-\frac{2 d c_p}{p^{-}}}$is
strictly positive and continuous on the interval $\left[0,
\frac{1}{\alpha(d)^{\frac{1}{d}}}\right]$ where $\delta(0)=1$. Hence
it attains its maximum at some $r_0 \in\left[0,
\frac{1}{\alpha(d)^{\frac{1}{d}}}\right]$. Then the first inequality
in \eqref{h4}  holds for $M=\delta\left(r_0\right)\geq
\delta(0)=1$.\end{itemize}
By the same argument used above, we show the second inequality in
\eqref{h4}. Then, taking M as the maximum,
we conclude the proof of \eqref{h4}.\\
Now, we are ready to prove \eqref{aq4}. Since $\mu \in
C^{0,1}(\overline{\Omega})$ and $|x-y| \leq 2 r$, one has
\begin{equation}\label{h6}
|\mu (x)-\mu (y)| \leq c_\mu|x-y| \leq 2 c_\mu  r.
\end{equation}
We start with the left hand side of \eqref{aq4}. Since $\xi
\in(0,1)$ and taking \eqref{h4} as well as \eqref{h6} into account,
we get
\begin{equation}\label{h7}
\begin{aligned}
(\xi t)^{p(x,\xi t)}+\mu (x)(\xi t)^{q(x,\xi t)} & \leq \xi^{p-} t^{p(x,\xi t)}+\mu (x) \xi^{p-} t^{q(x,\xi t)} \\
& \leq \xi^{p^{-}} M\left(t^{p(y,t)}+\mu (x) t^{q(y,t)}\right) \\
& \leq \xi^{p-} M\left(t^{p(y,t)}+\mu (y) t^{q(y,t)}+2 c_\mu  r t^{q(y,t)}\right) \\
& \leq \xi^{p-} M\left(t^{p(y,t)}+2 c_\mu  r t^{q(y,t)}\right)+\mu
(y) t^{q(y,t)}
\end{aligned}
\end{equation}
where the last inequality holds providing
$\xi<M^{-\frac{1}{p^{-}}}$. Continuing \eqref{h7} and by
\eqref{aq2}, we obtain
\begin{equation}\label{h8}
\begin{aligned}
(\xi t)^{p(x,\xi t)}+\mu (x)(\xi t)^{q(x,\xi t)} & \leq \xi^{p^{-}} M t^{p(y,t)}\left(1+2 c_\mu  r t^{q(y,t)-p(y,t)}\right)+\mu (y) t^{q(y,t)} \\
& \leq \xi^{p^{-}} M t^{p(y,t)}\left(1+2 c_\mu  r\left((q^+)^{\frac{1}{p^-}}\alpha(d)^{-\frac{1}{p^{-}}} r^{-\frac{d}{p^{-}}}\right)^{q^{+}-p^{-}}\right)\\&\ \  \ +\mu (y) t^{q(y,t)} \\
& =\xi^{p^{-}} M t^{p(y,t)}\left(1+2 c_\mu
\alpha(d)^{-\frac{q^{+}}{p^{-}}+1} (q^+)^{\frac{q^{+}}{p^{-}}-1}
r^{1+d-d \frac{q^{+}}{p^{-}}}\right)\\&\ \ \ +\mu (y) t^{q(y,t)}.
\end{aligned}
\end{equation}
Recall that $r<1$ and $1+d-d \frac{q^{+}}{p^{-}} >   0$
$(\mathrm{H}_4)$. It follows, in view of \eqref{h8}, that
\begin{equation}\label{YR1}
(\xi t)^{p(x,\xi t)}+\mu (x)(\xi t)^{q(x,\xi t)} \leq \xi^{p-} M
t^{p(y,t)}\left(1+2 c_\mu  (q^+)^{\frac{q^{+}}{p^{-}}-1}
\alpha(d)^{-\frac{q^{+}}{p^{-}}+1}\right)+\mu (y) t^{q(y,t)}.
\end{equation}
Choosing $\xi>0$ small enough in \eqref{YR1}, namely
$$
\xi<M^{-\frac{1}{p-}}\left(1+2 c_\mu  (q^+)^{\frac{q^{+}}{p^{-}}-1}
\alpha(d)^{-\frac{q^{+}}{p^{-}}+1}\right)^{-\frac{1}{p^{-}}},
$$
we deduce the proof of \eqref{aq4}. Note that the choice of $\xi$
depends only on $M, d, p^-, q^+$ and $\mu$.

Finally, let $t$ verifies \eqref{aq2} and \eqref{aq4}, for some $\xi
\in(0,1)$, for a.a. $x, y \in B \cap \Omega$ and for every ball $B$
with $|B| \leq 1$. Then, by \eqref{l2}, one has
 \begin{equation}\label{YR}
q^+ \mathcal{H}(y, t) \geq t h(y, t) \geq(\xi t) h(x, \xi t) \geq
\mathcal{H}(x, \xi t).
 \end{equation}
 Therefore, using \eqref{YR} and the fact that $\mathcal{H} (x,\cdot)$ is convex, we
 deduce that
$$
\mathcal{H}(y, t) \geq \frac{\mathcal{H}(x, \xi t)}{q^+} \geq
\mathcal{H}\left(x, \frac{\xi}{q^+} t\right).
$$
Hence, $\mathcal{H}$ satisfies $\left(\mathrm{A}_1^{\prime}\right)$ with $\beta=\dfrac{\xi}{q^+}$.\\
\vspace{0.33cm}
    \end{itemize}
    This completes the proof.
\end{proof}

\begin{prop} \label{AAA}
     Let hypotheses $\left(\mathrm{H}_1\right)$, $\left(\mathrm{H}_2'\right)$, $\left(\mathrm{H}_3\right)$ and $\left(\mathrm{H}_4\right)$ be satisfied, then $\mathcal{H}$ satisfies $(\mathrm{A}_0)$, $(\mathrm{A}_1)$, $(\mathrm{A}_2)$ and $(\mathrm{aDec})$.
\end{prop}
\begin{proof}
First, in view of the monotonocity of $t\rightarrow h(x,t)$, we
deduce that:
$$
 1 \leq  \int_1^2 h(x,1)ds \leq  \int_1^2 h(x,s)ds \leq \mathcal{H}(x, 2) \leq \int_0^2 h(x,2)ds \leq2^{q^+}(1 + \|\mu\|_\infty) \cdot 1.$$
Thus, $\mathcal{H}(x, 2) \approx 1$. Therefore, according to
Proposition \ref{eqs}, $\mathcal{H}$ satisfies $(\mathrm{A}_0)$.
Next, by the boundedness of $\O$ and  applying the third assertion
of Proposition \ref{eqs}, we show that $\mathcal{H}$ satisfies
$(\mathrm{A}_2)$ (see Proposition \ref{dlt2}). Again, invoking
Proposition \ref{eqs} and the fact that $\mathcal{H}$ satisfies the
$\Delta_2$-condition, we see that the assumption $(\mathrm{aDec})$
is fulfilled.  Finally, according to the definition of $\h$,
Proposition \ref{eqs}(ii) and condition $(\mathrm{A}_0)$, we need
only to show that $\mathcal{H}$ satisfies $(\mathrm{A}' _1)$.
Therefore, by Lemma \ref{lm0}, we conclude that $\h $ satisfies
$(\mathrm{A} _1)$. This ends the proof.
\end{proof}


We would like to note that, in our case, the explicit form of $\h$
is unknown. This limitation prevents us from applying \cite[Theorems
1.1 and 1.2]{Fan2012a} to establish compact embedding results.
Consequently, our approach will be based on the recent work by A.
Cianchi and L. Diening \cite{Cianchi2023}. First, let us review some
key definitions provided in \cite{Cianchi2023}.

\begin{dfns}
\begin{itemize}
  \item [$(\mathrm{i})$]

The Sobolev conjugate of $\mathcal{H}$ is defined as the generalized
N-function $\mathcal{H}_\ast$ obeying
$$
\mathcal{H}_\ast(x, t):=\widehat{\mathcal{H}}\left(x,
\mathcal{N}^{-1}(x, t)\right), \quad \text { for } x \in \O \text {
and } t \geq 0,
$$
where $\mathcal{N}: \O \times[0, +\infty) \rightarrow[0,+ \infty)$
is the function given by
$$
\mathcal{N}(x,
t):=\left(\int_0^t\left(\frac{\tau}{\widehat{\mathcal{H}}(x,
\tau)}\right)^{\frac{1}{d-1}} d \tau\right)^{\frac{1}{d^{\prime}}}
,\text { for } x \in \O \text { and } t \geq 0,
$$
with

$$\widehat{\mathcal{H}}(x,t):=\begin{cases}
    2\max\left\{ \h\left(x, \h^{-1}(x,1)t\right), 2t-1 \right\}-1
& \text{if } t \geq 1\\
\ & \ \\
    t  & \text{if } 0\leq t <1,
\end{cases}$$  for all $x \in \O$.

\item [$(\mathrm{ii})$] Assume that $D$ is an open set in $\mathbb{R}^d$ and $\varphi$ is a generalized N-function. The homogeneous Musielak-Orlicz-Sobolev space $V^{1, \varphi}(D)$ is defined as
$$
V^{1, \varphi}(D)=\left\{u \in W_{\text {loc }}^{1,1}(D):|\nabla u|
\in L^{\varphi}(D)\right\} .
$$
If $D$ is connected and $G$ is a bounded open set such that
$\overline{G} \subset D$, then the functional
$$
\|u\|_{L^1(G)}+\|\nabla u\|_{L^{\varphi}(D)}
$$
defines a norm on $V^{1, \varphi}(D)$.
\item [$(\mathrm{iii})$] The subspace of functions which vanish on $\partial D$ is suitably given by
$$
V_0^{1, \varphi}(D)=\left\{u \in W_{\text {loc }}^{1,1}(D) \text { :
the extension to } \mathbb{R}^d \text { of } u \text { by } 0 \text
{ outside } D \text { belongs to } V^{1,
\varphi}\left(\mathbb{R}^d\right)\right\},
$$
and can equipped with the norm
$$
\|\nabla u\|_{L^{\varphi}(D)} .
$$
\item [$(\mathrm{iv})$] An open set
$D \subset \R^d$ is called a bounded John domain if there exist
positive constants $\theta$ and $\delta$, and a point $x_0 \in \O$
such that for every $x \in D$ there exists a rectifiable curve $\psi
: [0, \delta] \rightarrow D$, parametrized by arclength, such that
$\psi(0) = x, \psi(\delta) = x_0$, and
$$\operatorname{ dist} (\psi(r), \partial D) \geq \theta r \text{ for } r \in [0, \delta].$$
\end{itemize}
\end{dfns}
\begin{remarks}\label{KK}
  \begin{itemize}
    \item [(1)] We can observe that \( \widehat{\h} \approx \h \). Therefore, according to \cite[Remark 3.8]{Cianchi2023}, we can substitute \( \h \) by \( \widehat{\h} \) in the definitions of \( \mathcal{N} \) and $\h_\ast$.

    \item [(2)] Every bounded Lipschitz domain is a John domain. Thus, $\O$ is a John domain.
  \end{itemize}
\end{remarks}
Now, we are ready to give some embedding results of the new spaces
$\ww$ and $\w$ into Musielak-Orlicz spaces.
\begin{theorem} \label{inject10}
   Let hypotheses $\left(\mathrm{H}_1\right)$, $\left(\mathrm{H}_2'\right)$, $\left(\mathrm{H}_3\right)$ and $\left(\mathrm{H}_4\right)$ be satisfied and let $\vartheta $  be a generalized N-function.
   \vspace{2mm}
   \begin{itemize}
  \item [$\mathrm{(i)}$] If $\vartheta\preceq \mathcal{H}_\ast$, then the continuous  embeddings $W^{1,\mathcal{H}} (\O) \hookrightarrow L^{\vartheta}(\Omega)$ and $W^{1,\mathcal{H}}_0 (\O) \hookrightarrow L^{\vartheta}(\Omega)$ hold.
      \vspace{2mm}
    \item [$\mathrm{(ii)}$] If $\vartheta\ll \mathcal{H}_\ast ,$  then the compact  embeddings $W^{1,\mathcal{H}} (\O) \hookrightarrow L^{\vartheta}(\Omega)$ and $W^{1,\mathcal{H}}_0 (\O) \hookrightarrow L^{\vartheta}(\Omega)$ hold. \end{itemize}
\end{theorem}
\vspace{2mm}
\begin{proof}
\begin{itemize}
    \item  [$(\mathrm{i})$] Let $G$ be a bounded set such that $\overline{G} \subset \O$. Since $\lh \hookrightarrow L^1(\O)$ (see \cite[Proposition 1.5]{Fan2012b}), there exists $C>0$ such that  $$\|u\|_{L^1(G)}+\|\nabla u\|_{L^{\mathcal{H}}(\O)}\leq \|u\|_{L^1(\O)}+\|\nabla u\|_{L^{\mathcal{H}}(\O)}\leq C \left( \|u\|_{L^{\mathcal{H}}(\O)}+\|\nabla u\|_{L^{\mathcal{H}}(\O)} \right).$$ Therefore, \begin{equation}\label{FA} W^{1,{\mathcal{H}}}(\O) \hookrightarrow V ^{1,{\mathcal{H}}}(\O)  \text{ and }W^{1,{\mathcal{H}}}_0(\O)\hookrightarrow V^{1,{\mathcal{H}}}_0(\O).\end{equation} By applying Proposition \ref{AAA} and \cite[Theorem 3.5]{Cianchi2023}, we can confirm that the continuous embedding $W^{1,\mathcal{H}}_0 (\O) \hookrightarrow L^{\vartheta}(\Omega)$ holds. Furthermore, in view of Remarks \ref{KK}-(2) and  \cite[Theorem 3.6.]{Cianchi2023}, we can also assert the continuous embedding $W^{1,\mathcal{H}} (\O) \hookrightarrow L^{\vartheta}(\Omega)$.
        \vspace{2mm}
     \item  [$(\mathrm{ii})$] In light of \eqref{FA} and \cite[Theorem 3.7]{Cianchi2023}, we conclude the proof of  $(\mathrm{ii})$.
\end{itemize}\end{proof}

As a consequence of Theorem \ref{inject10}, we have the following
Poincaré-type inequality:
\begin{prop}\label{prp3}
 Let the hypotheses $\left(\mathrm{H}_1\right)$, $\left(\mathrm{H}_2'\right)$, $\left(\mathrm{H}_3\right)$ and $\left(\mathrm{H}_4\right)$ be satisfied. Then, there exists a constant $C>0$ independent of $u$ such that
\begin{equation}\label{af}
\|u\|_{\lh} \leq C\|\nabla u\|_{\lh}, \quad \text{ for all } u \in
W_0^{1, \mathcal{H}}(\Omega).
\end{equation}

\end{prop}
\begin{proof}
Arguing by contradiction, suppose that the inequality \eqref{af} is
false. Then, there exists a sequence $\left\{v_n\right\}_{n \in
\mathbb{N}} \subseteq \w$ such that
$$
\left\|v_n\right\|_{\lh}>n\left\|\nabla v_n\right\|_{\lh}.
$$
Let $y_n:=\dfrac{v_n}{\|v_n\|_{\lh}}, \text { for all } n >1. $ It
is easy to check that
$$
1>\frac{1}{n}>\|\nabla y_n\|_{\lh}  \text { and } \|y_n\|_{\lh}=1, \
\text { for all } n >1 .
$$
Namely, the sequence $\left\{y_n\right\}_{n \in \mathbb{N}}$ is
bounded in $\w$. Then, we can establish the existence of a
subsequence $\left\{y_n\right\}_{n \in \mathbb{N}}$ and $y$ in  $\w$
such that
\begin{equation}\label{SO1}
y_n \rightharpoonup y \quad \text { in } \w.
\end{equation}
Through the property of weak lower semicontinuity of the map $v
\mapsto\|\nabla v\|_{\lh}$ on $\w$ (it is a convex  continuous
mapping), there holds
$$
\|\nabla y\|_{\lh} \leq \liminf _{n \rightarrow \infty}\left\|\nabla
y_n\right\|_{\lh} \leq \lim _{n \rightarrow \infty} \frac{1}{n}=0.
$$
Thus, $y\equiv c \in \mathbb{R}$ is a constant function. Therefore,
since $y\in W^{1,p^{-}}(\Omega)$, we deduce that $y \equiv 0 $. On
the other hand, in view of Theorem \ref{inject10} and \eqref{SO1},
we have
$$
y_n \longrightarrow y \quad \text { in } L^{\h}(\Omega), \text{ as }
n \longrightarrow +\infty.
$$
Hence, $\|y\|_{\lh}=\ds\lim _{n \rightarrow
\infty}\left\|y_n\right\|_{\lh}=1$, so $y \not\equiv 0$. Thus a
contradiction holds. This ends the proof.
\end{proof}
\begin{remark}

Based on Proposition \ref{prp3}, we can equip the space $\w$ with
the equivalent norm $\|\nabla u\|_{\h}.$
\end{remark}
Next, we discuss the density of the smooth functions in
$W^{1,\mathcal{H}}(\Omega)$.
\begin{theorem}\label{dnsty}
  Let hypotheses $\left(\mathrm{H}_1\right)$, $\left(\mathrm{H}_2'\right)$, $\left(\mathrm{H}_3\right)$ and $\left(\mathrm{H}_4\right)$ be satisfied. Then $C^\infty \cap W^{1, \mathcal{H}}(\O)$ is dense in $W^{1, \mathcal{H}}(\O).$
\end{theorem}
\begin{proof}
The proof of this theorem follows from Proposition \ref{AAA} and
\cite[Theorem 6.4.7]{Harjulehto2019}.
\end{proof}

\subsection{ Some properties of the new double phase operator}\label{3.2}
In this part, we present the new double phase operator, denoted by
$A$, associated to problem \eqref{t11}. For this purpose, we define
$A: W_0^{1, \mathcal{H}}(\Omega) \rightarrow \left(W_0^{1,
\mathcal{H}}(\Omega)\right)^*$ as follows:
\[
\begin{aligned}
\langle A(u), v\rangle&=\int_{\Omega}\left(|\nabla u|^{p(x, |\nabla
u|)-2} \nabla u \cdot \nabla v+\mu (x) |\nabla u|^{q(x, |\nabla
u|)-2} \nabla u \cdot \nabla v\right) \mathrm{d} x\\&=\int_{\Omega}
a(x,\vert \nabla u\vert)\nabla u\cdot\nabla v\ \mathrm{d} x,
\end{aligned}\]
for all $u, v \in W_0^{1, \mathcal{H}}(\Omega)$, where $\langle\cdot, \cdot\rangle$ signifies the duality pairing between $W_0^{1, \mathcal{H}}(\Omega)$ and its dual space $\left(W_0^{1, \mathcal{H}}(\Omega)\right)^*$.\\

The energy functional $I: W_0^{1, \mathcal{H}}(\Omega) \rightarrow
\mathbb{R}$ related to $A$ is defined by
\[
I(u):=\int_{\Omega}\mathcal{H}(x, |\nabla u|) \mathrm{d} x,\
\text{for all}\ u \in W_0^{1, \mathcal{H}}(\Omega).
\]
\begin{prop}\label{op1} Let hypotheses $\left(\mathrm{H}_1\right)$, $\left(\mathrm{H}_2'\right)$, $\left(\mathrm{H}_3\right)$ and $\left(\mathrm{H}_4\right)$  hold true, then the functional $I$ is well-defined and belongs to  $C^1\left(W_0^{1, \mathcal{H}}(\Omega),\mathbb{R}\right)$, with $I^{\prime}(u)=A(u)$.\end{prop}
\begin{proof}
First, let us observe that by Proposition \ref{rel}, we have
$$\min\left\lbrace \Vert u\Vert^{p^-}_{\w},\Vert u\Vert_{\w}^{q^+ }\right\rbrace\leq I(u)\leq \max\left\lbrace \Vert u\Vert^{p^-}_{\w},\Vert u\Vert_{\w}^{q^+ }\right\rbrace,\ $$ $\text{ for all }\ u\in W^{1,\mathcal{H}}_0(\Omega).$
Thus, the functional $I$ is well-defined.

Next, let $u,v\in W^{1,\mathcal{H}}_0(\Omega)$ and $t\in[-1,1]$.
Note that
\begin{equation}\label{C1}
   \frac{I(u+tv)-I(u)}{t}=\int_{\Omega}\frac{\mathcal{H}(x, |\nabla u+t\nabla v|)-\mathcal{H}(x, |\nabla u|)}{t}\ \mathrm{d} x.
\end{equation}
By the mean value theorem, for some $t\in \mathbb{R}$, there exists
$\theta(x, t)\in [0, 1]$, such that
\begin{align*}
    \mathcal{H}(x, |\nabla u+t\nabla v|)-\mathcal{H}(x, |\nabla u|)&=a(x, |\nabla u+\theta(x,t)t\nabla v|)(\nabla u+\theta(x,t)t\nabla v)\cdot t\nabla v.
\end{align*}
Therefore, for a.a. $x \in \O$, one has
\begin{align}\label{C2}
   \frac{\mathcal{H}(x, |\nabla u+t\nabla v|)-\mathcal{H}(x, |\nabla u|)}{t}&=a(x, |\nabla u+\theta(x,t)t\nabla v|)(\nabla u+\theta(x,t)t\nabla v)\cdot \nabla v\nonumber\\
    &\stackrel{t\to 0}{\longrightarrow }a(x, |\nabla u|)\nabla u\cdot \nabla v.
\end{align}
Using the generalized Young inequality \eqref{Yi} and  \eqref{L1},
for $0<\vert t\vert \leq t_0$, we obtain
\begin{align*}
   \frac{\mathcal{H}(x, |\nabla u+t\nabla v|)-\mathcal{H}(x, |\nabla u|)}{t}&=a(x, |\nabla u+\theta(x,t)t\nabla v|)(\nabla u+\theta(x,t)t\nabla v)\cdot \nabla v\\
  & \leq a(x, |\nabla u+\theta(x,t)t\nabla v|)\vert\nabla u+\theta(x,t)t\nabla v\vert\vert\nabla v\vert\\
   & \leq \widetilde{\mathcal{H}}(x,h(x,|\nabla u+\theta(x,t)t\nabla v|))+\mathcal{H}(x,\vert\nabla v\vert)\\
   & \leq \widetilde{\mathcal{H}}(x,h(x,|\nabla u\vert+t_0\vert \nabla v|))+\mathcal{H}(x,\vert\nabla v\vert)\\
   & \leq (q^+-1)\mathcal{H}(x,|\nabla u\vert+t_0\vert \nabla v|)+\mathcal{H}(x,\vert\nabla v\vert)\\
    & \leq \frac{(q^+-1)}{2}\left(\mathcal{H}(x,2|\nabla u\vert)+\mathcal{H}(x,2t_0\vert \nabla v|)\right)+\mathcal{H}(x,\vert\nabla v\vert).
\end{align*}
It follows, by \eqref{C1}, \eqref{C2} and the dominated convergence
theorem, that
\begin{align*}
\lim_{t\to 0} \frac{I(u+tv)-I(u)}{t}=\int_{\Omega}a(x, |\nabla
u|)\nabla u\cdot\nabla v\ \mathrm{d} x=\langle A(u),v\rangle.
\end{align*}
Thus, $I^{'}(u)=A(u)$.\\
For the $C^1$-property, let $\{u_n\}_{n \in \mathbb{N}}\subseteq \w$
be a sequence such that $u_n \to u$ in $W^{1,\mathcal{H}}_0(\Omega)$
and $v\in W^{1,\mathcal{H}}_0(\Omega)$ with $\Vert v\Vert_{\w}=1$.
We have
\begin{equation}\label{C3}
\vert \nabla u_n \vert \longrightarrow \vert \nabla u\vert \text{
in}\ L^{\mathcal{H}}(\Omega).
\end{equation}
Hence,
\begin{equation}\label{C7}
   \left\lbrace\vert \nabla u_n\vert\right\rbrace_{n\in \mathbb{N}}\ \text{converge in measure to}\  \vert\nabla u\vert\ \text{in}\ \Omega,
\end{equation}
and
\begin{equation}\label{C4}
   \left\lbrace\mathcal{H}(x,\vert \nabla u_n\vert)\right\rbrace_{n\in \mathbb{N}}\ \text{converge in measure to}\  \mathcal{H}(x,\vert\nabla u\vert)\ \text{in}\ \Omega.
\end{equation}
Then, according to the converse of Vitali’s theorem \cite[Lemma
21.6]{Bauer2001}, we can see that the sequence of functions
\begin{equation}\label{C5}
\left\lbrace\mathcal{H}(x,\vert \nabla u_n\vert)\right\rbrace_{n\in
\mathbb{N}}\ \text{ is uniformly integrable}.
\end{equation}
On the other side, by \eqref{L1} and the generalized Young
inequality\eqref{Yi}, we find that
\begin{align*}
  \int_{\Omega}a(x, |\nabla u_n|)\nabla u_n\cdot\nabla v\ \mathrm{d} x&\leq \int_{\Omega}\widetilde{\mathcal{H}}\left(x,a(x, |\nabla u_n|)|\nabla u_n|\right)+\mathcal{H}(x,\vert\nabla v\vert)\ \mathrm{d} x\\
  &\leq (q^+-1)\int_{\Omega}\mathcal{H}\left(x,|\nabla u_n|\right)\ \mathrm{d} x +\int_{\Omega}\mathcal{H}(x,\vert\nabla v\vert)\ \mathrm{d} x.
\end{align*}
It follows, by \eqref{C5}, that the sequence of functions
\begin{equation}\label{C6}
\left\lbrace a(x, |\nabla u_n|)\nabla u_n\cdot\nabla
v\right\rbrace_{n\in \mathbb{N}}\ \text{ is uniformly integrable},
\end{equation}
and, by \eqref{C7}, that
\begin{equation}\label{C8}
   \left\lbrace a(x,\vert \nabla u_n\vert)\nabla u_n\cdot \nabla v\right\rbrace_{n\in \mathbb{N}}\ \text{converge in measure to}\  a(x,\vert\nabla u\vert) \nabla u\cdot \nabla v\ \text{in}\ \Omega.
\end{equation}
Therefore, applying Vitali’s theorem \cite[Theorem
4.5.4]{Bogachev2007}, in combination with \eqref{C6} and \eqref{C8},
we obtain
$$ \langle A(u_n),v\rangle= \int_{\Omega}a(x, |\nabla u_n|)\nabla u_n\cdot\nabla v\ \mathrm{d} x  \stackrel{n\to +\infty}{\longrightarrow }\int_{\Omega}a(x, |\nabla u|)\nabla u\cdot\nabla v\ \mathrm{d} x=\langle A(u),v\rangle.$$
Thus, the proof is completed.
\end{proof}

We are now ready to discuss the main characteristics of the operator $A$ inspired by the work of Liu-Dai \cite{Liu2018} for the case of constant exponent.\\

\begin{theorem}\label{op2}~~~~~~Let hypotheses $\left(\mathrm{H}_1\right)$, $\left(\mathrm{H}_2'\right)$, $\left(\mathrm{H}_3\right)$ and $\left(\mathrm{H}_4\right)$  hold true. Then
\begin{itemize}
                \item [$(\mathrm{i})$] The operator $A: W_0^{1, \mathcal{H}}(\Omega) \rightarrow \left(W_0^{1, \mathcal{H}}(\Omega)\right)^*$ is continuous, bounded, and strictly monotone.
\item [$(\mathrm{ii})$] The operator $A: W_0^{1, \mathcal{H}}(\Omega) \rightarrow \left(W_0^{1, \mathcal{H}}(\Omega)\right)^*$ fulfills the $\left(\mathrm{S}_{+}\right)$-property, i.e.,
\[
u_n \rightharpoonup u \text { in } W_0^{1, \mathcal{H}}(\Omega)
\quad \text { and } \quad \limsup _{n \rightarrow
\infty}\left\langle A\left(u_n\right), u_n-u\right\rangle \leq 0,
\]
leads to $u_n \rightarrow u$ in $W_0^{1, \mathcal{H}}(\Omega)$.
\item [$(\mathrm{iii})$] The operator $A: W_0^{1, \mathcal{H}}(\Omega) \rightarrow \left(W_0^{1, \mathcal{H}}(\Omega)\right)^*$  forms a homeomorphism.
\item [$(\mathrm{iv})$] The operator $A: W_0^{1, \mathcal{H}}(\Omega) \rightarrow \left(W_0^{1, \mathcal{H}}(\Omega)\right)^*$ is strongly coercive, that is,
    $$\lim_{\|u\|_{\w}\rightarrow+\infty} \dfrac{\langle Au, u\rangle}{\|u\|_{\w}}\longrightarrow+\infty.$$
    \end{itemize}
\end{theorem}
\begin{proof}
\begin{itemize}
    \item [$(\mathrm{i})$]  According to Proposition \ref{op1}, we have that $A=I^{\prime}$ with $I$ being of class $C^1$, hence $A$ is continuous.

    From  Lemma \ref{simon} and Proposition \ref{rel}, we obtain
    \begin{align*}
        \left\langle A\left(u\right)-A(v), u-v\right\rangle&=\int_{\Omega}\left(|\nabla u|^{p(x, |\nabla u|)-2} \nabla u- |\nabla v|^{p(x, |\nabla v|)-2} \nabla v\right)\cdot \nabla (u-v) \ \mathrm{d} x\\
        & + \int_{\Omega}\mu (x) \left(|\nabla u|^{q(x, |\nabla u|)-2} \nabla u- |\nabla v|^{q(x, |\nabla v|)-2} \nabla v\right)\cdot \nabla (u-v) \ \mathrm{d} x\\
        &= \int_{\Omega} \left( a(x,|\nabla u|)\nabla u-a(x,|\nabla v|)\nabla v\right)\cdot\nabla (u-v) \ \mathrm{d} x\\
        & \geq 4\int_{\Omega} \mathcal{H}\left(x,\frac{\vert \nabla (u-v)\vert}{2} \right)\ \mathrm{d} x\\& \geq 4\min\left\lbrace \left\Vert \frac{u-v}{2}\right\Vert^{p^-}_{\w},\left\Vert \frac{u-v}{2}\right\Vert_{\w}^{q^+}\right\rbrace\\
        &>0,
    \end{align*}
   for all $u\not\equiv v$. Thus, the operator $A$ is strictly monotone.\\
   To prove the boundedness of $A$, let us take $u,v\in W^{1,\mathcal{H}}_0(\Omega)\setminus\lbrace 0\rbrace$. In light of \eqref{L1} and Young’s inequality \eqref{Yi}, we see that
    \begin{align*}
    \left\langle A(u), \frac{v}{\left\Vert v\right\Vert_{\w}}\right\rangle&=\int_{\Omega}  a(x,|\nabla u|) \nabla u\cdot\frac{\nabla v}{\left\Vert v\right\Vert_{\w}}\ \mathrm{d} x\\
    & \leq \int_{\Omega} \left[ \widetilde{\mathcal{H}}\left(x,a(x,|\nabla u|) |\nabla u|\right) + \mathcal{H}\left(x,\frac{\vert \nabla v\vert}{\left\Vert v\right\Vert_{\w}}\right)\right]\ \mathrm{d} x\\
    & \leq \int_{\Omega} \left[ (q^+-1)\mathcal{H}\left(x,|\nabla u|\right) + \mathcal{H}\left(x,\frac{\vert \nabla v\vert}{\left\Vert v\right\Vert_{\w}}\right)\right]\ \mathrm{d} x\\
     & \leq (q^+-1)\int_{\Omega} \mathcal{H}\left(x,|\nabla u|\right)\ \mathrm{d} x+\int_{\Omega} \mathcal{H}\left(x,\frac{\vert \nabla v\vert}{\left\Vert v\right\Vert_{\w}}\right)\ \mathrm{d} x\\
   &\leq (q^+-1)I(u)+1.
   \end{align*}
   Therefore, by Proposition \ref{rel}, it yields that
  $$\begin{aligned}\Vert A(u)\Vert_*=\sup_{v\in W^{1,\mathcal{H}}_0(\Omega)\setminus\lbrace 0\rbrace}\frac{\langle A(u),v\rangle}{\Vert v\Vert_{\w}}&\leq (q^+-1)I(u)+1\\&\leq (q^+-1)\max\left\lbrace \Vert u\Vert^{p^-}_{\w},\Vert u\Vert_{\w}^{q^+}\right\rbrace+1.\end{aligned}$$
  Hence, $A$ is bounded.
  \item [$(\mathrm{ii})$] Let $\lbrace u_n\rbrace_{n\in \mathbb{N}}\subseteq W^{1,\mathcal{H}}_0(\Omega)$ be a sequence such that
  \begin{equation}\label{S+}
     u_n \rightharpoonup u \text { in } W_0^{1, \mathcal{H}}(\Omega) \quad \text { and } \quad \limsup _{n \rightarrow \infty}\left\langle A\left(u_n\right), u_n-u\right\rangle \leq 0.
  \end{equation}
Exploiting the weak convergence of $u_n$ to $u$
 in $W_0^{1, \mathcal{H}}(\Omega)$, we get
 $$\lim _{n \rightarrow \infty}\left\langle A\left(u\right), u_n-u\right\rangle = 0.$$
 Thus, by \eqref{S+}, we prove that
 $$ \limsup _{n \rightarrow \infty}\left\langle A\left(u_n\right)-A(u), u_n-u\right\rangle \leq 0.$$
 Then, the strict monotonicity of the operator $A$ gives that
  $$ 0\leq \liminf _{n \rightarrow \infty}\left\langle A\left(u_n\right)-A(u), u_n-u\right\rangle\leq \limsup _{n \rightarrow \infty}\left\langle A\left(u_n\right)-A(u), u_n-u\right\rangle \leq 0.$$
  Hence, we see that
  \begin{equation}\label{S1}
      \lim_{n \rightarrow \infty}\left\langle A\left(u_n\right)-A(u), u_n-u\right\rangle=0=\lim _{n \rightarrow \infty}\left\langle A\left(u_n\right), u_n-u\right\rangle.
  \end{equation}
  On the one side, by \eqref{S1} and Lemma \ref{simon}, we find that
\begin{equation}\label{S2}
\lim_{n\to +\infty}\int_{\Omega} \mathcal{H}\left(x,\frac{\vert
\nabla u_n-\nabla u\vert}{2} \right)\ \mathrm{d} x=0.
\end{equation}
Then, from Proposition \ref{rel}, it yields that
\begin{equation}\label{S3}
 \vert\nabla u_n\vert \longrightarrow \vert\nabla u\vert \text{ in}\ L^{\mathcal{H}}(\Omega).
\end{equation}
Thus, $\lbrace \vert \nabla u_n\vert \rbrace_{n\in \mathbb{N}}$
converges in measure to $\vert\nabla u\vert$ in $\Omega$. Therefore,
there exists
a subsequence, still denoted by $\lbrace \vert\nabla u_n \vert\rbrace_{n\in \mathbb{N}}$, that converges to $\vert\nabla u \vert$ a.a. in $\Omega$.\\
Since $t\longmapsto \mathcal{H}(x,t)$ is convex for all $x\in
\Omega$, then the functional $I$ is convex. Hence,
\begin{equation*}
    I(u)\geq I(u_n)+\langle A(u_n),u-u_n\rangle.
\end{equation*}
This is equivalent to
\begin{equation}\label{S4}
  \langle A(u_n),u_n-u\rangle \geq I(u_n)-I(u).
\end{equation}
Combining \eqref{S1} and \eqref{S4}, we obtain
\begin{equation}\label{S5}
   \lim_{n\to +\infty}I(u_n)\leq I(u).
\end{equation}
On the other side, by Fatou's lemma, we have
\begin{equation}\label{S6}
  I(u) \leq\liminf_{n\to +\infty}I(u_n).
\end{equation}
From \eqref{S5} and \eqref{S6}, it yields that
\begin{equation}\label{S7}
\lim_{n\to +\infty}I(u_n)= I(u),
\end{equation}
which implies,  together with the fact that $\lbrace \vert\nabla u_n
\vert\rbrace_{n\in \mathbb{N}}$ converges in measure to $\vert\nabla
u \vert$ in $\Omega$, that
\begin{equation}\label{S8}
   \left\lbrace\mathcal{H}(x,\vert \nabla u_n\vert)\right\rbrace_{n\in \mathbb{N}}\ \text{converge in measure to}\  \mathcal{H}(x,\vert\nabla u\vert)\ \text{in}\ \Omega.
\end{equation}
According to the converse of Vitali’s theorem \cite[Lemma
21.6]{Bauer2001}, we can see that the sequence of functions
\begin{equation}\label{S9}
\left\lbrace\mathcal{H}(x,\vert \nabla u_n\vert)\right\rbrace_{n\in
\mathbb{N}}\ \text{ is uniformly integrable}.
\end{equation}
Using the monotonicity and  convexity of the function $t\longmapsto
\mathcal{H}(\cdot,t)$  and Proposition \ref{zoo}, for all $x\in
\Omega$, we get
\begin{align*}
    \mathcal{H}(x,\vert \nabla u_n-\nabla u\vert)&\leq \mathcal{H}(x,\vert \nabla u_n\vert+\vert\nabla u\vert)\\
    & \leq \frac{1}{2}\left( \mathcal{H}(x,2\vert\nabla u_n\vert)+ \mathcal{H}(x,2\vert \nabla u\vert)\right)\\
    & \leq 2^{q^+-1}\left(\mathcal{H}(x,\vert\nabla u_n\vert)+ \mathcal{H}(x,\vert \nabla u\vert)\right).
\end{align*}
Thus, from \eqref{S9}, we have the uniform integrability of the
sequence of functions
$$\left\lbrace\mathcal{H}(x,\vert \nabla u_n-\nabla u\vert)\right\rbrace_{n\in \mathbb{N}}.$$
Again, since  $\lbrace \vert \nabla u_n \vert\rbrace_{n\in
\mathbb{N}}$ converges in measure to $\vert\nabla u \vert$ in
$\Omega$, we see that
\begin{equation}\label{S10}
   \left\lbrace\mathcal{H}(x,\vert \nabla u_n-\nabla u\vert)\right\rbrace_{n\in \mathbb{N}}\ \text{converge in measure to}\  0\ \text{in}\ \Omega.
\end{equation}
Applying Vitali’s theorem \cite[Theorem 4.5.4]{Bogachev2007}, it
follows that
$$
\lim_{n\to +\infty}I(u_n-u)=\int_{\Omega}\mathcal{H}(x,\vert \nabla
u_n-\nabla u\vert)\ \mathrm{d} x=0.
$$
Finally, by Proposition \ref{rel}, we obtain
$$\Vert u_n-u\Vert_{\w}\longrightarrow 0\ , \text{ as }\ n\longrightarrow +\infty.$$
Namely, $u_n\longrightarrow u$ in $W^{1,\mathcal{H}}_0(\Omega)$.
Hence, the operator $A$ satisfies the $(S_+)$-property.
\item [$(\mathrm{iii})$] 
 First, by strict monotonicity, $A$ is an injection. Furthermore, by the Minty-Browder theorem, we find that $A$ is a surjection. Hence, $A$ has an inverse mapping $A^{-1}:\left(W^{1,\mathcal{H}}_0(\Omega)\right)^*\longrightarrow W^{1,\mathcal{H}}_0(\Omega)$. Therefore,
to complete the proof of $(\mathrm{iii})$, it remains to show that
$A^{-1}$ is continuous. To this end, let $\lbrace f_n\rbrace_{n\in
\mathbb{N}}\subseteq \left(W^{1,\mathcal{H}}_0(\Omega)\right)^*$ be
a sequence such that $f_n\longrightarrow f$ in
$\left(W^{1,\mathcal{H}}_0(\Omega)\right)^*$  and let
$u_n:=A^{-1}(f_n)$ as well as $u:=A^{-1}(f)$. By the strong
convergence of $\lbrace f_n\rbrace_{n\in \mathbb{N}}$ and the
boundedness of $A^{-1}$ we get that the sequence $u_n$ is bounded in
$W^{1,\mathcal{H}}_0(\Omega)$. Hence, without loss of generality, we
can assume that there exists a subsequence, still denoted by
$\lbrace u_n\rbrace_{n\in \mathbb{N}}$,
 such that
 $$u_n\rightharpoonup \widetilde{u}\ \text{in}\ W^{1,\mathcal{H}}_0(\Omega).$$
 It follows that
 \begin{align*}
  \lim _{n \rightarrow \infty}\left\langle A\left(u_n\right)-A(\widetilde{u}), u_n-\widetilde{u}\right\rangle  &= \lim _{n \rightarrow \infty}\left\langle f_n-f, u_n-\widetilde{u}\right\rangle =0
 \end{align*}
 which proves, together with the fact that $A$ satisfies the $S_+$-property, that
 $$u_n\longrightarrow \widetilde{u}\ \text{in}\ W^{1,\mathcal{H}}_0(\Omega).$$
  By the continuity of the operator $A$ we find that
  $$A(\widetilde{u})=\lim_{n\to +\infty}A(u_n)=\lim_{ n\to +\infty}f_n=f=A(u).$$
  Since $A$ is injective, it follows that $u \equiv \widetilde{u}$. Hence,
  $$u_n\longrightarrow u\ \text{in}\ W^{1,\mathcal{H}}_0(\Omega).$$
  Therefore, $A^{-1}$ is continuous, so the operator $A$ is a homeomorphism.\\
  \item [$(\mathrm{iv})$] First let us observe that, by Proposition \ref{rel}, one has
\begin{align*}
    \frac{\langle A(u),u\rangle}{\Vert u\Vert_{\w}}&=\frac{1}{\Vert u\Vert_{\w}}\int_{\Omega}\left(|\nabla u|^{p(x, |\nabla u|)} +\mu (x) |\nabla u|^{q(x, |\nabla u|)} \right) \mathrm{d} x\\
    & \geq \frac{1}{\Vert u\Vert_{\w}}\int_{\Omega}\mathcal{H}(x,\vert \nabla u\vert) \mathrm{d} x\\
    & \geq \frac{1}{\Vert u\Vert_{\w}}\min\left\lbrace \Vert u\Vert_{\w}^{p^-}, \Vert u\Vert_{\w}^{q^+}\right\rbrace\\
    & =\min\left\lbrace \Vert u\Vert_{\w}^{p^--1}, \Vert u\Vert_{\w}^{q^+-1}\right\rbrace.
\end{align*}
Hence, since $q^+ >p^->1$, the operator $A$ is strongly coercive.
 \end{itemize}
This ends the proof.
\end{proof}

\subsection{Applications}
Now, we are willing to tackle the problem \eqref{t11}. This
subsection is divided into two parts: the first part deals with the
variational case, and the second part focuses on the nonvariational
situation.
\subsubsection{Variational problem driven by the new double phase operator}
Let $\Omega$ be a bounded open subset of $\mathbb{R}^d$ $(d\geq 2)$
with a Lipschitz boundary $\partial\Omega$. We consider the
following Dirichlet problem:
\begin{equation}\label{P}
\left\lbrace\begin{array}{rll}
-\operatorname{div}\left(|\nabla u|^{p(x,|\nabla u|)-2} \nabla u+ \mu(x)|\nabla u|^{q(x,|\nabla u|)-2} \nabla u\right) &=f(x, u),  & \text { in } \Omega, \\
\ & \ & \\
u &=0, & \text { on } \partial \Omega,
\end{array}\right.\tag{$\mathcal{P}_1$}
\end{equation}
where $\mu \in L^\infty (\Omega)$ such that $\mu(x)\geq 0$ a.a. $x
\in \Omega$, $p, q, : \Omega \times \mathbb{R} \longrightarrow
\mathbb{R}$ are tow Carathéodory functions satisfying assumptions
$\left(\mathrm{H}_1\right)$, $\left(\mathrm{H}_2'\right)$,
$\left(\mathrm{H}_3\right)$ and $\left(\mathrm{H}_4\right)$. For the
reaction term $f$, we consider the following class of conditions:
\begin{enumerate}
    \item[$(\mathcal{F}_1)$]  $f: \Omega \times \mathbb{R} \longrightarrow \mathbb{R}$ is a Carathéodory function satisfying
    \begin{enumerate}
\item[$(f_1)$] $f(x,0)=0$ and there exist  a generalized $N$-function $\displaystyle{\mathcal{B}(x,t)=\int_{0}^t b(x,s) \ ds}$ with $x\mapsto b(x,\cdot)$ continuous on $\overline{\Omega}$, and a positive function $\widehat{a}(x)\in L^{\infty}(\Omega)$ such that
$$\vert f(x,t)\vert \leq \widehat{a}(x)(1+b(x,\vert t\vert)), \ \ \forall\ t\in \mathbb{R},\ \text{a.a.}\ x\in \overline{\Omega},$$
$$\mathcal{H}\ll \mathcal{B}\ll \mathcal{H}_*,$$
and
$$\displaystyle{p^-<b^-:=\inf\limits_{t>0}\frac{b(x,t)t}{\mathcal{B}(x,t)}\leq b^{+}:=\sup\limits_{t>0}\frac{b(x,t)t}{\mathcal{B}(x,t)}<+\infty,\ \text{a.a.}\ x\in \overline{\Omega}.}$$
\item[$(f_2)$] There exist a generalized N-function $\displaystyle{\mathcal{S}(x,t)=\int_{0}^t s(x,r) \ dr}$, and a positive constants $\widetilde{C}\geq 0$, $\delta\geq 0$ such that
$$\mathcal{S}\ll \mathcal{H},$$
$$\displaystyle{1<s^-:=\inf\limits_{t>0}\frac{s(x,t)t}{\mathcal{S}(x,t)}\leq s^{+}:=\sup\limits_{t>0}\frac{s(x,t)t}{\mathcal{S}(x,t)}<p^-},\ \text{a.a.}\ x\in \Omega$$
and
$$\widetilde{C} s(x,t)t \leq s^+F(x,t),\ \text{a.a.}\ x\in \Omega\ \text{and for all}\ 0<\vert t\vert \leq\delta.$$
\item[$(f_3)$] There exist $\eta_-<0$ and  $\eta_+>0$ such that
$$f(x,\eta_+)<0<f(x,\eta_-),\ \ \text{for a.a.}\ x\in \Omega.$$
\end{enumerate}
\end{enumerate}
\begin{dfn}\label{weak1}
    We say that $u\in W^{1,\mathcal{H}}_0(\Omega)$ is a weak solution of problem \eqref{P}, if
    $$\langle A(u),v\rangle=\int_{\Omega}f(x,u)v\ dx,\ \ \text{for all}\ v\in W^{1,\mathcal{H}}_0(\Omega).$$
\end{dfn}

The problem \eqref{P} has a variational structure and its
corresponding  energy functional is
$J:W^{1,\mathcal{H}}_0(\Omega)\longrightarrow \mathbb{R}$ defined by
\begin{align}\label{A1}
J(u):=I(u) -\int_{\Omega}F(x,u)dx,\ \ \text{for all}\ u\in
W^{1,\mathcal{H}}_0(\Omega).
\end{align}
According to Proposition \ref{op1}, we see that $J\in
C^{1}(W^{1,\mathcal{H}}_0(\Omega),\mathbb{R})$ and
\begin{equation}\label{A2}
    \langle J^{'}(u),v\rangle=\langle A(u),v\rangle-\int_{\Omega}f(x,u)v\ dx,\ \ \text{for all}\ u,v\in W^{1,\mathcal{H}}_0(\Omega).
\end{equation}
Note that, all weak solutions of problem \eqref{P} are critical
points of the functional $J$. In order to prove the existence of
nonnegative and nonpositive weak solutions to problem \eqref{P},
 let us introduce the following two Carath\'eodory functions $f_+:\Omega\times\mathbb{R}\longrightarrow \mathbb{R}$ defined by
\begin{equation}\label{A3}
f_+(x,t):=\left\lbrace
\begin{array}{ll}
 f(x,t^+) &\text{if}\ \ t\leq\eta_+\\
\ & \ \\
 f(x,\eta_+) & \text{if}\ \ t> \eta_+,
\end{array}
 \right.
\end{equation}
and $f_-:\Omega\times\mathbb{R}\longrightarrow\mathbb{R}$ defined by
\begin{equation}\label{A4}
f_-(x,t):=\left\lbrace
\begin{array}{ll}
 f(x,\eta_-) &\text{if}\ \ t<\eta_-\\
\ & \ \\
 f(x,t^-) & \text{if}\ \ t\geq \eta_-,
\end{array}
 \right.
\end{equation}
where $\eta_+$ and $\eta_-$ are mentioned in $(f_4)$. We set
$\displaystyle{F_\pm (x,s):=\int^{s}_{0}f_\pm (x,t)\ dt}$ and
consider the $C^1$-functionals
$J_\pm:W^{1,\mathcal{H}}_0(\Omega)\longrightarrow \mathbb{R}$
defined by
\begin{align}\label{A5}
J_\pm(u):=I(u) -\int_{\Omega}F_\pm(x,u)dx,\ \ \text{for all}\ u\in
W^{1,\mathcal{H}}_0(\Omega).
\end{align}
It is clear that $J_\pm\in
C^{1}(W^{1,\mathcal{H}}_0(\Omega),\mathbb{R})$ and
\begin{equation}\label{A6}
    \langle J_\pm^{'}(u),v\rangle=\langle A(u),v\rangle-\int_{\Omega}f_\pm(x,u)v\ dx,\ \ \text{for all}\ u,v\in W^{1,\mathcal{H}}_0(\Omega).
\end{equation}
The main result of this part is the following:
\begin{theorem}\label{thmA1}
    Let hypotheses $(\mathcal{F}_1)$ , $\left(\mathrm{H}_1\right)$, $\left(\mathrm{H}_2'\right)$, $\left(\mathrm{H}_3\right)$ and $\left(\mathrm{H}_4\right)$ be satisfied. Then, problem \eqref{P} admits at least one nonnegative and one nonpositive  weak solutions $u_0,v_0\in W^{1,\mathcal{H}}_0(\Omega)$.
\end{theorem}
\begin{proof} We start with the existence of a nonnegative solution.
Using \eqref{A3}, \eqref{A5} and Proposition \ref{rel}, we find that
\begin{align}\label{A7}
J_{+}(u)\geq \min\left\lbrace \Vert u\Vert^{p^-},\Vert u\Vert^{q^+}
\right\rbrace-C_1\vert \Omega\vert,\ \ \text{for all}\ u\in
W^{1,\mathcal{H}}_0(\Omega)
\end{align}
where $C_1$ is a positive constant. From \eqref{A7})  it is clear that $J_+$ is coercive.\\
Let $\lbrace u_n\rbrace_{n\in \mathbb{N}}\subseteq
W^{1,\mathcal{H}}_0(\Omega)$ such that $u_n\rightharpoonup u$ in
$W^{1,\mathcal{H}}_0(\Omega)$. We have
\begin{align}\label{A8}
\liminf\limits_{n\rightarrow +\infty}J_{+}(u_n)\geq
\liminf\limits_{n\rightarrow
+\infty}I(u_n)-\limsup\limits_{n\rightarrow
+\infty}\int_{\Omega}F_{+}(x,u_n)\ dx.
\end{align}
 Exploiting Fatou's lemma, Theorem \ref{inject10} and \eqref{A8}, we obtain
 $$\liminf\limits_{n\rightarrow +\infty}J_{+}(u_n)\geq J_{+}(u).$$
Therefore, $J_+$ is sequentially weakly lower semicontinuous. Then,
by the Weierstrass-Tonelli theorem we can find  $u_0\in
W^{1,\mathcal{H}}_0(\Omega)$ such that
\begin{equation}\label{A9}
J_{+}(u_0)=\min_{v\in W^{1,\mathcal{H}}_0(\Omega)}J_{+}(v).
\end{equation}
Let $u\in \text{int}(C^1(\overline{\Omega})_+)$, where
$\text{int}(C^1(\overline{\Omega})_+):=\left\{u\in
C^1(\overline{\O});\ u(x)>0 \text{ for all }x \in
\overline{\O}\right\}$, and choose $t\in (0,1)$ small enough such
that
$$0< tu(x)\leq \min\lbrace \delta,\eta_+\rbrace,\ \ \text{for all}\ \ x\in \overline{\Omega}.$$
Using Proposition \ref{rel} and assumption $(f_2)$, we get
\begin{align*}
J_{+}( tu)& =I(tu)-\int_{\Omega}F_{+}(x, tu)\ dx\\
& \leq \vert t\vert^{p^-}I(u)-\int_{\Omega}F(x, tu)\ dx\\
& \leq \vert t\vert^{p^-}I(u)-\widetilde{C}\frac{s^-}{s^+}\vert
t\vert^{s^+}\int_{\Omega}\mathcal{S}(x,\vert u\vert)\ dx.
\end{align*}
Since $s^+<p^-$, we can choose $t\in(0,1)$ sufficiently small such
that $J_{+}(  tu)<0$. Hence, by \eqref{A9}, we get $$J_+(u_0)\leq
J_+(tu)<0=J_+(0).$$ Therefore, $u_0\not\equiv 0$. Recall that $u_0$
is a global minimizer  of  $J_+$, then
\begin{align}\label{A10}
\langle J_+^{'}(u_0),v\rangle=\langle
A(u_0),v\rangle-\int_{\Omega}f_+(x,u_0)v\ dx=0,\ \ \text{for all}\
v\in \w.
\end{align}
In \eqref{A10}, we act with $v\equiv u_0^-\equiv \max(-u,0)$ , we
obtain
\begin{align*}
 \int_{\Omega} a(x,\vert\nabla u_0\vert)\nabla u_0.\nabla u_0^- dx & =\int_{\Omega}f_+(x,u_0)u_0^-dx,
\end{align*}
it follows, by \eqref{l2} and the truncation \eqref{A3}, that
\begin{align*}
p^-I(u_0^-)\leq \int_{\Omega} a(x,\vert\nabla u_0^-\vert)\vert\nabla
u_0^-\vert^2 dx \leq \int_{\Omega} a(x,\vert\nabla u_0\vert)\nabla
u_0\nabla u_0^- dx\leq0.
\end{align*}
Thus, $u_0^- \equiv 0$. Hence,  $u_0\not\equiv 0$ and $u_0\geq 0$.
Again, in \eqref{A10}, we act with $u\equiv(u_0-\eta_+)^+$ and using
$(f_3)$, we infer that
\begin{align}\label{A11}
\int_{\Omega}a(x,\vert\nabla u_0\vert)\nabla u_0.\nabla (u_0-\eta_+)^+ dx & = \int_{\Omega}f_+(x,u_0)(u_0-\eta_+)^+dx\nonumber\\
&=\int_{\Omega}f(x,\eta_+)(u_0-\eta_+)^+dx\nonumber\\
&\leq 0.
\end{align}
Exploiting \eqref{A11}, Lemma \ref{simon}, and Proposition
\ref{rel}, we get
\begin{align*}
0& \geq\int_{u_0\geq \eta_+}a(x,\vert\nabla u_0\vert)\nabla u_0.\nabla (u_0-\eta_+) dx \nonumber\\
& \geq 4\int_{u_0\geq \eta_+}\mathcal{H}\left(x, \frac{\vert \nabla( u_0-\eta_+)\vert}{2}\right) \nonumber\\
& =  4\int_{\Omega}\mathcal{H}\left(x, \frac{\vert \nabla (u_0-\eta_+)^+\vert}{2}\right)\nonumber\\
& \geq 4\min\left\lbrace \left\Vert
\frac{(u_0-\eta_+)^+}{2}\right\Vert_{\w} ^{p^-},\left\Vert
\frac{(u_0-\eta_+)^+}{2}\right\Vert_{\w} ^{q^+} \right\rbrace.
\end{align*}
Therefore, we infer that $(u_0-\eta_+)^+\equiv0$. Namely,
$u_0\in[0,\eta_+]$. Hence, from the truncation \eqref{A3}, we
conclude that $u_0$ is a nonnegative bounded weak solution for
problem \eqref{P}. Similarly, using the functional $J_-$ and the
truncation \eqref{A4}, we prove that problem \eqref{P} has a
nonpositive weak solution $v_0\in [\eta_-,0]$. Thus, the proof is
completed.
\end{proof}
\subsubsection{Nonvariational problem driven by the new double phase operator}
Let $\Omega$ be a bounded domain in  $\mathbb{R}^d$ $(d\geq 2)$
with a Lipschitz boundary $\partial\Omega$. We consider the
following Dirichlet problem:
\begin{equation}\label{P2}
\left\lbrace\begin{array}{rll}
-\operatorname{div}\left(|\nabla u|^{p(x,|\nabla u|)-2} \nabla u+ \mu(x)|\nabla u|^{q(x,|\nabla u|)-2} \nabla u\right) &=\vert u\vert^{r(x,\vert \nabla u\vert)-2}u +f(x),  & \text { in } \Omega, \\
&\ & \ \\
u &=0, & \text { on } \partial \Omega,
\end{array}\tag{$\mathcal{P}_2$} \right.
\end{equation}
where $\mu \in L^\infty (\Omega)$ such that $\mu(x)\geq 0$ a.a. $x
\in \Omega$, and $p, q, r: \Omega \times \mathbb{R} \longrightarrow
\mathbb{R}$ are three Carathéodory functions satisfying assumption
$\left(\mathrm{H}_1\right)$, $\left(\mathrm{H}_2'\right)$,
$\left(\mathrm{H}_3\right)$,  $\left(\mathrm{H}_4\right)$ and
\begin{itemize}
  \item[$(\mathcal{R})$] $\displaystyle{: 1\leq r^-:=\inf_{(x,t)\in \Omega\times \R}r(x,t)\leq r(x,\cdot)\leq r^+:=\sup_{(x,t)\in \Omega\times \R}r(x,t)\leq p^-}$.
  \item[$(\mathcal{F}_2)$] $f\in L^2(\O)$ such that 
  $f\not\equiv 0.$
\end{itemize}

\begin{dfns}\label{weak2}
\begin{itemize}
  \item [(1)] We say that $u\in W^{1,\mathcal{H}}_0(\Omega)$ is a weak solution of problem \eqref{P2}, if
    $$\langle A(u),v\rangle-\langle\mathcal{N}_r(u),v\rangle=0,\ \ \text{for all}\ v\in W^{1,\mathcal{H}}_0(\Omega),$$
     where $\displaystyle{\langle\mathcal{N}_r(u),v\rangle:=\int_{\Omega}\vert u\vert^{r(x,\vert \nabla u\vert)-2}uv\ dx}+\int_\O f(x)vdx$.
      \item [(2)] The operator $A$ is called  pseudomonotone if $u_n\rightharpoonup u$ in $\w$  and $\ds \limsup_{n \rightarrow +\infty}\langle Au_n, u_n - u\rangle \leq 0$ imply
          $$
          \ds  A(u_n) \rightharpoonup A(u) \ \mbox{and} \
          \langle A(u_n),u_n\rangle \rightarrow \langle A(u),u\rangle.
          $$
      \end{itemize}
\end{dfns}
Note that, due to the behavior of the nonlinear term, the equation
\eqref{P2} has a nonvariational structure. Therefore, our existence
result is based on the surjectivity theorem for pseudomonotone
operators given below, see Papageorgiou-Winkert \cite[Theorem
6.1.57]{papa1}.

\begin{theorem}\label{surjt}
 Let $X$ be a real, reflexive Banach space and $A$ is a pseudomonotone, bounded, and strongly coercive operator from $X$ to its dual space $X^\ast$, then there exists a solution to the equation $Au = b$ for any $b$ in $X^\ast$.
\end{theorem}
The main result of this part reads as follows:
\begin{theorem}\label{thmA2}
Let hypotheses $(\mathcal{F}_2)$,  $\left(\mathrm{H}_1\right)$,
$\left(\mathrm{H}_2'\right)$, $\left(\mathrm{H}_3\right)$,
$\left(\mathrm{H}_4\right)$ and $(\mathcal{R})$ be satisfied. Then,
problem \eqref{P2} admits at least one nontrivial weak solution
$u\in W^{1,\mathcal{H}}_0(\Omega)$.
\end{theorem}
\begin{proof}
    For $u\in W^{1,\mathcal{H}}_0(\Omega)$, we define $\mathcal{A}:W^{1,\mathcal{H}}_0(\Omega)\longrightarrow \left(W^{1,\mathcal{H}}_0(\Omega)\right)^*$ by
    $$\mathcal{A}(u):=A(u)-\mathcal{N}_r(u),\ \text{for all}\ u,v\in W^{1,\mathcal{H}}_0(\Omega).$$
    In order to prove that problem \eqref{P2} admits a nontrivial weak solution, we shall prove that the operator $\mathcal{A}$ is bounded, continuous, strong coercive, and pseudomonotone.
    To this end, using Theorem \ref{inject10} and assumptions  $(\mathcal{F}_2)$ and
 $(\mathcal{R})$, the operator $\mathcal{A}$ is well-defined.  Moreover, it is clear that, from Theorem \ref{op2}, the operator $\mathcal{A}$ is bounded and continuous.
 So, it remains to show the strong coercivity and the pseudomonotonicity proprieties of $\mathcal{A}$. For this purpose, we consider the sequence  $\lbrace u_n\rbrace_{n\geq 1}\subseteq W^{1,\mathcal{H}}_0(\Omega)$ such that
    \begin{equation}\label{A21}
        u_n\rightharpoonup u\ \text{in}\  W^{1,\mathcal{H}}_0(\Omega)\ \text{and}\ \limsup _{n \rightarrow \infty}\left\langle \mathcal{A}\left(u_n\right), u_n-u\right\rangle \leq 0.
    \end{equation}
In view of Theorem \ref{inject10}, we have
\begin{equation}\label{A22}
       \left\lbrace  \begin{aligned}u_n\longrightarrow u\ \text{in}\  L^{r^-}(\Omega),&\\
       u_n\longrightarrow u\ \text{in}\  L^{r^+}(\Omega)&\\
       \text{ and }u_n\longrightarrow u\ \text{in}\  L^{2}(\Omega).
       \end{aligned}\right.
    \end{equation} \
   Then,  using again Theorem \ref{inject10} and  H\"older inequality, we obtain
    $$\begin{aligned}
        \left\vert \langle\mathcal{N}_r(u_n),u_n-u\rangle\right\vert &=\left\vert \int_{\Omega}\vert u_n\vert^{r(x,\vert \nabla u_n\vert)-2}u_n(u_n-u)\ dx+\int_\O f(x)(u_n-u)dx\right\vert\\
        & \leq \int_\O |u_n|^{r^--1}|u_n-u|dx +\int_\O |u_n|^{r^+-1}|u_n-u|dx \\&\ \ \ + \int_{\Omega}\vert f(x) \vert \vert u_n-u\vert dx\\
        & \leq c\small \left( \|u_n\|^{r^--1}_{\w} \|u_n-u \|_{L^{r^{-}}(\Omega)}+\|u_n\|^{r^+-1}_{\w} \|u_n-u \|_{L^{r^{+}}(\Omega)}\right)\\& \ \ \  +c\|f\|_{L^2(\O)}\|u_n-u\|_{L^2(\O)} .
    \end{aligned}$$
    From this along with \eqref{A22}, we conclude that
    $$\lim_{n\to +\infty}\langle\mathcal{N}_r(u_n),u_n-u\rangle=0.$$
  It follows, by \eqref{A21}, that
    $$\limsup _{n \rightarrow \infty}\left\langle A\left(u_n\right), u_n-u\right\rangle=\limsup _{n \rightarrow \infty}\left\langle \mathcal{A}\left(u_n\right), u_n-u\right\rangle \leq 0.$$
    Since $A$ fulfills the $(\text{S}_+)$-property, (see Theorem \ref{op2}(ii)), then $u_n\longrightarrow u$ in $W^{1,\mathcal{H}}_0(\Omega)$. Therefore, by the continuity of
    the operator $\mathcal{A}$, we conclude that $\mathcal{A}(u_n)\longrightarrow \mathcal{A}(u)$ in $\left(W^{1,\mathcal{H}}_0(\Omega)\right)^*$ and
    $\langle A(u_n),u_n\rangle \rightarrow \langle A(u),u\rangle$ in $\mathbb{R}$. Thus, the operator $\mathcal{A}$ is pseudomonotone.\\
    Let us now prove that $\mathcal{A}$ is strongly coercive.  For $u\in \w$, we have
\begin{align*}
    \langle \mathcal{A}(u),u\rangle&= \langle A(u),u\rangle- \langle \mathcal{N}_r(u),u\rangle\\
    &=\int_{\Omega}\left(|\nabla u|^{p(x, |\nabla u|)} +\mu (x) |\nabla u|^{q(x, |\nabla u|)} \right) \mathrm{d} x-\int_{\Omega}\vert u\vert^{r(x,\vert \nabla u\vert)}\ dx -\int_{\O}^{} f(x)udx.
\end{align*}
Then, by $(\mathcal{R})$, Proposition \ref{rel} and for $u\in \w$
with $\|u\|_{\w} \geq1$, we have
$$\langle \mathcal{A}(u),u\rangle \geq   \|u\|_{\w}^{q^+ } -c_1\|u\|_{\w}^{r^+}- c_2\|f\|_{L^2(\O)}\|u\|_{\w},  $$ where $c_1,\ c_2$ are two positive constants.
   It follows, since $q^+ >p^->r^+$, that
operator $\mathcal{A}$ is strongly coercive. \\
Therefore, we have proved that the operator
$\mathcal{A}:W^{1,\mathcal{H}}_0(\Omega)\longrightarrow
\left(W^{1,\mathcal{H}}_0(\Omega)\right)^*$ is bounded,
pseudomonotone, and strongly coercive. Hence, by Theorem \ref{surjt}
and condition $(\mathcal{F}_2)$ , problem \eqref{P2} has a
nontrivial weak solution $u\in W^{1,\mathcal{H}}_0(\Omega)$. Thus,
the proof is completed.
\end{proof}

\section{Double phase problems with    exponents dependent on  the solution}\label{last}
In this section, we present the second type of the new class of
double-phase problem with exponents that depend on the solution
itself:
\begin{equation}\label{P1}
\left\lbrace\begin{array}{rll}
-\operatorname{div}\left(|\nabla u|^{p(x,u)-2} \nabla u+ \mu(x)|\nabla u|^{q(x,u)-2} \nabla u\right) &=f(x, u),  & \text { in } \Omega, \\
\ & \ \\
u &=0, & \text { on } \partial \Omega,
\end{array}\right.\tag{$\mathcal{P}_3$}
\end{equation}
where $ \Omega \subset \R^d,\ d \geq 2$ is a bounded domain with Lipschitz boundary $\partial \Omega.$\\

An important distinction emerges in the analysis of equations
\eqref{P} and \eqref{P1}. Specifically, in equation \eqref{P}, the
exponents are functions of the solution’s gradient. Conversely, in
equation \eqref{P1}, the exponents are functions of the solution
itself. Notably, the new operator introduced in equation \eqref{P1}
cannot be considered as $\Delta_\h$ (see \eqref{Ph}). Consequently,
the reasoning utilized in proving Theorems \ref{thmA1} and
\ref{thmA2} no longer applies. Hence, we adopt a different approach,
drawing inspiration
from the methodologies outlined in papers \cite{Anouar2024, Chipot2019}.\\

Equation \eqref{P1} is taken under the following assumptions:
\begin{itemize}
  \item [$(\mathrm{H}_1)$] $p,q: \O\times \mathbb{R} \rightarrow \mathbb{R}$ are two continuous functions such that
\vspace{2mm}
 \begin{itemize}
     \item [(i)]
$$
d<p^{-}=\ds\inf _{(x,t) \in \O \times \mathbb{R}} p(x,t) \leq p(x,t)
\leq p^{+}=\ds\sup _{(x,t) \in \O\times\mathbb{R}} p(x,t)<+\infty,
$$
$$q^-=\ds\inf _{(x,t) \in \O \times \mathbb{R}} q(x,t) \leq q(x,t) \leq q^{+}=\ds\sup _{(x,t) \in \O\times\mathbb{R}} q(x,t)<+\infty$$ and $$p(x,t)\leq q(x,t),\ \forall t \in \mathbb{R},\ \foral x\in \O.$$
\vspace{2mm}
\item[(ii)] There exist $c_p, \ c_q>0$ such that
$$|p(x,t_1)-p(y,t_2)|\leq c_p|t_1-t_2| \text{ and } |q(x,t_1)-q(y,t_2)|\leq c_q|t_1-t_2|,$$
$ \text{ for a.a. }x, y \in \O \text{ and for all }  t_1,t_2 \in \R.
$
\end{itemize}
  \item [$(\mathrm{H}_2)$]
 $f:\O \times \R\rightarrow \R $  is a Carathéodory function such that  $f(.,0)< 0$. Moreover, we assume that  there exist $1\leq r(x)<p^-$ and some positive constant $c$ such that $$|f(x,t)| \leq c \big(1+|t|^{r(x)-1}\big), \text{ for a.a. }x\in \O \text{ and all }t\in \R.$$
 \item [$(\mathrm{H}'_2)$]
 $f:\O \times \R\rightarrow \R $  is a Carathéodory function such that $f(.,0)< 0$. Moreover, we assume that there exist $1\leq r(x)<(p^-)^\ast$ and some positive constant $c$ such that $$|f(x,t)| \leq c \big(1+|t|^{r(x)-1}\big), \text{ for a.a. }x\in \O \text{ and all }t\in \R.$$
  \item [$(\mathrm{H}_3)$] $\mu \in L^\infty (\O)$ such that $\mu (x)\geq 0 \text{ a.a. } x \in \O.$
\item [$(\mathrm{H}_4)$] There exist $\eta_- < 0$ and $\eta_+ > 0$ such that
$$
f(x, \eta_+) < 0 < f(x, \eta_-),\text{ for a.a. }  x \in \O.
$$
\end{itemize}

 In the rest, we suppose that $(\mathrm{H}_1)$ is fulfilled and $u:\O\rightarrow \R$ is a continuous function. We define the functional space $E_u$ by
$$ E_u:= \bigg\{ v\in W^{1,1}_{0}(\O); \int_{\O}^{} |\nabla v|^{p(x,u)}+ \mu (x)| \nabla v|^{q(x,u)}dx < +\infty \bigg\},$$

with respect to the norm
$$\|v\|_{E_u}= \inf \bigg\{ \lambda >0;\ \int_\O \bigg( \left| \frac{v}{\l}\right|+  \big|\frac{ \nabla v}{\lambda}\big|^{p(x,u)}+\mu (x) \big|\frac{\nabla v|}{\lambda}\big|^{q(x,u)}\bigg)dx \leq 1 \bigg\}.$$
\begin{dfn}\label{d1}
A function $ u\in E_u$  is said to be a weak solution to equation
\eqref{P1} if it satisfies
$$\int_{\O} |\nabla u|^{p(x,u)-2}\nabla u \nabla w + \mu (x)|\nabla u|^{q(x,u)-2}\nabla u\nabla w d x = \int_{\O} f(x, u) w d x, \ \forall \ w \in E_u.$$
\end{dfn}
Our main results in the present section are given by the following
Theorems.
\begin{theorem} \label{th1}
Assume that $(\mathrm{H}_1)-(\mathrm{H}_3)$ hold. Then, problem
\eqref{P1} admits at least one nontrivial weak solution $u\in E_u$
in the sense of Definition \ref{d1}.
\end{theorem}

In proving Theorem \ref{th1}, it is important to note that we
employed a distinct approach for the auxiliary problem \eqref{P3},
differing from the methods used in \cite{bah, Chipot2019}. In these
references, the authors addressed problem \eqref{P1} under the
condition that $\mu \equiv 0$.

\begin{theorem} \label{th2}
Assume that $(\mathrm{H}_1)$, $(\mathrm{H}'_2)$, $(\mathrm{H}_3)$
and $(\mathrm{H}_4)$
 hold. Then, there exist   at least one nonpositive weak solution and at least one nonnegative weak solution to the problem \eqref{P1} in the sense of Definition \ref{d1}.
\end{theorem}
Theorem \ref{th2} initiates the investigation into the existence of
multiple solutions for equation \eqref{P1}. The primary tools
utilized include the surjectivity theorem (see Theorem \ref{surjt})
combined with the truncation method.

\subsection{Functional settings and further properties of $E_u$}

 If we fix a continuous function $u:\O\rightarrow \R$ and we set $p_1(x)=p(x,u(x))$ and $q_1(x)=q(x,u(x))$, then $E_u $ can be rewritten as follows:
$$
E_u:= \bigg\{ v\in W^{1,1}_{0}(\O);  |\nabla v| \in L^{\h_1}(\O)
\bigg\},
$$
with $\mathcal{H}_1:\O\times [0,+\infty[\rightarrow[0,+\infty[$ is a
N-function defined by $$\mathcal{H}_1 (x,t)=t^{p_1(x)}+\mu
(x)t^{q_1(x)}.$$ The functional space $E_u$ is equipped with the
following equivalent norm
$$\|v\|_{E_u}:= \|v\|_{L^1(\O)}  + \|\nabla v\|_{L^{\h_1}(\O)}.$$
The space $(E_u,\| .\|_{E_u})$ is a reflexive Banach space.
Moreover, we have:


\begin{prop}\label{prp2}
Assume that $(\mathrm{H}_1) $ holds. Then, there exists a constant
$C>0$ independent of $v$ such that
\begin{equation}\label{pt²}
  \|v\|_{L^1(\O)} \leq  C \|\nabla v\|_{L^{\h_1}(\O)}, \text{ for all }v\in E_u.
\end{equation}

\end{prop}
\begin{proof}
Noting that  the embedding $E_u \hookrightarrow L^1(\O)$   is
compact, then the proof is similar to  one of Proposition
\ref{prp3}.
\end{proof}
\begin{remark}
As a consequence of Proposition \ref{prp2} we  equip the space $E_u$
with the equivalent norm $\|\nabla v\|_{\h_1}.$
\end{remark}
\begin{prop}\label{GG}
  Let $v \in E_u$, $\left\{v_n\right\}_{n\in \mathbb{N}}\subseteq E_u$, then the following assertions hold true
  \begin{itemize}
    \item [(i)]$ \|v\|_{E_u} \leq 1\Rightarrow \|v\|_{E_u}^{q^+} \leq \ds\int_{\O}^{} |\nabla v |^{p(x,u)} + \mu (x) |\nabla v|^{q(x,u)} dx \leq \|v\|_{E_u}^{p^-}.$
    \item [(ii)] $\|v\|_{E_u} > 1\Rightarrow \|v\|_{E_u}^{p^-} \leq \ds\int_{\O}^{} |\nabla v |^{p(x,u)} + \mu (x) |\nabla v|^{q(x,u)} dx \leq \|v\|_{E_u}^{q^+}.$
    \item [(iii)]$\|v_n-v\|_{E_u}\rightarrow 0\Leftrightarrow  \ds \int_{\O}^{} |\nabla (v_n-v) |^{p(x,u)} + \mu (x) |\nabla (v_n-v)|^{q(x,u)} dx\rightarrow0.$
  \end{itemize}

\end{prop}

Now, we show that the space of infinitely differentiable functions
with compact support, $\mathcal{C}_0^\infty (\O)$, is dense in the
space $E_u$. This property will be essential to prove our existence
result. First, we recall an important inclusion in the following
remark, (see \cite[Remark 2.3]{CrespoBlanco2022}).

\begin{remark}\label{rem1}
Note that for a bounded domain $\Omega \subset \mathbb{R}^d$ and
$\gamma>d$ we have the following inclusions
$$
C^{0,1}(\overline{\Omega}) \subset W^{1, \gamma}(\Omega) \subset
C^{0,1-\frac{d}{\gamma}}(\overline{\Omega}) \subset C^{0,
\frac{1}{|\log t |}}(\overline{\Omega}) ,
$$
where $C^{0, \frac{1}{|\log t \mid}}(\overline{\Omega})$ is  the set
of all functions $\pi: \overline{\Omega} \rightarrow \mathbb{R}$
that are log-Hölder continuous, that is, there exists $C>0$ such
that
$$
|\pi(x)-\pi(y)| \leq \frac{C}{|\log | x-y||} \quad \text { for all }
x, y \in \overline{\Omega} \text { with }|x-y|<\frac{1}{2} .
$$
\end{remark}
\begin{prop}\label{p1}
  Let the assumptions $(\mathrm{H}_1)$ and $(\mathrm{H}_3)$ be satisfied. Moreover, we suppose that $u\in E_u$. Then, the space
  $\mathcal{C}_0^\infty (\O)$ is dense in $E_u.$
\end{prop}

\begin{proof}
We start the the proof by the following claim:
\\
\vspace{2mm}
\textbf{Claim. }  $p_1(\cdot)$ and $q_1(\cdot) $ are log-Hölder continuous.\\
Indeed; let $u \in E_u$. From the fact that  $p^- > d$ and Remark
\ref{rem1}, there exists a constant $C$ depending on $p$ and $d$
such that
\begin{equation}\label{1}
|u(x)-u(y)| \leq C|| u \|_{W^{1,
p^{-}}\left(\O\right)}|x-y|^{1-\frac{d}{p^{-}}}, \forall x, y \in
\overline{\O }.\end{equation} On the other side, by hypothesis
$(H_1)$, there is a constant $C_1>0$ such that
\begin{equation}\label{2}
\left| p_1(x)-p_1(y) \right|=|p(x,u(x))-p(y,u(y))| \leq
C_1|u(x)-u(y)|,\  \forall x, y \in \overline{\O}.
\end{equation}
Combining \eqref{1} with \eqref{2}, we deduce that
$$
\left| p_1(x)-p_1(y) \right|=\left|p(x,u(x))-p(y,u(y))\right| \leq
C_1 C\|u\|_{W^{1, p^{-}}(\O)}|x-y|^{1-\frac{d}{p^{-}}}, \ \forall x,
y \in \overline{\O }.
$$
Therefore, using Remark \ref{rem1}, $p_1(\cdot)$ is log-Hölder
continuous; that is, there is a constant $L>0$ such that
\begin{equation}\label{4}
\left| p_1(x)-p_1(y) \right| \leq \frac{-L}{\log |x-y|},\  \forall \
x, y \in \overline{\O},\ |x-y|<\frac{1}{2} .
\end{equation}
Similarly, we can easily show that $q_1(\cdot)$ is log-Hölder
continuous. This proves the claim. Then, the proof follows by
combining  the above claim with the argument used in the proof of
\cite[Theorem 2.6]{Fan2001}.
\end{proof}

\vspace{2mm}

\subsection{The existence of solution: Proof of Theorem \ref{th1}}
Notice that $E_u $
  is the natural functional space for solving problem
  \eqref{P1} using both classical variational and nonvariational methods. However, the dependence of the exponents on the solution itself presents
   a significant obstacle in establishing the fundamental topological properties of $ E_u $. Consequently, classical methods cannot be directly applied to equation
\eqref{P1}. Therefore, we adopt a perturbation approach as
introduced in \cite{Chipot2019}. Our proof is divided into two
parts. In the first part, we address an auxiliary problem defined on
a classical function space.
  In the second part, we employ a limiting process to establish the existence of a weak solution as per Definition
 \ref{d1}.
\subsubsection{The approximate problem}
 First, we recall the following classical monotonicity inequalities which will be crucial in the remainder of the work.
 \begin{lm}\cite[Lemma 2.1]{Misawa2023} \label{smn}
For every $p \in(1, \infty)$ there are positive constants
$C_j=C_j(d, p)(j=1,2,3)$ such that for all $\xi, \eta \in \R^d (d
\geq 1)$
$$
\left| | \xi|^{p-2} \xi-|\eta|^{p-2} \eta\right|\leq
C_1(|\xi|+|\eta|)^{p-2}\left| \xi-\eta \right|
$$
and
$$
\left(|\xi|^{p-2} \xi-|\eta|^{p-2} \eta\right) \cdot(\xi-\eta) \geq
C_2(|\xi|+|\eta|)^{p-2}|\xi-\eta|^2,
$$
where the symbol $\cdot$ denotes the inner product on $\R^d$. In
particular, when $p \geq 2$
$$
\left(|\xi|^{p-2} \xi-|\eta|^{p-2} \eta\right) \cdot(\xi-\eta) \geq
C_3|\xi-\eta|^p.
$$
\end{lm}

Now, for each $\e >0$, we define  the following auxiliary problem:
\begin{equation}\label{P3}
\left\lbrace\begin{array}{rll}
\mathcal{T}(u)-\e \operatorname{div}\left(|\nabla u|^{p^+-2} \nabla u+ \mu(x)|\nabla u|^{q^+-2} \nabla u\right) &=f(x, u),  & \text { in } \Omega, \\
\ & \ & \\
u &=0, & \text { on } \partial \Omega,
\end{array}\right.\tag{$\mathcal{P}_\e$}
\end{equation}
with $\mathcal{T}(u):=-\operatorname{div}\left(|\nabla u|^{p(x,
u)-2} \nabla u+ \mu(x)|\nabla u|^{q(x, u)-2} \nabla u\right).$
 \begin{theorem}\label{l1} Let $(\mathrm{H}_1)-(\mathrm{H}_3)$ be satisfied. Then, for each $\epsilon>0$, there exists a function $u_\epsilon \in W_0^{1, \mathcal{H}_{+}}\left(\Omega\right)$ such that
\begin{equation}\label{mb}
\begin{aligned}
& \int_{\Omega} \left|\nabla u_\epsilon\right|^{p\left(x,u_\epsilon\right)-2} \nabla u_\epsilon \nabla v d x+\int_{\Omega} \mu (x)\left|\nabla u_\epsilon\right|^{q(x,u_\epsilon)-2} \nabla u_\epsilon \nabla v d x \\
& +\epsilon\left(\int_{\Omega} \left|\nabla u_\epsilon\right|^{p^{+}-2} \nabla u_\epsilon \nabla v d x+\int_{\Omega} \mu (x)\left|\nabla u_\epsilon\right|^{q^{+}-2} \nabla u_\epsilon \nabla v d x\right) \\
& =\int_{\Omega} f\left(x, u_{\epsilon}\right)  v d x, \text{ for all } v \in W_0^{1, \mathcal{H}_{+}}\left(\Omega\right), \\
&
\end{aligned}
\end{equation}
where $\mathcal{H}_{+}: \Omega \times \R_+$ is defined by the correspondence $(x,t)\mapsto t^{p^+}+\mu (x)t^{q^{+}}$, for $ x \in \O $ and $ t\geq 0.$ 
\end{theorem}


\begin{proof}

 We define the operator $A^+:W^{1,\mathcal{H}_+}_0(\Omega)  \longrightarrow (W^{1,\mathcal{H}_+}_0(\Omega)) ^*$ as follows:
\begin{equation*}
\begin{aligned}
    &\langle A^+u,v \rangle = \int_{\Omega} |\nabla u|^{p(x,u) - 2}  \nabla u \nabla v + \mu (x) |\nabla u|^{q(x,u) - 2}  \nabla u \nabla v dx\\&  + \epsilon \int_{\Omega} |\nabla u|^{p^+ - 2}  \nabla u \nabla v + \mu (x) |\nabla u|^{q^+ - 2}  \nabla u \nabla v dx - \int_{\Omega} f(x, u) v dx, \text{ for all } v \in W^{1,\mathcal{H}_+}_0(\Omega).  \end{aligned}
\end{equation*}

\textbf{\underline{Step 1.}} The operator $A^+$ is strongly
coercive. By  Propositions \ref{zoo}, \ref{WW} and assumption
$(\mathrm{H}_2)$, we have
    $$\begin{aligned}
        \langle A^+u,u\rangle  =&\int_{\Omega} |\nabla u|^{p(x,u) }+ \mu (x) |\nabla u|^{q(x,u) } dx + \epsilon \int_{\Omega} |\nabla u|^{p^+}+ \mu (x) |\nabla u|^{q^+}dx - \int_{\Omega} f(x, u)u  dx\\
         \geq&\epsilon  \int_{\Omega} |\nabla u|^{p^+}+ \mu (x) |\nabla u|^{q^+}dx - \int_{\Omega} f(x, u)u  dx\\
        \geq & \epsilon  \min \big(  || u\|^{p^+}_{W_0^{1, \mathcal{H}_{+}}\left(\Omega\right)} , || u\|^{q^+} _{W_0^{1, \mathcal{H}_{+}}\left(\Omega\right)} \big)-\int_{\Omega} (1+|u|^{r(x)-1}) |u|dx\\
        \geq & \epsilon  \min \big(  || u\|^{p^+}_{W_0^{1, \mathcal{H}_{+}}\left(\Omega\right)} , || u\|^{q^+} _{W_0^{1, \mathcal{H}_{+}}\left(\Omega\right)} \big)-c\int_{\Omega} |u|^{r^+} dx-c\int_{\Omega} |u|^{r^-}dx- c\int_{\Omega} |u| dx \\
         \geq & \epsilon  \min \big(  || u\|^{p^+}_{W_0^{1, \mathcal{H}_{+}}\left(\Omega\right)} , || u\|^{q^+} _{W_0^{1, \mathcal{H}_{+}}\left(\Omega\right)} \big)-c_0|| u\|^{r^+} _{W_0^{1, \mathcal{H}_{+}}\left(\Omega\right)} -c_1|| u\|^{r^-} _{W_0^{1, \mathcal{H}_{+}}\left(\Omega\right)} -c_2|| u\|^{} _{W_0^{1, \mathcal{H}_{+}}\left(\Omega\right)},
    \end{aligned}$$

   which implies,
    since $q^+ >r^+$, that $A^+$ is strongly coercive in $   W_0^{1, \mathcal{H}_{+}}\left(\Omega\right)   $.\\

\textbf{\underline{Step 2.}}
 The operator $A^+$ satisfies the
 $(\text{S}_+)$-property. Let $\left\{ u_n\right\}_{n \in \mathbb{N}},  \subseteq W_0^{1, \mathcal{H}_{+}}\left(\Omega\right)$ suth that $ u_n \rightharpoonup u $
 in $W_0^{1, \mathcal{H}_{+}}\left(\Omega\right)$ and $\ds\limsup_{n \rightarrow +\infty} \langle A^+u_n, u_n-u \rangle \leq 0.$ We observe that
\begin{equation}\label{c0}
     \begin{aligned}
     \langle A^+u_n, u_n-u \rangle =& \int_{\O} \left|\nabla u_n\right|^{p(x,u_n)-2} \nabla u_n \nabla (u_n-u) +\mu (x)\left|\nabla u_n\right|^{q(x,u_n)-2} \nabla u_n \nabla (u_n-u) d x \\
& +\epsilon\left(\int_{\O} \left|\nabla u_n\right|^{p^{+}-2} \nabla u_n \nabla (u_n-u) +\mu (x)\left|\nabla u_n\right|^{q^{+}-2} \nabla u_n \nabla (u_n-u) d x\right) \\
& -\int_{\O} f\left(x, u_n\right) (u_n-u) d x.
    \end{aligned}
\end{equation}
\textbf{Claim.} We have \begin{equation}\label{j8} \int_{\O}
\left|\nabla u_n\right|^{p(x,u_n)-2} \nabla u_n \nabla (u_n-u) +\mu
(x)\left|\nabla u_n\right|^{q(x,u_n)-2} \nabla u_n \nabla (u_n-u) d
x\geq o_n(1).\end{equation} Indeed:  In light of Lemma \ref{smn}, we
see that
     \begin{equation}\label{j5}
\begin{aligned}
&\int_{\O} \left|\nabla u_n\right|^{p(x,u_n)-2} \nabla u_n \nabla (u_n-u) +\mu (x)\left|\nabla u_n\right|^{q(x,u_n)-2} \nabla u_n \nabla (u_n-u) d x\\
&= \int_{\O}\big( \left|\nabla u_n\right|^{p(x,u_n)-2} \nabla u_n -\left|\nabla u\right|^{p(x,u_n)-2} \nabla u\big) \nabla (u_n-u)dx  \\ &+\int_{\O} \mu (x)\big (\left|\nabla u_n\right|^{q(x,u_n)-2} \nabla u_n -\left|\nabla u\right|^{q(x,u_n)-2} \nabla u \big) \nabla (u_n-u) d x\\
&+ \int_{\O} \left|\nabla u\right|^{p(x,u_n)-2} \nabla u \nabla (u_n-u) +\mu (x)\left|\nabla u\right|^{q(x,u_n)-2} \nabla u \nabla (u_n-u) d x\\
&\geq  \int_{\O} \left|\nabla u\right|^{p(x,u_n)-2} \nabla u \nabla (u_n-u) +\mu (x)\left|\nabla u\right|^{q(x,u_n)-2} \nabla u \nabla (u_n-u) d x\\
&= \int_{\O} \left|\nabla u\right|^{p(x,u)-2} \nabla u \nabla (u_n-u) +\mu (x)\left|\nabla u\right|^{q(x,u)-2} \nabla u \nabla (u_n-u) d x\\
&+\int_{\O}\big( \left|\nabla u\right|^{p(x,u_n)-2} \nabla u -\left|\nabla u\right|^{p(x,u)-2} \nabla u\big) \nabla (u_n-u)dx \\&+\int_{\O}\mu (x)\big (\left|\nabla u\right|^{q(x,u_n)-2} \nabla u -\left|\nabla u\right|^{q(x,u)-2} \nabla u \big) \nabla (u_n-u) d x.\\
\end{aligned}
\end{equation}
On the other hand, since   $u_n\rightharpoonup u$ in $W_0^{1,
\mathcal{H}_{+}}\left(\Omega\right)$, we have
\begin{equation}\label{j6}
   \int_{\O} \left|\nabla u\right|^{p(x,u)-2} \nabla u \nabla (u_n-u) +\mu (x)\left|\nabla u\right|^{q(x,u)-2} \nabla u \nabla (u_n-u) d x \rightarrow0, \text{ as } n\rightarrow +\infty.
\end{equation}
Moreover, beying $u_n\rightharpoonup u$ in $W^{1,\h_+}_0(\O)$ and
using  the Lebesgue dominated convergence theorem, we prove that
\begin{equation}\label{j7}
\begin{aligned}
 & \int_{\O}\big( \left|\nabla u\right|^{p(x,u_n)-2} \nabla u -\left|\nabla u\right|^{p(x,u)-2} \nabla u\big) \nabla (u_n-u) dx \\&+\int_{\O} \mu (x)\big (\left|\nabla u\right|^{q(x,u_n)-2} \nabla u -\left|\nabla u\right|^{q(x,u)-2} \nabla u \big) \nabla (u_n-u) d x\rightarrow 0, \text{ as } n \rightarrow +\infty.
\end{aligned}\end{equation}
Thus, by \eqref{j5}, \eqref{j6} and \eqref{j7}, we get the proof of the Claim.\\
Now, invoking $(\mathrm{H}_2)$ and Proposition \ref{prp}, we obtain
\begin{equation}\label{j9}
  \lim_{n \rightarrow +\infty}\int_{\O}^{} f(x,u_n)(u_n-u) dx =0.
\end{equation}
Combining the above claim with \eqref{c0} and \eqref{j9}, it yields
that
\begin{equation}\label{c1}
\begin{aligned}
&o_n(1) +\epsilon\small\int_{\O} \left(\left|\nabla
u_n\right|^{p^{+}-2} \nabla u_n \nabla (u_n-u) +\mu (x)\left|\nabla
u_n\right|^{q^{+}-2} \nabla u_n \nabla (u_n-u) \right)d x \\& \leq
\langle A^+ u_n, u_n-u \rangle.
\end{aligned}
\end{equation}
Moreover, due to Lemma \ref{smn}, we know that
$$
\begin{aligned}
&\int_{\O} \left|\nabla u_n\right|^{p^{+}-2} \nabla u_n \nabla (u_n-u) +\mu (x)\left|\nabla u_n\right|^{q^{+}-2} \nabla u_n \nabla (u_n-u) d x\\
& =\int_{\O} \big(\left|\nabla u_n\right|^{p^{+}-2} \nabla u_n-\left|\nabla u\right|^{p^{+}-2} \nabla u \big) \nabla (u_n-u)dx \\&\ \ +\int_{\O}^{}\mu (x)\big(\left|\nabla u_n\right|^{q^{+}-2} \nabla u_n-\left|\nabla u\right|^{q^{+}-2} \nabla u\big) \nabla (u_n-u) d x+o_n(1)\\
& \geq \int_{\O} \big|\nabla (u_n-u) \big|^{p^+} +\mu (x) \big|\nabla (u_n-u) \big|^{q^+} dx+o_n(1)= \rho_{\mathcal{H}_+}(u_n-u)+o_n(1).\\
\end{aligned}
$$
Then, the inequality \eqref{c1} becomes
\begin{equation}\label{c2}
   \rho_{\mathcal{H}_+}(u_n-u)+o_n(1)   \leq   \langle A^+u_n, u_n-u \rangle.
\end{equation}
Therefore, passing to the limit as $n \rightarrow +\infty$
in\eqref{c2}, we get
$$0 \leq \lim_{n \rightarrow  +\infty} \rho_{\mathcal{H}_+}(u_n-u) \leq  \lim_{n \rightarrow  +\infty}  \langle A^+u_n, u_n-u \rangle \leq \limsup_{n \rightarrow  +\infty}  \langle A^+u_n, u_n-u \rangle \leq0,$$
which proves, using Proposition \ref{zoo}, the Step 2.\\

 \textbf{\underline{Step 3.}}
 \vspace{3.3mm}
The operator $A^+$ is pseudomonotone.\\
\textbf{Claim.} $A^+$ is continuous from $W_0^{1,
\mathcal{H}_{+}}\left(\Omega\right)$ strong into $\big( W_0^{1,
\mathcal{H}_{+}}\left(\Omega\right)\big)^\ast$ weak.
 Let $\left\{ u_n\right\}_{n \in \mathbb{N}} \ \subseteq W^{1,\mathcal{H}_+}_0(\Omega)$ such that
$u_n \rightarrow u$ in $W^{1,\mathcal{H}_+}_0(\Omega)$. Then, for
every $v\in W^{1,\mathcal{H}_+}_0(\Omega)$, one has
\begin{equation}\label{MR}
\begin{aligned}
 \small\langle A^+u_n-A^+u,v \rangle =& \int_{\O} \big ( \left|\nabla u_n\right|^{p(x,u_n)-2} \nabla u_n -\left|\nabla u\right|^{p\left(x,u\right)-2} \nabla u\big) \nabla vdx \\
 &+ \int_{\O}\mu (x)\big (\left|\nabla u_n\right|^{q(x,u_n)-2} \nabla u_n -\left|\nabla u\right|^{q(x,u)-2} \nabla u \big) \nabla v d x\\
 &+ \int_{\O} \big ( \left|\nabla u_n\right|^{p^+-2} \nabla u_n -\left|\nabla u\right|^{p^+-2} \nabla u\big) \nabla v dx \\
 &+ \int_{\O}\mu (x)\big (\left|\nabla u_n\right|^{q^+-2} \nabla u_n -\left|\nabla u\right|^{q^+-2} \nabla u \big) \nabla v d x\\
 &-\int_{\O}^{} \big( f(x,u_n)-f(x,u)\big) v dx .
\end{aligned}
\end{equation}
In what follows, we prove that for every $v \in W^{1,\h_+}_0(\O)$
\begin{equation}\label{YB}
\begin{aligned}&\int_{\O}\left( \big ( \left|\nabla u_n\right|^{p(x,u_n)-2} \nabla u_n -\left|\nabla u\right|^{p(x,u)-2} \nabla u\big) \nabla v\right) dx\\&+\int_{\O}^{}\left(\mu (x)\big (\left|\nabla u_n\right|^{q(x,u_n)-2} \nabla u_n -\left|\nabla u\right|^{q(x,u)-2} \nabla u \big) \nabla v\right)dx\\& \ \  \longrightarrow 0, \text{ as } n \rightarrow +\infty.\end{aligned}
\end{equation}
First, for every $v \in W^{1,\h_+}_0(\O)$, we write
  \begin{equation}\label{MR1}
  \begin{aligned}
 &\int_{\O} \big ( \left|\nabla u_n\right|^{p(x,u_n)-2} \nabla u_n -\left|\nabla u\right|^{p(x,u)-2} \nabla u\big) \nabla vdx\\&=\int_{\O} \big ( \left|\nabla u_n\right|^{p(x,u_n)-2} \nabla u_n -\left|\nabla u\right|^{p(x,u_n)-2} \nabla u\big) \nabla vdx\\
  &\ \ + \int_{\O} \big ( \left|\nabla u\right|^{p(x,u_n)-2} \nabla u_n -\left|\nabla u\right|^{p(x,u)-2} \nabla u\big) \nabla vdx.\\
  \end{aligned}
 \end{equation}
By Hölder inequality and Lemma \ref{smn}, we have
  \begin{equation}\label{MR2}
\small\small\begin{aligned}
& \int_{\O} \bigg|\big ( \left|\nabla u_n\right|^{p(x,u_n)-2} \nabla u_n -\left|\nabla u\right|^{p(x,u_n)-2} \nabla u\big) \nabla v\bigg|dx\\
  &\leq  \bigg( \int_{\O} \left| \left|\nabla u_n\right|^{p(x,u_n)-2} \nabla u_n -\left|\nabla u\right|^{p(x,u_n)-2} \nabla u\right|^{\frac{p^+}{p^+-1}}dx \bigg)^{\frac{p^+-1}{p^+}} \|\nabla v\|_{L^{p^+}(\O) } \\
  &\leq c \bigg( \int_{\O} \big( \left|\nabla u_n\right|+\left|\nabla u\right|\big) ^{(p(x,u_n) - 2) \frac{p^+}{p^+-1}} \big|\nabla (u_n-u) \big|^{\frac{p^+}{p^+-1}}dx \bigg)^{\frac{p^+-1}{p^+}} \|\nabla v\|_{L^{p^+}(\O) } \\
   &\leq c \bigg( \int_{\{\left|\nabla u_n\right|+\left|\nabla u\right|\geq 1\}} \big | \left|\nabla u_n\right|+\left|\nabla u\right|\big) ^{(p^+ - 2) \frac{p^+}{p^+-1}} \big|\nabla (u_n-u) \big|^{\frac{p^+}{p^+-1}}dx \bigg)^{\frac{p^+-1}{p^+}} \|\nabla v\|_{L^{p^+}(\O) }\\
  &\ \  + c \bigg( \int_{\{\left|\nabla u_n\right|+\left|\nabla u\right|\leq 1\}} \big|\nabla (u_n-u) \big|^{\frac{p^+}{p^+-1}}dx \bigg)^{\frac{p^+-1}{p^+}} \|\nabla v\|_{L^{p^+}(\O) } \\
     &\leq c \bigg( \int_{\O} \big ( \left|\nabla u_n\right|+\left|\nabla u\right|\big )^{p^+} dx \bigg)^{\frac{p^+-2}{p^+}}  \|\nabla (u_n-u) \|_{L^{p^+}(\O) } \|\nabla v\|_{L^{p^+}(\O) }
   \\&+ c  \|\nabla (u_n-u) \|_{L^{\frac{p^+}{p^+-1}}(\O)}\|\nabla \|_{L^{p^+}(\O) } = o_n(1).
\end{aligned}
  \end{equation}
On the other hand, using the Lebesgue dominated convergence theorem,
we infer that
  \begin{equation}\label{MR3}
  \int_{\O} \big ( \left|\nabla u\right|^{p(x,u_n)-2} \nabla u_n -\left|\nabla u\right|^{p(x,u)-2} \nabla u\big) \nabla vdx =o_n(1).
  \end{equation}
  Therefore, combining \eqref{MR1}, \eqref{MR2} and \eqref{MR3}, we get
  \begin{equation}\label{MR4}
      \int_{\O} \big ( \left|\nabla u_n\right|^{p(x,u_n)-2} \nabla u_n -\left|\nabla u\right|^{p(x,u)-2} \nabla u\big) \nabla vdx \longrightarrow 0, \text{ as } n \rightarrow +\infty.
  \end{equation}
  Similarly, we can show that
  \begin{equation}\label{MR5}
      \int_{\O}
  \mu (x)\big (\left|\nabla u_n\right|^{q(x,u_n)-2} \nabla u_n -\left|\nabla u\right|^{q(x,u)-2} \nabla u \big) \nabla vdx \longrightarrow 0, \text{ as } n \rightarrow +\infty.
  \end{equation}
  Consequently, by \eqref{MR4} and \eqref{MR5}, we get the proof of \eqref{YB}.\\
 Recall that $u_n \rightarrow u$ in $W^{1,\h_+}_0(\O)$ and using condition $(\mathrm{H_2})$, Proposition \ref{prp},  Hölder inequality and   the Lebesgue dominated  convergence theorem, we conclude that
 \begin{equation}\label{MR6}
   \begin{aligned}&  \int_{\O} \big ( \left|\nabla u_n\right|^{p^+-2} \nabla u_n -\left|\nabla u\right|^{p^+-2} \nabla u\big) \nabla v
 + \mu (x)\big (\left|\nabla u_n\right|^{q^+-2} \nabla u_n -\left|\nabla u\right|^{q^+-2} \nabla u \big) \nabla v d x\\& \longrightarrow 0, \text{ as } n \rightarrow +\infty,\end{aligned}
 \end{equation}
 and
\begin{equation}\label{MR7}
    \lim_{n \rightarrow +\infty}  \int_{\O}^{}  (f(x,u_n)-f(x,u))vdx =0.
  \end{equation}
Thus, from \eqref{MR}, \eqref{YB}, \eqref{MR6} and \eqref{MR7}, we
conclude   that
   $$\small\langle A^+u_n-A^+u,v \rangle \longrightarrow 0, \text{ as } n \rightarrow+\infty, \text{ for all } v \in W^{1,\h_+}_0(\O).$$
   This proves the claim.\\ It follows, by combining  Step 2 and the above claim, that $A^+$ is pseudomonotone.\\

\textbf{\underline{Step 4.}} The operator $A^+$ is bounded. First,
we introduce the Nemytskij operator $N_f: W^{1,\h_+}_0(\O)
\rightarrow (W^{1,\h_+}_0(\O))^\ast$ associated to $f$, by
$$\langle N_f(u), v \rangle =\int_\O f(x,u)v dx, \ \text{for all } v \in W^{1,\h_+}_0(\O).$$
From $(\mathrm{H}_2)$ it is clear that $N_f$ is well-defined and
bounded. It remains  to prove  that $A^+-N_f$ is bounded. For this
aim, taking $u, v \in W_0^{1, \mathcal{H}_+}(\Omega)
\backslash\{0\}$, we have:
\begin{equation*}
\begin{aligned}
  \left\vert\langle A^+(u)-N_f(u) ,v  \rangle\right\vert \leq& \int_{\O}^{} |\nabla u|^{p(x,u)-1} |\nabla v| + \mu (x) |\nabla u|^{q(x,u)-1} |\nabla v| dx\\
                                          & + \epsilon \int_{\O}^{} |\nabla u|^{p^+-1} |\nabla v| + \mu (x) |\nabla u|^{q^+-1} |\nabla v| dx\\
                                          &\leq (1+\e)  \int_{\O}^{} |\nabla u|^{p^+-1} |\nabla v| + \mu (x) |\nabla u|^{q^+-1} |\nabla v| dx\\
                                          &\ \ + \int_{\O}^{} |\nabla u|^{p^--1} |\nabla v| + \mu (x) |\nabla u|^{q^--1} |\nabla v| dx.
\end{aligned}
\end{equation*}
\\
Hence, $A^+$ is bounded. This ends the proof of Step 4.\vspace{2mm}\\
Consequently,  by Steps 1, 2, 3 and 4, and Theorem \ref{surjt} there
exists $u_\e \in W^{1,\mathcal{H}_+}_0(\O)$ satisfies \eqref{mb}.
Moreover, from condition $(\mathrm{H}_2)$ ($f(.,0) \neq0$), we can
show that $u_\e \not\equiv 0$. Thus, the proof is completed.
\end{proof}
\subsubsection{Passage to the limit }
 Choosing $\epsilon =\dfrac{1}{n}$ in Theorem \ref{l1}, we deduce that there exists $u_n \in W_0^{1, \mathcal{H}_{+}}\left(\Omega\right)$ such that
\begin{equation}\label{e15}
\begin{aligned}
& \int_{\O}\left( \left|\nabla u_n\right|^{p(x,u_n)-2} \nabla u_n \nabla v +\mu (x)\left|\nabla u_n\right|^{q(x,u_n)-2} \nabla u_n \nabla v\right) d x \\
& +\frac{1}{n}\int_{\O}\left( \left|\nabla u_n\right|^{p^{+}-2} \nabla u_n \nabla v +\mu (x)\left|\nabla u_n\right|^{q^{+}-2} \nabla u_n \nabla v \right) dx\\
= & \int_{\O} f\left(x, u_n\right) v d x, \text{ for all } v \in
W_0^{1, \mathcal{H}_{+}}\left(\Omega\right)  .
\end{aligned}
\end{equation}
Let $\mathcal{H_-}: \O \times [0,+\infty) \rightarrow [0,+\infty)$ be defined as $$\mathcal{H_-} (x,t):= t^{p^-}+\mu (x)t^{q^-}, \text{  for all } (x,t) \in \O\times [0,+\infty).$$\\

\textbf{\underline{Step 1.}} $\left\{ u_n\right\}_{n \in
\mathbb{N}}$ is bounded in $ W^{1 ,\mathcal{H_-}}_0 (\O).$ Taking $v
= u_n$ in \eqref{e15}, one has
\begin{equation}\label{e16}
\begin{aligned}
& \int_{\O} \left(\left|\nabla u_n\right|^{p(x,u_n)}  +\mu
(x)\left|\nabla u_n\right|^{q(x,u_n)} \right) d x +\frac{1}{n}
\int_{\O} \left|\nabla u_n\right|^{p^+}dx
\\
 & \leq\int_{\O} f\left(x, u_n\right) u_n d x ,
\end{aligned}
\end{equation}
and
\begin{equation}\label{e160}
\begin{aligned}
& \int_{\O} \left(\left|\nabla u_n\right|^{p(x,u_n)}  +\mu
(x)\left|\nabla u_n\right|^{q(x,u_n)} \right) d x +\frac{1}{n}
\int_{\O} \mu(x) \left|\nabla u_n\right|^{q^+}dx
\\
 & \leq\int_{\O} f\left(x, u_n\right) u_n d x ,
\end{aligned}
\end{equation}
 By  $(\mathrm{H}_2)$,  Proposition \ref{prp} and the Poincaré-type inequality in $W^{1,p^-}_0(\O)$, we obtain
 $$\begin{aligned}
\bigg|  \int_{\O} f\left(x, u_n\right) u_n d x \bigg| &\leq c
\int_{\O}\left( |u_n| + |u_n|^{r(x)}\right)dx
 \leq c \int_{\O} |u_n| ^{r^+}dx +c|\O|\\
 &\leq c'\bigg(\int_{\O} | u_n| ^{p^-}dx \bigg)^{\frac{r^+}{p^-}}+c|\O|&\\
  & \leq c''\bigg(\int_{\O} |\nabla u_n| ^{p^-} dx \bigg)^{\frac{r^+}{p^-}} +c|\O|\\
 & \leq c''\bigg(\int_{\O} |\nabla u_n| ^{p(x,u_n)} dx +|\O| \bigg)^{\frac{r^+}{p^-}} +c|\O|\\
 \end{aligned}$$

It follows, using \eqref{e16} and \eqref{e160}, that
\begin{equation}\label{e70}
\left.
\begin{array}{ll}
&\ds \int_{\O} \left(\left|\nabla u_n\right|^{p(x,u_n)}  +\mu
(x)\left|\nabla u_n\right|^{q(x,u_n)} \right) d x +\frac{1}{n}
\int_{\O} \left|\nabla u_n\right|^{p^+}dx
\\
&\\
  &\ds\int_{\O} \left(\left|\nabla u_n\right|^{p(x,u_n)}  +\mu (x)\left|\nabla u_n\right|^{q(x,u_n)} \right) d x +\frac{1}{n}
\int_{\O} \mu(x) \left|\nabla u_n\right|^{q^+}dx
  \end{array}
\right\} \leq c_2, \ \foral n \geq 1,
\end{equation}
where $c_2>0$ is a positive constant.\\
Thus there exists $c_3 > 0$ independent of $n$, such that,
\begin{equation}\label{e18}
\begin{aligned}
          \int_{\O} \left( \left|\nabla u_n\right|^{p^-}  +\mu (x)\left|\nabla u_n\right|^{q^-} \right) d x & \leq \int_{\left\{|\nabla u_n|>1\right\}} \left(\left|\nabla u_n\right|^{p(x,u_n)}  +\mu (x)\left|\nabla u_n\right|^{q(x,u_n)}\right)  d x  \\
           &+ \int_{\left\{|\nabla u_n|\leq1\right\}}^{} (1+\mu (x))dx \\
           & \leq \int_{\O} \left(\left|\nabla u_n\right|^{p(x,u_n)}  +\mu (x)\left|\nabla u_n\right|^{q(x,u_n)}\right)  d x +|\O |+\|\mu\|_\infty \\
           & \leq c_3.
        \end{aligned}\end{equation}
        This proves step 1.


By Step 1 and the reflexivity of the Sobolev space $ W^{1 ,\mathcal{H_-}}_0 (\O)$, there exists $u\in  W^{1 ,\mathcal{H_-}}_0 (\O)$ such that, up to a subsequence, $ u_n \rightharpoonup u$ weakly in  $W^{1 ,\mathcal{H_-}}_0 (\O)$  and $ u_n(x) \rightarrow u(x) \text{ a.a. }x \in \O.$\\

\textbf{\underline{Step 2.}} We prove that $u \in E_u$, namely,
\begin{equation}\label{e19}
  \int_{\O}^{}\left(|\nabla u|^{p(x,u)} +\mu (x) |\nabla u|^{q(x,u)}\right) dx < \infty.
\end{equation}
To achieve that, set $p_n(x)=p(x,u_n(x)), \ p_u(x)=p(x,u(x)),\
q_n(x)=q(x,u_n(x)),$ and $q_u(x)=q(x,u(x)).$ From Step 1, up to a
subsequence, we have
\begin{equation} \label{chp}
\begin{cases}
   u_n \rightharpoonup u \text{ in } W^{1 ,\mathcal{H_-}}_0 (\O) ,\\
   \nabla u_n \rightharpoonup \nabla u \text{ in } (L^{\h_-}(\O))^d,\\
   p_n(x) \rightarrow p_u(x), \ q_n(x)\rightarrow q_u(x) \text{ a.a. } x \in \O.
\end{cases}
\end{equation}
On the other hand, by \eqref{e70}, we infer that
\begin{equation}\label{e17}
    \int_{\O} \left(\left|\nabla u_n\right|^{p(x,u_n)}  +\mu (x)\left|\nabla u_n\right|^{q(x,u_n)} \right) d x \leq  c_2.
\end{equation}
Hence, using \eqref{chp}, \eqref{e17} and Lemma \ref{lem}, we
conclude the proof of Step 2. \vspace{2mm}

 \textbf{\underline{Step 3.}} $u$ satisfies the following inequality
 \begin{equation}\label{e35}
  \begin{aligned}
    \int_{\O}^{} f(x,u)(u-v) dx \geq   &  \int_{\O}\left|\nabla v\right|^{p(x,u)-2} \nabla v\nabla (u-v)dx \\
   & +  \int_{\O}\mu (x) \left|\nabla v\right|^{q(x,u)-2} \nabla v\nabla (u-v)  dx, \  \text{ for all } v \in \mathcal{C}^\infty_0(\O).
  \end{aligned}
\end{equation}
\textbf{Claim 1.} We prove that
\begin{equation}\label{e27}
  \begin{aligned}
    \int_{\O}^{} f(x,u)(u-v) dx \geq   &  \lim_{n \rightarrow +\infty}\left(\int_{\O}\left|\nabla v\right|^{p(x,u_n)-2} \nabla v\nabla (u_n-v)dx\right. \\
   &  \left. + \int_{\O}\mu (x) \left|\nabla v\right|^{q(x,u_n)-2} \nabla v\nabla (u_n-v)  dx\right), \foral v \in  W^{1 ,\mathcal{H_+}}_0 (\O).
  \end{aligned}
\end{equation}
Let $v \in W^{1 ,\mathcal{H_+}}_0 (\O)$. We have, for all $n\geq1$,

$$\begin{aligned}
 & \int_{\O}\big( \left|\nabla u_n\right|^{p(x,u_n)-2} \nabla u_n -\left|\nabla v\right|^{p(x,u_n)-2} \nabla v\big) \nabla (u_n-v) dx \\&+\int_{\O} \mu (x)\big (\left|\nabla u_n\right|^{q(x,u_n)-2} \nabla u_n -\left|\nabla v\right|^{q(x,u_n)-2} \nabla v \big) \nabla (u_n-v) d x\\
 & +\frac{1}{n}\int_{\O}\big(  \left|\nabla u_n\right|^{p^{+}-2} \nabla u_n - \left|\nabla v\right|^{p^{+}-2} \nabla v  \big) \nabla (u_n-v)dx\\ &+\frac{1}{n}\int_{\O} \mu (x)\big (\left|\nabla u_n \right|^{q^{+}-2} \nabla u_n -|\nabla v|^{q^+-2}\nabla v \big) \nabla (u_n-v) d x \\
 &\geq 0, \text{ for all } n\geq 1.
\end{aligned}$$
Thus, $$
\begin{aligned}
   & \int_{\O}\big(\left|\nabla u_n\right|^{p(x,u_n)-2} \nabla u_n + \frac{1}{n} \left|\nabla u_n\right|^{p^{+}-2} \nabla u_n\big) \nabla (u_n-v)dx \\
   & + \int_{\O}\mu (x) \big(\left|\nabla u_n\right|^{q(x,u_n)-2} \nabla u_n + \frac{1}{n} \left|\nabla u_n\right|^{q^{+}-2} \nabla u_n\big) \nabla (u_n-v)dx \\
  \geq  &  \int_{\O}\big(\left|\nabla v\right|^{p(x,u_n)-2} \nabla v+ \frac{1}{n} \left|\nabla v\right|^{p^{+}-2} \nabla v\big) \nabla (u_n-v)dx\\
   & + \int_{\O}\mu (x) \big(\left|\nabla v\right|^{q(x,u_n)-2} \nabla v + \frac{1}{n} \left|\nabla v\right|^{q^{+}-2} \nabla v\big) \nabla (u_n-v)dx.
\end{aligned}
$$
It follow, in view of \eqref{e15}, that
\begin{equation}\label{e20}
  \begin{aligned}
    \int_{\O}^{} f(x,u_n)(u_n-v) dx \geq   &  \int_{\O}\big(\left|\nabla v\right|^{p(x,u_n)-2} \nabla v+ \frac{1}{n} \left|\nabla v\right|^{p^{+}-2} \nabla v\big) \nabla (u_n-v)dx\\
   & + \int_{\O}\mu (x) \big(\left|\nabla v\right|^{q(x,u_n)-2} \nabla v + \frac{1}{n} \left|\nabla v\right|^{q^{+}-2} \nabla v\big) \nabla (u_n-v)dx.
  \end{aligned}
\end{equation}
On the other side, by \eqref{e70} and Hölder inequality, we have
\begin{equation}\label{e21}
  \begin{aligned}
    &\bigg| \frac{1}{n} \int_{\O}  \left|\nabla v\right|^{p^{+}-2} \nabla v \nabla (u_n-v)dx  \bigg| \\&\leq \frac{1}{n}  \int_{\O}  \left|\nabla v\right|^{p^{+}-1} | \nabla (u_n-v)|dx \\
    &\leq  \frac{1}{n} \bigg( \int_{\O}  \left|\nabla v\right|^{p^{+} } \bigg)^{\frac{p^+-1}{p^+}} \bigg(\int_{\O} | \nabla (u_n-v)|^{p^+} \bigg)^{\frac{1}{p^+}} \\
     & =  \left( \frac{1}{n}\right)^{\frac{p^+-1}{p^+}}       \bigg( \int_{\O}  \left|\nabla v\right|^{p^{+} } \bigg)^{\frac{p^+-1}{p^+}} \bigg(\frac{1}{n}\int_{\O} | \nabla (u_n-v)|^{p^+} \bigg)^{\frac{1}{p^+}} \\
      & \leq \left( \frac{1}{n}\right)^{\frac{p^+-1}{p^+}}   \bigg(\frac{1}{n}\int_{\O} | \nabla (u_n-v)|^{p^+} \bigg)^{\frac{1}{p^+}} \|v\|_{W^{1 ,\mathcal{H_+}}_0 (\O)}^{p^+-1}.\\
     &\leq c_2  \left( \frac{1}{n}\right)^{\frac{p^+-1}{p^+}}  \|v\|_{W^{1 ,\mathcal{H_+}}_0 (\O)}^{p^+-1},
  \end{aligned}
\end{equation}
and
\begin{equation}\label{e22}
   \bigg| \frac{1}{n} \int_{\O} \mu (x) \left|\nabla v\right|^{q^{+}-2} \nabla v \nabla (u_n-v)dx  \bigg|
   \leq   c_2  \left( \frac{1}{n}\right)^{\frac{q^+-1}{q^+}}  \|v\|_{W^{1 ,\mathcal{H_+}}_0 (\O)}^{q^+-1}.
\end{equation}
Then, combining \eqref{e21} and \eqref{e22}, we get
\begin{equation}\label{e23}
   \frac{1}{n} \int_{\O}  \left|\nabla v\right|^{p^{+}-2} \nabla v \nabla (u_n-v)dx+\frac{1}{n} \int_{\O} \mu (x) \left|\nabla v\right|^{q^{+}-2} \nabla v \nabla (u_n-v)dx \rightarrow0,\text{ as } n\rightarrow +\infty.
\end{equation}
Moreover, from $(\mathrm{H}_2)$, Proposition \ref{prp} and the
Lebesgue dominated convergence Theorem, one has
\begin{equation}\label{e24}
                  \int_{\O}^{} f(x,u_n)(u_n-v)dx \rightarrow\int_{\O}^{} f(x,u)(u-v)dx, \text{ as } n \rightarrow +\infty.
                \end{equation}
Consequently, from \eqref{e20}, \eqref{e23} and \eqref{e24}, we conclude the proof of Claim 1.\\
\textbf{Claim 2.} We prove that
\begin{equation}\label{KS}
\begin{aligned}
&\lim_{n \rightarrow +\infty }\int_{\O}\left(\left(\left|\nabla
v\right|^{p(x,u_n)-2} \nabla v+\mu (x) \left|\nabla
v\right|^{q(x,u_n)-2} \nabla v\right)\nabla (u_n-v)\right)dx\\&\ =
\int_{\O}\left(\left(\left|\nabla v\right|^{p(x,u)-2} \nabla v+\mu
(x) \left|\nabla v\right|^{q(x,u)-2} \nabla v\right)\nabla
(u-v)\right)dx, \ \foral v \in  W^{1 ,\mathcal{H_+}}_0 (\O).
\end{aligned}
\end{equation}
We know that
\begin{equation}\label{e28}
  \begin{aligned}
&\int_{\O}\left(\left(\left|\nabla v\right|^{p(x,u_n)-2} \nabla
v+\mu (x) \left|\nabla v\right|^{q(x,u_n)-2} \nabla v\right)\nabla
(u_n-v)\right)dx\\& =\int_{\O}\left(\left(\left|\nabla
v\right|^{p(x,u)-2} \nabla v+\mu (x) \left|\nabla
v\right|^{q(x,u)-2} \nabla v\right)\nabla (u_n-v)\right)dx\\& \ \
+\int_{\O}\left(\big(\left|\nabla v\right|^{p(x,u_n)-2} \nabla
v-\left|\nabla v\right|^{p(x,u)-2} \nabla v
\big)\right)\nabla(u_n-v)dx\\&+\int_{\O}^{}\left(\mu (x)
\big(\left|\nabla v\right|^{q(x,u_n)-2} \nabla v-\left|\nabla
v\right|^{q(x,u)-2} \nabla v \big) \nabla (u_n-v)\right)dx.
\end{aligned}
\end{equation}
Using the weak convergence of $\left(u_n-v\right)_{n\in \mathbb{N}}$
to $(u-v)$ in $W_0^{1, \mathcal{H}_{+}}\left(\Omega\right)$, we get
 \begin{equation}\label{e30}\begin{aligned}
&\int_{\O}^{}|\nabla v|^{p(x,u)-2}  \nabla v \nabla (u_n-v)+\mu (x)|\nabla v|^{q(x,u)-2}  \nabla v \nabla (u_n-v) dx \\
&\longrightarrow \int_{\O}^{}|\nabla v|^{p(x,u)-2}  \nabla v \nabla
(u -v)+\mu (x)|\nabla v|^{q(x,u)-2}  \nabla v \nabla (u-v)
dx.\end{aligned}\end{equation}
Moreover, we have
\begin{equation}\label{e31}
\begin{aligned}
& \left|\int_{\O} \left(|\nabla v|^{p(x,u_n)-2} \nabla v-|\nabla v|^{p(x,u)-2} \nabla v\right) \nabla\left(u_n-v\right) d x\right| \\
\leq & \left(\int_{\O} \left|| \nabla v|^{p(x,u_n)-2} \nabla v-|\nabla v|^{p(x,u)-2} \nabla v\right|^{\frac{p^{+}}{p^{+}-1}} d x\right)^{\frac{p^{+}-1}{p^{+}}} \\
& \times\left(\int_{\O}\left|\nabla\left(u_n-v\right)\right|^{p^{+}} d x\right)^{\frac{1}{p^{+}}} \\
\leq & c_4\left(\int_{\O} \left|| \nabla v|^{p(x,u_n)-2} \nabla
v-|\nabla v|^{p(x,u)-2} \nabla v\right|^{\frac{p^{+}}{p^{+}-1}} d
x\right)^{\frac{p^{+}-1}{p^{+}}} \longrightarrow 0, \text{ as }
n\rightarrow +\infty.
\end{aligned}
\end{equation}
In a similar way, we have
\begin{equation}\label{e32}
   \int_{\O} \mu (x) \left(|\nabla v|^{q(x,u_n)-2} \nabla v-|\nabla v|^{q(x,u)-2} \nabla v\right) \nabla\left(u_n-v\right) d x \longrightarrow 0, \text{ as } n \rightarrow +\infty.
\end{equation}
Therefore, from  \eqref{e31}-\eqref{e32}, we deduce that
\begin{equation}\label{e33}
   \begin{aligned}
      & \int_{\O} \left(|\nabla v|^{p(x,u_n)-2} \nabla v-|\nabla v|^{p(x,u)-2} \nabla v\right) \nabla\left(u_n-v\right) d x \\
      &  +\int_{\O} \mu (x) \left(|\nabla v|^{q(x,u_n)-2} \nabla v-|\nabla v|^{q(x,u)-2} \nabla v\right) \nabla\left(u_n-v\right) d x\\&
        \longrightarrow 0,\text{ as } n \rightarrow +\infty.
   \end{aligned}
\end{equation}
Then, by \eqref{e28},\eqref{e30} and \eqref{e33}, we get the proof of Claim 2.\vspace{2mm}\\
Consequently, the proof of Step 3 follows from  Claims 1 and
2.\vspace{2mm}

\textbf{\underline{Step 4.}} $u$ is a weak solution of \eqref{P1}.\\
From $(\mathrm{H}_2)$, one deduce that the functional
$$w\mapsto \int_{\O}^{} f(x,u)(u-w)dx,$$
is continuous on $(E_u,\|.\|_{E_u})$. Using the density of
$\mathcal{C}^\infty _0 (\O)$ in $E_u$ (see Proposition \ref{p1}),
the inequality \eqref{e35} can be extended to the whole space $E_u.$
For $s > 0$ and $w\in  E_u$, by choosing $v = u -sw $ as a test
function in \eqref{e35}, we obtain
$$\begin{aligned}
s \int_{\O} f(x, u) w d x \geq & s \int_{\O} |\nabla u-s \nabla w|^{p(x,u)-2}(\nabla u-s \nabla w) \nabla w d x \\
& +s \int_{\O} \mu (x)|\nabla u-s\nabla w|^{q(x,u)-2}(\nabla
u-s\nabla w)\nabla w d x .
\end{aligned}$$
Therefore,
$$\begin{aligned}
& \int_{\O} f(x, u) w d x -  \int_{\O} |\nabla u-s \nabla w|^{p(x,u)-2}(\nabla u-s \nabla w) \nabla w d x \\
& + \int_{\O} \mu (x)|\nabla u-s\nabla w|^{q(x,u)-2}(\nabla
u-s\nabla w)\nabla w d x\geq 0 .
\end{aligned}$$
Tends $s$ to $0^+$ in this last inequality, we get

$$\begin{aligned}
\int_{\O} f(x, u) w d x \geq &  \int_{\O} |\nabla u|^{p(x,u)-2}\nabla u \nabla w d x \\
 &+ \int_{\O} \mu (x)|\nabla u|^{q(x,u)-2}\nabla u\nabla w d x .
\end{aligned}$$
Obviously, the same inequality exists with $( -w)$ instead of $w$.
Consequently,
$$\int_{\O} |\nabla u|^{p(x,u)-2}\nabla u \nabla w d x
 + \int_{\O} \mu (x)|\nabla u|^{q(x,u)-2}\nabla u\nabla w d x = \int_{\O} f(x, u) w d x, \ \text{ for all } w\in E_u.$$
 This ends the proof of Theorem \ref{th1}.
\subsection{The multiplicity of solutions: Proof of Theorem \ref{th2}}

We will now define the Carathéodory functions $f_+ :\Omega \times
\mathbb{R} \rightarrow \mathbb{R}$ as follows

\begin{equation}\label{s1}
f_{+}(x, t):= \begin{cases} f\left(x, t^{+}\right) & \text {if } t \leq \eta_{+} \\
\ & \  \\
f\left(x, \eta_{+}\right) & \text {if } t>\eta_{+}\end{cases}
\end{equation}
and $f_{-}: \Omega \times \mathbb{R} \rightarrow \mathbb{R}$ defined
by
\begin{equation}\label{s2}
f_{-}(x, t):= \begin{cases} f\left(x, \eta_{-}\right) & \text {if } t<\eta_{-} \\
\ & \  \\
f\left(x, t^{-}\right) & \text {if } t \geq \eta_{-},\end{cases}
\end{equation}
where $\eta_{+}$and $\eta_{-}$are defined in $(\mathrm{H}_4)$.
\newline
Note that, from condition $(\mathrm{H}_2)'$,  it is easy to see that
$f_+$ satisfies conditions $(\mathrm{H}_2)$. Therefore, invoking
Theorem \ref{th1}, there exists at least one weak solution for the
following problem:

\begin{equation}\label{s3}
 \left\lbrace
\begin{array}{rll}-\operatorname{div}\left(|\nabla u|^{p(x,u)-2} \nabla u+ \mu (x)|\nabla u|^{q(x,u)-2} \nabla u\right)  & =f_+(x, u) &
\text { in } \Omega, \\
\ & \ &  \\
u  & =0 & \text { on } \partial \Omega.\end{array}\right.
\end{equation}
Let $u_0 \in E_{u_0} $ be a weak solution of \eqref{s3}, that is,
\begin{equation}\label{s4}
\begin{aligned}&\int_{\O} \left(|\nabla u_0|^{p(x,u_0)-2}\nabla u_0 \nabla w +\mu (x)|\nabla u_0|^{q(x,u_0)-2}\nabla u_0\nabla w \right)d x\\& = \int_{\O} f_+(x, u_0) w d x, \ \text{ for all } w\in E_{u_0}.\end{aligned}
\end{equation}
We act with $w = u_0^-=\max \left\{0,-u\right\}$ in equation
\eqref{s4}, we get
$$\begin{aligned}
&\int_{\O} \left(|\nabla u_0|^{p(x,u_0)-2}\nabla u_0 \nabla u_0^- +
\mu (x)|\nabla u_0|^{q(x,u_0)-2}\nabla u_0\nabla u_0^- \right)d x
\\&= \int_{\O} f_+(x, u_0) u_0^- d x.\end{aligned}
$$
Using the truncation \eqref{s1}, we can conclude that
$$0\leq \int_{\O} |\nabla u_0^-|^{p(x,u_{0})}+\mu (x)|\nabla u_0^-|^{q(x,u_{0})}dx \leq 0.$$
It follows, since $\mu (x)\geq 0\text{ for a.a. } x \in \O$ and using Proposition \ref{GG}, that $u_0\geq 0.$\\

Again, we act with $w = (u_0 - \eta ^+)^+=\max\left\{u_0 - \eta
^+,0\right\}$ in \eqref{s4} and applying $(\mathrm{H}_4)$, we
conclude that
\begin{equation}\label{s5}
\begin{aligned}&\small\int_{\O}\left( |\nabla u_0|^{p(x,u_0)-2}\nabla u_0 \nabla (u_0-\eta _+)^+ + \mu (x)|\nabla u_0|^{q(x,u_0)-2}\nabla u_0\nabla (u_0-\eta _+)^+\right) d x \\&=\small \int_{\O} f_+(x, u_0) (u_0-\eta _+)^+d x\\
&=\small\int_{\O} f(x, \eta_+) (u_0-\eta_+)^+d x\\
&\leq 0.
\end{aligned}
\end{equation}
Recall that
$$
\begin{aligned}
&|\nabla u_0|^{p(x,u_0)-2}\nabla u_0 \nabla (u_0-\eta _+)^++\mu
(x)|\nabla u_0|^{q(x,u_0)-2}\nabla u_0 \nabla (u_0-\eta _+)^+\\&=
\big(|\nabla u_0|^{p(x,u_0)-2}\nabla u_0-|\nabla \eta
_+|^{p(x,u_0)-2}\nabla\eta _+\big)  \nabla (u_0-\eta _+)^+\\&\ + \mu
(x)\big(|\nabla u_0|^{q(x,u_0)-2}\nabla u_0-|\nabla \eta
_+|^{q(x,u_0)-2}\nabla\eta _+\big)  \nabla (u_0-\eta _+)^+.
\end{aligned}
$$
Then, using Lemma  \ref{smn} and \eqref{s5}, we infer that

$$\begin{aligned}
    0 & \geq \int_{u_0 \geq \eta _+}^{}\big(|\nabla u_0|^{p(x,u_0)-2}\nabla u_0-|\nabla \eta _+|^{p(x,u_0)-2}\nabla\eta _+\big)  \nabla (u_0-\eta _+)^+ dx\\
    &  \geq \int_{u_0 \geq \eta _+}^{} |\nabla (u_0-\eta _+)|^{p(x,u_0) }+\mu (x) |\nabla (u_0-\eta _+)|^{q(x,u_0)} dx\\
     & = \int_{\O} |\nabla (u_0-\eta _+)^+|^{p(x,u_0) }+\mu (x) |\nabla (u_0-\eta _+)^+|^{q(x,u_0)}dx,
  \end{aligned}$$
which implies, again using Proposition \ref{GG}, that $(u_0 -\eta_+)^+ \equiv 0$, namely $u_0 \in [0, \eta_+]$. By using the truncation \eqref{s1}, we can conclude that $u_0$  is a nonnegative bounded weak solution for problem \eqref{P1}.\\
Similarly, we take $v_0$ a weak solution for the following problem
$$\left\lbrace
\begin{array}{rll}
-\operatorname{div}\left(|\nabla u|^{p(x,u)-2} \nabla u+ \mu (x)|\nabla u|^{q(x,u)-2} \nabla u\right)  & =f_-(x, u) & \text { in } \Omega, \\
\ & \ & \\
u  & =0 & \text { on } \partial \Omega,\end{array}\right.
$$
satisfying,

\begin{equation}\label{s7}
\begin{aligned}&\int_{\O}\left( |\nabla v_0|^{p(x,v_0)-2}\nabla u_0 \nabla w + \mu (x)|\nabla v_0|^{q(x,v_0)-2}\nabla v_0\nabla w\right) d x \\&= \int_{\O} f_-(x, v_0) w d x, \  w\in E_{v_0}.\end{aligned}
\end{equation}
 If we take $w =v_0^+$ in \eqref{s7} it follows that $v_0\leq 0$. In addition, if we choose $w = (v_0-\eta_-)^-$, we infer that $ v_0\in [\eta_-,0].$ This ends the proof of Theorem \ref{th2}.



%


%
%

\vspace{1cm}
\textbf{Declarations:}\\
\noindent{\bf  Data Availability:} Data sharing not applicable to this article as no datasets were generated or analysed during
the current study.

\noindent{\bf Conflict of interest}: On behalf of all authors, the corresponding author states that there is no Conflict of
interest.


 \bigskip
\noindent \textsc{\textsc{Ala Eddine Bahrouni}} \\
Mathematics Department, Faculty of Sciences, University of Monastir,
5019 Monastir, Tunisia\\
 ({\tt ala.bahrouni@fsm.rnu.tn})

\bigskip
\noindent \textsc{\textsc{Anouar Bahrouni}} \\
Mathematics Department, Faculty of Sciences, University of Monastir,
5019 Monastir, Tunisia\\
 ({\tt Anouar.Bahrouni@fsm.rnu.tn; bahrounianouar@yahoo.fr})

\bigskip
\noindent \textsc{\textsc{Hlel Missaoui}} \\
Mathematics Department, Faculty of Sciences, University of Monastir,
5019 Monastir, Tunisia\\
 ({\tt hlelmissaoui55@gmail.com})

\end{document}